%% file: paper1.tex
\documentclass[leqno]{amsart}
\usepackage{egastyle}
\usepackage{mymacros}
\usepackage[all,tips]{xy}
\usepackage{url}
\usepackage{bm}
\usepackage{verbatim}
\usepackage{epsfig}
\usepackage{graphicx}

\CompileMatrices

\input{paper2-export.aux}

\newwrite\exportaux
\gdef\exportauxname{paper1-export.aux}
\immediate\openout\exportaux\exportauxname

\makeatletter
\let\oldlabel\label
\def\label#1{\@bsphack\immediate\write\exportaux%
  {\string\newlabel{#1}{{I.\@currentlabel}{\thepage}}}\@esphack%
  \oldlabel{#1}}
\makeatother

\def\cA{{\mathcal A}}
\def\cB{{\mathcal B}}
\def\cC{{\mathcal C}}
\def\cD{{\mathcal D}}

\def\cG{{\mathcal G}}

\def\cM{{\mathcal M}}

\def\cT{{\mathcal T}}
\def\cU{{\mathcal U}}

\def\cX{{\mathcal X}}
\def\cZ{{\mathcal Z}}

\def\B{\bB}
\def\A{\bA}
\def\P{\bP}
\def\G{\bG}

\newcommand{\CC}{\mathbf{C}}

\newcommand{\PP}{\mathbf{P}}

\newcommand{\QQ}{\mathbf{Q}}
\newcommand{\RR}{\mathbf{R}}
\newcommand{\ZZ}{\mathbf{Z}}

\renewcommand{\OO}{{\ms O}}

\DeclareMathOperator\val{val}
\DeclareMathOperator\DV{DV}

\DeclareMathOperator\Jac{Jac}

\DeclareMathOperator\Int{Int}
\DeclareMathOperator\pr{pr}
\DeclareMathOperator\Isom{Isom}
\renewcommand\red{\operatorname{red}}
\newcommand\reg{{\operatorname{reg}}}

\newcommand\an{{\operatorname{an}}}
\newcommand\sep{{\operatorname{sep}}}

\renewcommand\tilde\widetilde

\newcommand{\p}{{}'}

\title[Lifting harmonic morphisms I]{Lifting harmonic morphisms I: 
  metrized complexes and Berkovich skeleta}

\author{Omid Amini} 
\email{oamini@math.ens.fr}
\address{CNRS-DMA, \'Ecole Normale Sup\'erieure, 45 Rue d'Ulm, Paris}

\author{Matthew Baker} 
\email{mbaker@math.gatech.edu}
\address{School of Mathematics, Georgia Institute of Technology, Atlanta GA 30332-0160, USA}

\author{Erwan Brugall\'e} 
\email{brugalle@math.jussieu.fr}
\address{Universit\'e Pierre et Marie Curie, Paris 6, 4 Place Jussieu, 75 005 Paris, France}

\author{Joseph Rabinoff}
\email{jrabinoff@math.gatech.edu}
\address{School of Mathematics, Georgia Institute of Technology, Atlanta GA 30332-0160, USA}

\begin{document}

\begin{abstract}  
  Let $K$ be an algebraically closed, complete non-Archimedean field.  The
  purpose of this paper is to carefully study the extent to which finite
  morphisms of algebraic $K$-curves are controlled by certain
  combinatorial objects, called \emph{skeleta}.  A skeleton is a metric
  graph embedded in the Berkovich analytification of $X$.  A skeleton has
  the natural structure of a metrized complex of curves.  We prove that a finite
  morphism of $K$-curves gives rise to a finite harmonic morphism of a
  suitable choice of skeleta.  We use this to give analytic proofs of
  stronger ``skeletonized'' versions of some foundational results of
  Liu-Lorenzini, Coleman, and Liu on simultaneous semistable reduction of
  curves.  We then consider the inverse problem of lifting finite harmonic
  morphisms of metrized complexes to morphisms of curves over $K$.  We
  prove that every tamely ramified finite harmonic morphism of
  $\Lambda$-metrized complexes of $k$-curves lifts to a finite morphism of
  $K$-curves.  If in addition the ramification points are marked, we
  obtain a complete classification of all such lifts along with their
  automorphisms.  This generalizes and provides new analytic proofs of
  earlier results of Sa\"idi and Wewers.  As an application, we discuss
  the relationship between harmonic morphisms of metric graphs and induced
  maps between component groups of N{\'e}ron models, providing a negative
  answer to a question of Ribet motivated by number theory.

  This article is the first in a series of two.  The second article
  contains several applications of our lifting results to questions about
  lifting morphisms of tropical curves.

\end{abstract}

\thanks{The authors thank Amaury Thuillier and Antoine Ducros for their help with some of the descent arguments in 
\parref{par:descent}.
We are grateful to Andrew Obus for a number of useful comments based on a careful reading of the first arXiv version of this manuscript. 
M.B. was partially supported by NSF grant DMS-1201473.
E.B. was partially supported by the ANR-09-BLAN-0039-01.}

\maketitle

{\small Throughout this paper, unless explicitly stated otherwise,
$K$ denotes a complete algebraically closed non-Archimedean field with nontrivial valuation 
$\val:K\to\R\cup\{\infty\}$.  Its valuation ring is denoted $R$, its
maximal ideal is $\fm_R$, and the residue field is $k = R/\fm_R$.  We
denote the value group of $K$ by $\Lambda = \val(K^\times)\subset\R$.}

\section{Introduction}

This article is the first in a series of two.  The second, entitled
\emph{Lifting harmonic morphisms II: tropical curves and metrized complexes}, 
will be cited as~\cite{abbr:lifting2}; references of the form
``Theorem~II.1.1'' will refer to Theorem~1.1 in~\cite{abbr:lifting2}. 

\paragraph
The purpose of this paper is to carefully study the extent to which finite
morphisms of algebraic $K$-curves are controlled by certain combinatorial
objects, called \emph{skeleta}.  Let $X$ be a smooth, proper, connected
$K$-curve.  Roughly speaking, a \emph{triangulation} $(X,V\cup D)$ of $X$ (with
respect to a finite set of punctures $D\subset X(K)$) is a finite set $V$
of points in the Berkovich analytification $X^\an$ of $X$ whose removal
partitions $X^\an$ into open balls and finitely many open annuli (with the
punctures belonging to distinct open balls).  Triangulations of $(X,D)$
are naturally in one-to-one correspondence with semistable models $\fX$ of
$(X,D)$: see Section~\ref{sec:simultaneous.ss.reduction}.  A triangulation
$(X,V\cup D)$ gives rise to a skeleton $\Sigma(X,V\cup D)$ of $X$.  The
skeleton of a triangulated punctured curve is the dual graph of the special fiber
$\fX_k$ of the corresponding semistable model, equipped with a canonical
metric, along with completed infinite rays in the directions of the punctures.
There are many skeleta in $X^\an$, although if $X\setminus D$ is
hyperbolic (i.e.,\ has negative Euler characteristic 
$\chi(X, D) = 2-2g(X)-\#D$), then there exists a
unique minimal skeleton.  A skeleton of $(X,D)$ is by definition a subset
$\Sigma\subset X^\an$ of the form $\Sigma(X,V\cup D)$ for some
triangulation of $X$. 

\paragraph[Skeletal simultaneous semistable reduction]
First we will prove that finite morphisms of curves induce morphisms of
skeleta: 

\begin{introlabel}
\begin{thm} \label{thm:intro1}
  Let $\phi:X'\to X$ be a finite morphism of smooth, proper, connected
  $K$-curves and let $D\subset X(K)$ be a finite set.  There exists a
  skeleton $\Sigma\subset X$ such that
  $\Sigma':=\phi\inv(\Sigma)$ is a skeleton of $X'$.  For any such
  $\Sigma$ the map $\phi:\Sigma'\to\Sigma$ is a 
  finite harmonic morphism of metric graphs. 
\end{thm}
\end{introlabel}

\smallskip
See Corollaries~\ref{cor:stable.hull.skel}
and~\ref{cor:morphism.to.harmonic}.  
Harmonic morphisms of graphs are defined in
Section~\ref{sec:definitions}.  
Due to the close relationship 
between semistable models and skeleta, Theorem~\ref{thm:intro1} should
be interpreted as a very strong simultaneous semistable reduction
theorem.  In fact, we will show how to formally derive
from Theorem~\ref{thm:intro1} the simultaneous semistable reduction theorems of
Liu--Lorenzini~\cite{liu_lorenzini:models}, 
Coleman~\cite{coleman:stable_maps}, and Liu~\cite{liu:stable_hull}.
Our version of these results hold over any non-Archimedean field $K_0$,
not assumed to be discretely valued.  As an example, the following weak form
of Liu's theorem is a consequence of Theorem~\ref{thm:intro1}:

\begin{cor*}
  Let $X,X'$ be smooth, proper, geometrically connected curves over a
  non-Archimedean field $K_0$ and let $\phi:X'\to X$ be a finite
  morphism.  Then there exists a finite, separable extension $K_1$ of
  $K_0$ and semistable models $\fX,\fX'$ of the curves $X_{K_1},X'_{K_1}$,
  respectively, such that $\phi_{K_1}$ extends to a \textbf{finite} morphism
  $\fX'\to\fX$.
\end{cor*}

\smallskip
We discuss simultaneous semistable reduction theorems in
Section~\ref{sec:simultaneous.ss.reduction}.  We wish to emphasize that
Theorem~\ref{thm:intro1} does not follow from any classical simultaneous semistable
reduction theorem, in that there exist finite morphisms $\phi:\fX'\to\fX$
of semistable models such that the inverse image of the corresponding
skeleton of $X$ is not a skeleton of $X'$.
(See, however, \cite{cornelissen:graph_li-yau} where a skeletal version of
Liu's theorem is derived from Liu's method of proof in the discretely valued case.)
Our proof of Theorem~\ref{thm:intro1} is entirely analytic, resting on an
analysis of morphisms between open annuli and open balls; in particular it
makes almost no reference to semistable reduction theory, and is
therefore quite different from Liu and Liu--Lorenzini's approach.

A skeletal simultaneous semistable reduction theorem is important
precisely when one wants to obtain a well-behaved morphism of graphs from
a finite morphism of curves.  This is crucial for obtaining the spectral lower
bound on the gonality in~\cite{cornelissen:graph_li-yau}, for instance. 

\paragraph[Lifting harmonic morphisms]
Our second goal is in a sense inverse to the first: we wish
to lift finite harmonic morphisms of metric graphs to finite
morphisms of curves.  More precisely, let $(X,D)$ be a punctured
$K$-curve as above, let $\Sigma$ be a skeleton, and let 
$\bar\phi:\Sigma'\to\Sigma$ be a finite harmonic morphism of metric graphs.
It is natural to ask whether there exists a curve $X'$ and a finite morphism 
$\phi:X'\to X$ such that $\phi\inv(\Sigma)$ is a skeleton of $X'$ and
$\phi\inv(\Sigma)\cong\Sigma'$ as metric graphs over $\Sigma$.
In general the answer is ``no'': there are subtle ``Hurwitz obstructions''
to such a lift, as described below.  One solution is to enrich 
$\Sigma,\Sigma'$ with the extra structure of 
\emph{metrized complexes of curves}.  A metrized complex of curves is
basically a metric graph $\Gamma$ and for every (finite) vertex 
$p\in\Gamma$, the data of a smooth, proper $k$-curve 
$C_p$, along with an identification of the edges adjacent to $p$ with 
distinct $k$-points of $C_p$.  There is a notion of a finite harmonic morphism of
metrized complexes of curves, which consists of a finite harmonic morphism
$\phi:\Gamma'\to\Gamma$ of underlying metric graphs, and for every vertex 
$p'\in\Gamma'$ with $p = \phi(p')$, a finite morphism 
$\phi_{p'}:C_{p'}\to C_p$, satisfying various compatibility conditions.

Now let $\Sigma = \Sigma(X,V\cup D)$ for a triangulated punctured curve 
$(X,V\cup D)$.  We will show that
$\Sigma$ is naturally a metrized complex of curves, and that finite
harmonic morphisms $\Sigma'\to\Sigma$ of metrized complexes do lift to
finite morphisms of curves, under a mild tameness hypothesis.
Moreover we have very precise control over the set of such lifts entirely
in terms of the morphism of metrized complexes.  The main lifting theorem,
stated somewhat imprecisely, is as follows: 

\begin{introlabel}
\begin{thm}[Lifting theorem] \label{thm:lifting.intro}
  Let $(X,V\cup D)$ be a triangulated punctured $K$-curve, let $\Sigma$ be its
  skeleton, and let $\Sigma'\to\Sigma$ be a tame covering of metrized
  complexes of curves.  Then there exists a curve $X'$ and a finite morphism
  $\phi:X'\to X$, branched only over $D$, such that $\phi\inv(\Sigma)$ is
  a skeleton of $(X',\phi\inv(D))$ and $\phi\inv(\Sigma)\cong\Sigma'$ as metrized
  complexes of curves over $\Sigma$.  There are finitely many such lifts
  $X'$ up to $X$-isomorphism.  There are explicit descriptions of the set
  of lifts and the automorphism group of each lift (as a cover of $X$) in
  terms of the morphism $\Sigma'\to\Sigma$.
\end{thm}
\end{introlabel}

\smallskip
See Theorems~\ref{thm:lifting2}, \ref{thm:lifting}, and~\ref{thm:lifting3}
for the precise statements.  Tame coverings are defined in
Definition~\ref{def:tame.covering.metrized.complexes}.  
The hardest part of the proof of Theorem~\ref{thm:lifting.intro} is a local
lifting result, Theorem~\ref{thm:star.lifts}, in which the target curve
$X$ is replaced by a neighborhood of a vertex $p\in\Sigma$.  
Theorem~\ref{thm:star.lifts} eventually reduces to classical results about
the tamely ramified \'etale fundamental group.

It is worth mentioning that given a metrized complex of curves $\cC$ with
edge lengths contained in $\Lambda$, it is not hard to construct a
triangulated punctured $K$-curve $(X,V\cup D)$  such that 
$\Sigma(X,V\cup D)\cong\cC$.  See Theorem~\ref{thm:graph.to.curve}.

\paragraph
Our lifting results extend several theorems in the literature, as well as
making them more precise.  Sa\"idi~\cite[Th\'eor\'eme~3.7]{Sai97} proves
a version of the lifting theorem for semistable formal curves (without
punctures) over a complete discrete valuation ring.  His methods also make
use of the tamely ramified \'etale fundamental group and analytic gluing
arguments.  Wewers~\cite{wewers:tame_covers} works more generally with
marked curves over a complete Noetherian local ring using
deformation-theoretic arguments, proving that every tamely ramified
admissible cover of a marked semistable curve over the residue field of a
complete Noetherian local ring lifts.  Wewers classifies the possible
lifts (fixing, as we do, a lift of the target) in terms of certain
non-canonical ``deformation data'' depending on compatible choices of
formal coordinates at the nodes on both the source and target curves.  

Our results have several advantages compared to these previous approaches.
For one, we are able to work over non-Noetherian rank-$1$ valuation rings,
and descend the lifting theorems \textit{a postiori} to an essentially
arbitrary non-Archimedean field: see~\parref{par:descent.lifting}.  The real
novelty, however, is the systematic use of metrized complexes of curves in
the formulation of the lifting theorems, in particular the
classification of the set of lifts $X'$ of $\Sigma'\to\Sigma$.  
For example, in our situation, Wewers' non-canonical deformation data are replaced by
certain \textit{canonical}
gluing data, which means that the automorphism group
of the morphism of metrized complexes acts naturally on the set of gluing
data (but not on the set of deformation data).  This is what allows us to 
classify lifts up to isomorphism as covers of the target curve, as well as
determine the automorphism group of such a lift.  Sa\"idi's lifting
theorem is more similar to ours, with generators of inertia groups over
nodal points playing the role of gluing data.  It is not clear to the authors
whether his methods can be extended to give a classification of lifts.

The question of lifting (and classification of
all possible liftings) in the wildly ramified case is more subtle and
cannot be guaranteed in general.  See Remark~\ref{rem:wild.lifting}.

\paragraph[Homomorphisms of component groups]
We will discuss in Section~\ref{sec:applications1} the relationship
between harmonic morphisms of metric graphs and component groups of
N{\'e}ron models, connecting our lifting results to natural questions in
arithmetic algebraic geometry.  The idea is as follows.  Let $X$ be a smooth,
proper, geometrically connected curve over a discretely valued field
$K_0$, and suppose that $X$ admits a semistable model $\fX$ over the
valuation ring $R_0$ of $K_0$.  Let $J$ be the Jacobian of $X$ and let
$\Phi_X$ be the group of  connected components of the special fiber of the
N\'eron model of $J$.  

Given a metric graph $\Gamma$ with integer edge lengths, there is
a naturally associated abelian group $\Jac_\reg(\Gamma)$, called its
\emph{regularized Jacobian}.  The association 
$\Gamma\mapsto\Jac_\reg(\Gamma)$ satisfies Picard and Albanese
funtoriality: that is, given a finite harmonic morphism
$\Gamma'\to\Gamma$ there are natural push-forward and pull-back
homomorphisms
$\Jac_\reg(\Gamma')\to\Jac_\reg(\Gamma)$ and
$\Jac_\reg(\Gamma)\to\Jac_\reg(\Gamma')$.

A series of theorems of Raynaud, translated into our language, imply that
$\Phi_X\cong\Jac_\reg(\Sigma_\fX)$ canonically, where $\Sigma_\fX$ is the
skeleton induced by the semistable model $\fX$.  We observe that this
isomorphism respects both Picard and Albanese functoriality.

For modular curves such as $X_0(N)$ with $N$ prime, the minimal regular model has a special fiber consisting of two projective lines intersecting transversely.  Maps from such curves to elliptic curves, and their induced maps on component groups of N{\'e}ron models, play an important role in arithmetic geometry; for example, Ribet's work establishing that the Shimura-Taniyama conjecture implies Fermat's Last Theorem was motivated by the failure of such induced maps on component groups to be surjective in general.  (The cokernel of such maps controls congruences between modular forms.)  
Ribet noticed that if one considers instead the case where the target curve has
genus at least 2, then in all examples he knew of, the induced map on component groups was surjective, and 
in personal communication with the second author, he asked whether this was a general property of curves whose reduction looks like that of $X_0(N)$.  
As a concrete application of Theorem~\ref{thm:lifting.intro}, the next proposition provides a negative answer to this question.

\begin{introlabel}
\begin{prop} \label{prop:intro2}
  There exists  a finite morphism $f : X' \to X$
  of semistable curves over
  a discretely valued field $K_0$ with $g(X)\geq 2$ such that:
  \begin{itemize}
  \item the special fiber of the minimal regular model of $X'$ consists of
    two projective lines intersecting transversely;
  \item  the induced map $f_* : \Phi_{X'} \to \Phi_X$ on component groups
    of N{\'e}ron models of Jacobians is not  surjective.
  \end{itemize}
\end{prop}
\end{introlabel}

\smallskip
In order to prove Proposition~\ref{prop:intro2}, we proceed as follows.  
As component groups can be calculated in terms of metric graphs, one first
finds a finite harmonic morphism of metric graphs $\phi:\Gamma'\to\Gamma$
such that $\Gamma'$ has two vertices and at
least three edges, all of which have length one, and such that
$\phi_*:\Jac_\reg(\Gamma')\to\Jac_\reg(\Gamma)$ is not surjective.  One
then enriches $\phi$ to a morphism of metrized complexes of curves, with
rational residue curves, and
uses Theorem~\ref{thm:graph.to.curve} and Theorem~\ref{thm:lifting.intro} (along
with a descent argument) to
lift this to a finite morphism of curves 
$X'\to X$.  Such a morphism satisfies the conditions of
Proposition~\ref{prop:intro2}.

\paragraph[Applications to tropical lifting theorems]
Tropical geometry has many applications to algebraic
geometry, in particular enumerative algebraic geometry.  The basic
strategy is to associate a combinatorial ``tropical'' object to an
algebraic object, and prove that counting the former is somehow equivalent
to counting the latter.  An important ingredient in any such argument is a
tropical lifting theorem, which is a precise description of the algebraic
lifts of a given tropical object.  

Many of the intended applications of our lifting results are contained in
the second paper~\cite{abbr:lifting2}, in which we prove a number of
tropical lifting theorems.  The skeleton of a curve, viewed as a metric
graph, is a tropical object associated to the curve; our lifting results
can therefore be interpreted as tropical lifting theorems for morphisms of
metric graphs, in the following sense.  Recall that $\Lambda\subset\R$ is
the value group of $K$.  A \emph{$\Lambda$-metric graph} is a metric
graph $\Gamma$ whose edge lengths are contained in $\Lambda$.
We say that a finite harmonic morphism 
$\bar\phi:\Gamma'\to\Gamma$ of $\Lambda$-metric graphs is \emph{liftable}
provided that there exists a finite morphism of $K$-curves $\phi:X'\to X$,
a set of punctures $D\subset X(K)$, and a skeleton $\Sigma$ of $(X,D)$
such that $\Sigma\cong\Gamma$ as metric graphs, $\phi\inv(\Sigma)$ is a
skeleton of $(X',\phi\inv(D))$, and $\Sigma'\cong\Gamma'$ as metric graphs
over $\Sigma\cong\Gamma$.  Theorem~\ref{thm:lifting.intro} (along with
Theorem~\ref{thm:graph.to.curve}) shows that the
obstruction to lifting a morphism of metric graphs
$\bar\phi:\Gamma'\to\Gamma$ to a morphism of curves is the same as the
obstruction to enriching $\phi$ with the structure of a morphism of
metrized complexes of curves.  The following more precise statement is an
immediate consequence of Proposition~\ref{prop:lifting} and
Theorem~\ref{thm:graph.to.curve}: 

\begin{introlabel} 
\begin{cor} \label{cor:intro.lifting}
  Assume $\chr(k) = 0$.  Let $\bar\phi:\Gamma'\to\Gamma$ be a finite harmonic
  morphism of $\Lambda$-metric graphs.
  Suppose that $\phi$ can be enriched to a morphism of metrized complexes
  of curves.  Then $\phi$ is liftable.
\end{cor}
\end{introlabel}

\smallskip

We remark that we actually take advantage of the stronger form of
Theorem~\ref{thm:lifting.intro}, in particular the calculation of the
automorphism group of a lift, when lifting group actions on metrized
complexes in~\cite{abbr:lifting2}.

The problem of enriching a morphism of metric graphs to a morphism of
metrized complexes is essentially equivalent to the problem of assigning $k$-curves
$C_p,C_{p'}$ to vertices $p\in\Gamma,p'\in\Gamma'$ and
finding morphisms of $k$-curves $C_{p'}\to C_p$ with prescribed
ramification profiles.  The latter is a question about whether certain
Hurwitz numbers are nonzero. 
See Proposition~\ref{prop:lift augmented to complex}. 
Allowing arbitrary choices of residue curves $C_p,C_{p'}$, this can always
be done:

\begin{thm*}
  Suppose that $\chr(k) = 0$.
  Any finite harmonic morphism $\phi:\Gamma'\to \Gamma$ 
  of $\Lambda$-metric graphs is liftable to a morphism of curves.
\end{thm*}

\smallskip
See Theorem~\ref{thm:lifting harm}.
Imposing requirements on the lifted curves $X,X'$ give rise to variants of this
tropical lifting problem.  The most basic such requirement is to specify
the genera of the curves.  An \emph{augmented metric graph} is a metric
graph $\Gamma$ along with an integer $g(p)\in\Z_{\geq 0}$ for every finite
vertex $p$ of $\Gamma$, called the \emph{genus}.
The \emph{genus} of an augmented graph $\Gamma$ is defined to be
\[ g(\Gamma) = h_1(\Gamma) + \sum_p g(p), \]
where $h_1(\Gamma)$ is the first Betti number.  
Any metrized complex of
curves has an underlying augmented metric graph, where $g(p) := g(C_p)$,
the genus of the residue curve.
If $\Sigma$ is a skeleton of a curve $X$, viewed
as an augmented metric graph, then we have $g(\Sigma) = g(X)$.

A finite harmonic morphism of augmented metric graphs $\Gamma'\to\Gamma$
is \emph{liftable} if the morphism of underlying metric graphs is liftable
to a finite morphism of curves $X'\to X$ such that the isomorphism of
$\Gamma$ (resp.\ $\Gamma'$) with a skeleton of $X$ (resp.\ $X'$) respects
the augmented metric graph structure.  In this case, $g(X)=g(\Gamma)$ and
$g(X')=g(\Gamma')$.  By our lifting theorem, $\Gamma'\to\Gamma$ is
liftable to a morphism of curves if and only if it can be enhanced to a
morphism of metrized complexes such that $g(C_p) = g(p)$ for every finite
vertex.  Now the obstruction to lifting is nontrivial: there exist finite
harmonic morphisms of metric graphs which cannot be enhanced to a morphism
of metrized complexes.  This corresponds to the emptiness of certain
Hurwitz spaces over $k$.  See Section~\ref{sec:examples} for the following
striking example of this phenomenon:

\begin{thm*}
  There exists an augmented metric graph $\Gamma$ admitting a finite
  morphism of degree $d > 0$ to an augmented metric graph $T$ of genus
  zero, but which is not isomorphic to a skeleton of a $d$-gonal curve
  $X$. 
\end{thm*}

\smallskip 

In~\cite{abbr:lifting2} we also give applications of the lifting theorem
to the following kinds of tropical lifting problems, among others:
\begin{itemize}
\item When $T$ is an augmented metric graph of genus zero, we study a
  variant of the lifting problem in which the genus of the source curve is
  prescribed, but the degree of the morphism is not.   
  We use this to prove that linear equivalence of divisors on a tropical
  curve $C$  coincides with the equivalence relation generated by
  declaring that the fibers of every finite harmonic morphism from
  $C$ to the tropical projective line are equivalent.  

\item We study liftability of metrized complexes equipped with a finite
  group action, and as an application classify all (augmented) metric
  graphs arising as the skeleton of a hyperelliptic curve.
  As mentioned above, this study actually takes advantage of the more
  precise form of Theorem~\ref{thm:lifting.intro}, not just the existence
  of a lift.
\end{itemize}

\paragraph[Organization of the paper]
In Section~\ref{sec:definitions} we recall some basic definitions of
graphs with additional structures, including metrized complexes of curves,
and define harmonic morphisms between them.  We provide somewhat more
detail than is strictly necessary for the exposition in this paper, but
will prove necessary in~\cite{abbr:lifting2}.  

In Section~\ref{sec:skeleta} we review the structure theory of analytic
curves, and define the skeleton of a curve.  We show that the skeleton is
naturally a metrized complex of curves, and we prove that any metrized
complex arises in this way.  In Section~\ref{sec:morphisms.skeleta} we
develop ``relative'' versions of the results of
Section~\ref{sec:skeleta}.  We prove that a finite morphism of curves
induces a finite morphism between suitable choices of skeleta, and that
the map on skeleta is a finite harmonic morphism of metrized complexes of
curves.  In Section~\ref{sec:simultaneous.ss.reduction} we show how
the previous two sections can be used to derive various simultaneous
semistable reduction theorems.

In Section~\ref{sec:local lifting} we prove a local lifting theorem, which
essentially says that a tame covering of metrized complexes of
curves can be lifted in a unique way to a tame covering of analytic
$K$-curves in a neighborhood of a vertex.  This is the technical heart of
the lifting theorem.  In Section~\ref{sec:global.lifting} we globalize the
considerations of the previous section, giving the classification of lifts
of a tame covering of metrized complexes of curves to a tame covering of
triangulated punctured $K$-curves.

Section~\ref{sec:applications1} contains an example application, in which
we construct a morphism of curves over a discretely valued field such that
the covariantly associated morphism of component groups of Jacobians has a
prescribed behavior.  We use this to answer a question of Ribet.

\paragraph[Related work] 
Another framework for the foundations of tropical geometry has been proposed by 
Kontsevich-Soibelman~\cite{kontsevich_soibelman:SYZ, kontsevich_soibelman:affinestructures} 
and Mikhalkin \cite{mikhalkin:tropical_geometry,Mikhalkin:enumerative},
in which tropical objects are associated to real one-parameter families of complex varieties. 
We refer to~\cite{kontsevich_soibelman:SYZ, kontsevich_soibelman:affinestructures} 
for some conjectures relating this framework to Berkovich spaces.
In this setting,  the notion of metrized complex of curves is similar  to the notion of
phase-tropical curves, and Proposition \ref{prop:lifting} is a
consequence of Riemann's Existence Theorem. We refer the interested
reader to the forthcoming paper \cite{Mik08} for more details 
(see also \cite{Br13} where this
statement is implicit, as well as Corollary
\ref{cor:lift iff H not 0}).  

\smallskip

Harmonic morphisms of finite graphs were introduced by Urakawa in \cite{urakawa} and further explored in \cite{BN09}.
Harmonic morphisms of metric graphs 
have been introduced independently by several different people.
Except for the integrality condition on the slopes, they appear already 
in Anand's paper~\cite{anand}.
The definition we use is the same as the one given
in~\cite{mikhalkin:tropical_geometry, Cha12}.

\smallskip

Finite harmonic morphisms of metrized complexes can be regarded as a metrized version of
the notion of admissible cover due to Harris and Mumford~\cite{HM82} (where, in addition,
arbitrary ramification above smooth points is allowed).
Recall that for two semistable curves $Y'$ and $Y$ over a field $k$,
a finite surjective morphism $\phi:Y' \to Y$
is an {\em admissible cover} if $\phi\inv(Y^{\rm sing})=Y'^{\rm sing}$
and for each singular point $y'$ of $Y'$, the ramification indices at $y'$ along the two branches intersecting at $y'$ coincide.
(In addition, one usually imposes that all the other ramifications of $\phi$ are simple).
An admissible cover naturally gives rise to a finite harmonic morphism of metrized complexes:
denoting by $\cC$ the regularization of $Y$ (the metrized complex associated to $Y$
in which each edge has length one), define $\cC'$ as the metrized complex
obtained from $Y'$ by letting the length of the edge associated to
the double point $y'$ be $1/r_{y'}$ (where $r_{y'}$ is the ramification
index of $\phi$ at $y'$ along either of the two branches).
The morphism $\phi$ of semistable curves naturally extends to a finite harmonic morphism $\phi:\cC'
\to \cC$: on each edge $e'$ of $\cC'$ corresponding to a double point $y'$ of
$Y'$, the restriction of $\phi$ to $e'$ is linear with slope (or ``expansion factor'') $r_{y'}$.
Conversely, a finite harmonic morphism of metrized complexes of curves gives rise to an admissible cover of semistable curves
(without the supplementary condition on simple ramification) by forgetting the metrics on both sides, remembering only the expansion factor along each edge. 

\smallskip

A related but somewhat different Berkovich-theoretic point of view on
simultaneous semistable reduction for curves can be found in
Welliaveetil's recent preprint \cite{welliaveetil:RH_for_skeleta};
harmonic morphisms of metrized complexes of curves play an implicit role
in his paper.

\section{Metric graphs and metrized complexes of curves}
\label{sec:definitions}

In this section we recall several definitions of graphs with some
additional structures and morphisms between them.  A number of these
definitions are now standard in tropical geometry; we refer for example to
\cite{BakerNorine}, \cite{mikhalkin:tropical_geometry}, \cite{Br13}, and
\cite{AB12}.  We also provide a number of examples.
Some of the definitions in this section are not necessary in the 
generality or the form in which they are presented for the purposes of
this paper, but will be useful in~\cite{abbr:lifting2}.

Throughout this section, we fix $\Lambda$ a non-trivial subgroup of
$(\RR,+)$.

\paragraph[Metric graphs]
Given $r\in\ZZ_{\ge 1}$, we define $S_r\subset\CC$ to be a ``star with $r$ branches'', i.e., 
a topological space homeomorphic to the convex hull in $\RR^2$ of $(0,0)$ and any $r$ points, no two of which lie
on a common line through the origin.
We also define $S_0=\{0\}$.

A \emph{finite topological graph $\Gamma$} is the topological realization of a
finite graph. That is to say, $\Gamma$ is a compact 1-dimensional
topological space
such that for any point $p\in\Gamma$, there exists a neighborhood $U_p$
of $p$ in $\Gamma$ homeomorphic to some $S_r$; moreover there are only
finitely many
points $p$ with $U_p$ homeomorphic to $S_r$ with $r\ne 2$.

The unique integer $r$ such that $U_p$ is homeomorphic to $S_r$ is called
the \emph{valence} of $p$ and denoted $\val(p)$. 
A point of valence different from $2$ is called an \emph{essential vertex} of
$\Gamma$.  The set of \emph{tangent directions} at $p$ is 
$T_p(\Gamma) = \varinjlim_{U_p} \pi_0(U_p\setminus\{p\})$, where the limit
is taken over all neighborhoods of $p$ isomorphic to a star with $r$
branches.  
The set $T_p(\Gamma)$ has $\val(p)$ elements.

\begin{defn}
A \emph{metric graph} is a finite connected topological graph $\Gamma$
equipped with a complete inner metric on $\Gamma \setminus V_\infty(\Gamma)$,
where $V_\infty(\Gamma)\subsetneq\Gamma$ is some (finite) set of $1$-valent
vertices of $\Gamma$ 
called \emph{infinite vertices} of $\Gamma$.
(An \emph{inner metric} is a metric for which the distance between two points $x$ and $y$ is the minimum of the lengths of all paths between $x$ and $y$.)

Let $\Gamma$ be a metric graph containing a finite essential vertex.  We
say that $\Gamma$ is a \emph{$\Lambda$-metric graph} if the
distance between any two finite essential vertices of $\Gamma$ lies in
$\Lambda$.  A \emph{$\Lambda$-point}
of $\Gamma$ is a point of $\Gamma$ whose distance to any finite
essential vertex of $\Gamma$ lies in $\Lambda$. 

Let $\Gamma$ be a metric graph such that all essential vertices 
are infinite, so $\Gamma$ is homeomorphic to a circle or
a completed line.  We say that $\Gamma$ is a \emph{$\Lambda$-metric graph} 
if it is a completed line or a circle whose circumference is in $\Gamma$.  
Choose a point $v\in\Gamma$ of valency equal to $2$. 
A \emph{$\Lambda$-point} of $\Gamma$ is a point of $\Gamma$ whose distance
to $v$ lies in $\Lambda$. 
\end{defn}

\smallskip

One can equip $\Gamma$ with a degenerate metric in which the infinite vertices are at infinite distance from any other point of $\Gamma$.
When no confusion is possible about the subgroup $\Lambda$, we will sometimes write simply \emph{metric graph} instead of 
\emph{$\Lambda$-metric graph}.

\begin{defn}
Let $\Gamma$ be a $\Lambda$-metric graph. A \emph{vertex set}
$V(\Gamma)$ is a finite subset of the $\Lambda$-points
of $\Gamma$ containing all essential vertices.
An element of a fixed vertex set $V(\Gamma)$ is called a \emph{vertex} of
$\Gamma$, and the closure of a connected component of
$\Gamma\setminus V(\Gamma)$ is called an \emph{edge} of $\Gamma$. 
An edge which is homeomorphic to a circle is called a \emph{loop edge}.
An edge adjacent to an infinite vertex is called an \emph{infinite edge}.
We denote by $V_f(\Gamma)$ the set of finite vertices of $\Gamma$, and
by
$E_f(\Gamma)$ the set of finite edges of $\Gamma$. 
\end{defn}

\smallskip
Fix a vertex set $V(\Gamma)$.
We denote by $E(\Gamma)$ the set of edges of $\Gamma$.
Since $\Gamma$ is a metric graph, we can associate to each edge $e$ of
$\Gamma$ its length $\ell(e)\in\Lambda\cup\{+\infty\}$. Since the
metric on $\Gamma\setminus V_\infty(\Gamma)$ is complete, an edge $e$
is infinite if and only if $\ell(e)=+\infty$.
The notion of vertices and edges of $\Gamma$ depends, of course, on the
choice of a vertex set; we will fix such a choice without comment whenever
there is no danger of confusion.

\begin{defn}\label{def:morph metric}
Fix vertex sets $V(\Gamma')$ and $V(\Gamma)$ for two $\Lambda$-metric graphs $\Gamma'$ and $\Gamma$, respectively, and
let $\phi:\Gamma'\to\Gamma$ be a  continuous map.
\begin{itemize}
\item The map $\phi$ is called a \emph{$(V(\Gamma'),V(\Gamma))$-morphism of $\Lambda$-metric
  graphs} if we have $\phi(V(\Gamma'))\subset V(\Gamma)$, $\phi^{-1}(E(\Gamma))\subset E(\Gamma')$, and
the restriction
of $\phi$ to any edge $e'$ of $\Gamma'$ is a dilation by some factor
$d_{e'}(\phi)\in\Z_{\geq 0}$.

\item The map $\phi$ is called a \emph{morphism of $\Lambda$-metric graphs} if 
there exists a vertex set $V(\Gamma')$ of $\Gamma'$ and a 
vertex set $V(\Gamma)$ of $\Gamma$ such that
$\phi$ is a $(V(\Gamma'),V(\Gamma))$-morphism of $\Lambda$-metric
  graphs.

\item The map $\phi$ is said to be \emph{finite} if 
  $d_{e'}(\phi)>0$ for any edge $e'\in E(\Gamma')$.
\end{itemize}
\end{defn}

\smallskip
An edge $e'$ of $\Gamma'$ is mapped to a vertex of $\Gamma$ if and only if
$d_{e'}(\phi)=0$. Such an edge is said to be \emph{contracted} by $\phi$.
A morphism $\phi:\Gamma'\to\Gamma$ is finite if and only if there are no
contracted edges, which holds if and only if $\phi^{-1}(p)$ is a finite
set for any point $p\in\Gamma$.  If $\phi$ is finite, then 
$p'\in V_f(\Gamma')$ if and only if $\phi(p')\in V_f(\Gamma)$.

The integer $d_{e'}(\phi)\in\Z_{\geq 0}$ in Definition \ref{def:morph
  metric}
is called the \emph{degree} of $\phi$ along $e'$ (it is also sometimes 
called the \emph{weight} of $e'$ or  \emph{expansion factor} along $e'$ in
the literature).
Since $\ell(\phi(e')) = d_{e'}(\phi)\cdot\ell(e')$, it follows in particular 
that if $d_{e'}(\phi) \geq 1$ then $e'$ is infinite if and only if $\phi(e')$ is infinite.
Let $p'\in V(\Gamma')$, let $v'\in T_{p'}(\Gamma')$, and
let $e'\in E(\Gamma')$ be the edge in the direction of $v'$.  We define
the \emph{directional derivative of $\phi$ in the direction $v'$} to be
$d_{v'}(\phi) \coloneq d_{e'}(\phi)$.
If we set $p = \phi(p')$, then $\phi$ induces a map
\[ d\phi(p')~:~\big\{ v'\in T_{p'}(\Gamma')~:~d_{v'}(\phi)\neq 0 \big\}
\To T_p(\Gamma) \]
in the obvious way.

\begin{eg}
Figure \ref{fig:ex morph metric} depicts four examples of morphisms of metric graphs
$\phi:\Gamma'\to\Gamma$. We use the following
conventions in our pictures: black dots represent vertices of $\Gamma'$
and $\Gamma$; we label an edge with its degree if and only if the degree is different from $0$ and $1$;
we do not specify the lengths of edges of $\Gamma'$ and $\Gamma$.

The morphisms depicted in Figure \ref{fig:ex morph metric}(a), (b), and (c)
are finite, while the one depicted in Figure \ref{fig:ex morph metric}(d) is not.
\begin{figure}[h]
\begin{tabular}{ccccccc}
\scalebox{.31}{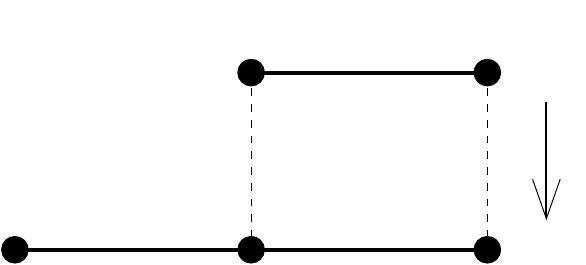}&  \hspace{3ex}&\scalebox{.31}{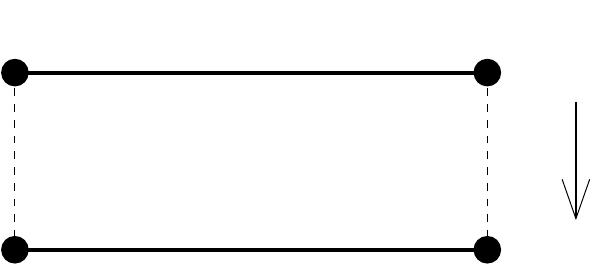}&  \hspace{3ex}&
\scalebox{.31}{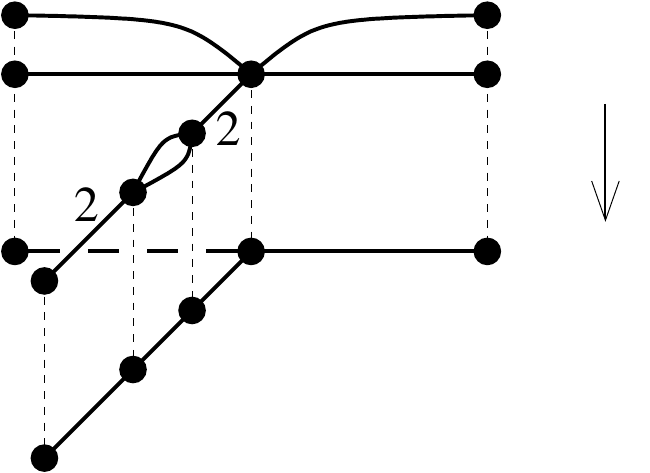}&  \hspace{3ex}&
\scalebox{.31}{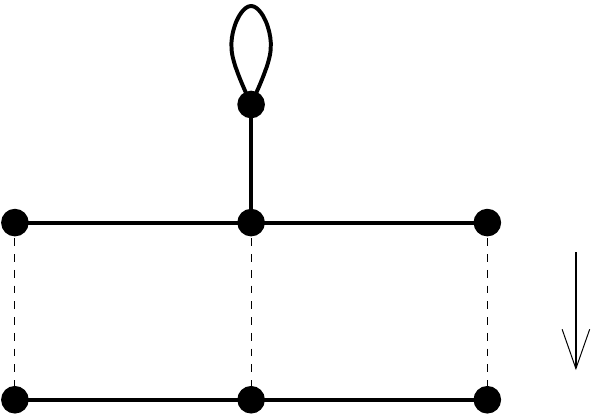}
\\ \\a)   && b) 
&& c) && d)
\end{tabular}
\caption{}\label{fig:ex morph metric}
\end{figure}
\end{eg}

\begin{defn}
Let $\phi:\Gamma'\to\Gamma$ be a morphism of $\Lambda$-metric graphs,
let $p'\in\Gamma'$, and let $p = \phi(p')$.  The morphism $\phi$
is \emph{harmonic at $p'$} provided that, for every tangent
direction $v\in T_{p}(\Gamma)$, the number 
\[ d_{p'}(\phi):=\sum_{\substack{v'\in T_{p'}(\Gamma')\\v'\mapsto v}} 
d_{v'}(\phi) \]
is independent of $v$.
The number $d_{p'}(\phi)$  is called the \emph{degree of $\phi$ at $p'$}.  

We say that $\phi$ is
\emph{harmonic} if it is \emph{surjective} and 
harmonic at all $p'\in \Gamma'$; in this case
the number $\deg(\phi) = \sum_{p'\mapsto p} d_{p'}(\phi)$ is
independent of $p\in\Gamma$, and is called the \emph{degree} of
$\phi$.  
 \end{defn}

\begin{rem}
  In the above situation, if $\Gamma$ consists of a single vertex $v$ 
  and $\phi: \Gamma'\to\Gamma$ is a  morphism of metric graphs, then
  $\phi$ is by definition harmonic.  
  In this case 
  the quantity $d_{p'}(\phi)$ is not defined, so we
  include the choice of a positive integer $d_{p'}(\phi)$ for each
  $p'\in V_f(\Gamma')$
  as part of the data of the morphism $\phi$. 
  Note that if $\phi$ is finite and $\Gamma = \{p\}$, 
  then $\Gamma' = \{p' \}$.

  If both $\Gamma'$ and $\Gamma$ have at least one edge, then a 
  non-constant morphism which is harmonic at all $p'\in V(\Gamma')$ is
  automatically surjective.
\end{rem}

\begin{eg} \label{eg:harmonic.examples}
The morphism depicted in Figure \ref{fig:ex morph metric}(a) is not
harmonic,
while the ones  depicted in Figure \ref{fig:ex morph metric}(b), (c), and (d)
are (for suitable choices of lengths).
\end{eg}

\paragraph[Harmonic morphisms and harmonic functions]
\label{par:harmonic.germs}
Given a metric graph $\Gamma$ and a non-empty open set $U$ in $\Gamma$,
a function $f:U \rightarrow \RR$ is said to be \emph{harmonic on $U$} if
$f$ is piecewise affine with integer slopes and
for all $x \in U$, the sum of the slopes of $f$ along all outgoing tangent directions at $x$ is equal to $0$.

Since we will not use it elsewhere in the paper, we omit the proof of the following result 
(which is very similar to the proof of \cite[Proposition 2.6]{BN09}):

\begin{prop}
\label{prop:harmonic.germs}
A morphism $\phi: \Gamma' \to \Gamma$ of metric graphs is harmonic if and only if for every open set $U \subseteq \Gamma$
and every harmonic function $f: U \to \RR$, the pullback function $\phi^* f : \phi^{-1}(U) \to \RR$ is also harmonic.
\end{prop}
\smallskip

Equivalently, a morphism of metric graphs is harmonic if and only if germs of harmonic functions pull back to germs of harmonic functions.

\paragraph
\label{par:pushforward.pullback}
For a metric graph $\Gamma$, we let $\Div(\Gamma)$ denote the free abelian group on $\Gamma$.
Given a harmonic morphism $\phi: \Gamma'\to\Gamma$ of metric graphs, 
there are natural pull-back and push-forward homomorphisms
$\phi^*:\Div(\Gamma)\to\Div(\Gamma')$ and
$\phi_*:\Div(\Gamma')\to\Div(\Gamma)$ 
 defined by
\[ \phi^*(p) = \sum_{p'\mapsto p} d_{p'}(\phi)\,(p') \quad
\mbox{and}\quad  \phi_*(p') = (\phi(p')) \]
and extending linearly.  It is clear that for $D\in\Div(\Gamma)$ we have
$\deg(\phi^*(D)) = \deg(\phi)\cdot\deg(D)$ and
$\deg(\phi_*(D)) = \deg(D)$.

\paragraph[Augmented metric graphs]
Here we consider $\Lambda$-metric
graphs together with the data of a non-negative integer at each point. For
a finite point $p$, this integer should be thought as the genus
of a ``virtual algebraic curve'' lying above $p$.

\begin{defn}
An \emph{augmented $\Lambda$-metric graph} is a $\Lambda$-metric
graph $\Gamma$ along with a function
$g:\Gamma\to\Z_{\geq 0}$, called the genus function,
such that $g(p)=0$ for all points of $\Gamma$ except possibly for finitely
many (finite) $\Lambda$-points $p\in\Gamma$.
The \emph{essential vertices} of $\Gamma$ are
the points $p$ for which $\val(p)\ne 2$ or
$g(p)>0$.
A \emph{vertex set} of $\Gamma$ is a vertex set  $V(\Gamma)$ of the underlying
metric graph which contains all essential vertices of $\Gamma$ as an augmented graph.

An augmented metric graph is said to be \emph{totally degenerate} if the genus
function is identically zero. 

The \emph{genus} of
$\Gamma$ is defined to be
\[ g(\Gamma) = h_1(\Gamma) + \sum_{p\in \Gamma} g(p), \]
where $h_1(\Gamma) = \dim_\Q H_1(\Gamma,\Q)$ is the first Betti number of $\Gamma$.
The \emph{canonical divisor} on $\Gamma$ is 
\begin{equation} \label{eq:canonical.divisor.on.augmented.graphs} 
K_\Gamma = \sum_{p\in \Gamma} (\val(p) + 2g(p) - 2)\, (p). 
\end{equation}

A \emph{harmonic morphism of  augmented $\Lambda$-metric
graphs} $\phi:\Gamma' \to\Gamma$ is a map which is
 a harmonic morphism between the
underlying metric graphs of $\Gamma'$ and $\Gamma$.
\end{defn}

\smallskip
Note that both summations in the above definition are in fact over essential vertices of $\Gamma$.
The degree of the canonical divisor of an
augmented $\Lambda$-metric graph $\Gamma$
 is $\deg(K_\Gamma) = 2 g(\Gamma) - 2$.
An augmented metric graph of genus 0 will also be called a
\emph{rational} augmented metric graph.

\paragraph
Let  $\phi:\Gamma'\to\Gamma$ be a harmonic morphism of
augmented $\Lambda$-metric 
graphs.
The \emph{ramification divisor} of $\phi$ is the divisor 
$R = \sum R_{p'} (p')$, 
where
\begin{equation} \label{eq:Rvprime}
  R_{p'} = d_{p'}(\phi)\cdot\big(2-2g(\phi(p'))\big) - \big(2-2g(p')\big)
  - \sum_{v' \in T_{p'}(\Gamma')} \big(d_{v'}(\phi)-1\big). 
\end{equation}
  
Note that if $p'\in V_\infty(\Gamma')$, then $R_{p'}\geq0$ if $d_{p'}(\phi)>0$, and 
$R_{p'}=-1$ if $d_{p'}(\phi)=0$.
The definition of $R_{v'}$ in \eqref{eq:Rvprime} is rigged so that we have the following graph-theoretic analogue of the \emph{Riemann--Hurwitz formula}:
\begin{equation} \label{eq:graph.RH}
  K_{\Gamma'} = \phi^*(K_\Gamma) + R.
\end{equation}
In particular, the sum $R = \sum R_{p'} (p')$ is in fact finite.

\smallskip

\begin{defn}
Let $\phi:\Gamma'\to\Gamma$ be a finite harmonic morphism of 
augmented $\Lambda$-metric
graphs.
We say that $\phi$ is \emph{\'etale} at a point $p'\in \Gamma'$
provided that $R_{p'}=0$.
We say that $\phi$ is
\emph{generically \'etale} if $R$ is supported on the set of infinite
vertices of $\Gamma'$ and that $\phi$ is \emph{\'etale} if $R=0$.
\end{defn}

\paragraph[Metrized complexes of curves] \label{par:metrized.complex.defn}
Recall that $k$ is an algebraically closed field.  A metrized complex of curves is, roughly speaking,
an augmented metric graph $\Gamma$ together with a
chosen vertex set $V(\Gamma)$ and a marked algebraic $k$-curve of
genus $g(p)$ for each finite vertex $p\in V(\Gamma)$.  More precisely:

\begin{defn}\label{def:metr cplx}
  A \emph{$\Lambda$-metrized complex of $k$-curves} consists of the following data:
\begin{itemize}
\item An augmented $\Lambda$-metric graph $\Gamma$ equipped with the choice of a distinguished vertex set $V(\Gamma)$.
\item For every finite vertex $p$ of $\Gamma$, a smooth, proper, connected $k$-curve $C_p$ of genus $g(p)$.
\item An injective function $\red_p: T_p(\Gamma)\to C_p(k)$, 
  called the \emph{reduction map}.
  We call the image of $\red_p$ the set of {\em marked points} on $C_p$.
\end{itemize}
\end{defn}

\begin{eg} \label{ex:metrized.complex}
Figure~\ref{fig:ex cplx} depicts a particular metrized complex of curves over $\CC$. 
At each finite vertex $p$ there is an associated compact Riemann surface $C_p$, together
with a finite set of marked points (in red). 
\begin{figure}[h]
\includegraphics[width=6cm,
  angle=0]{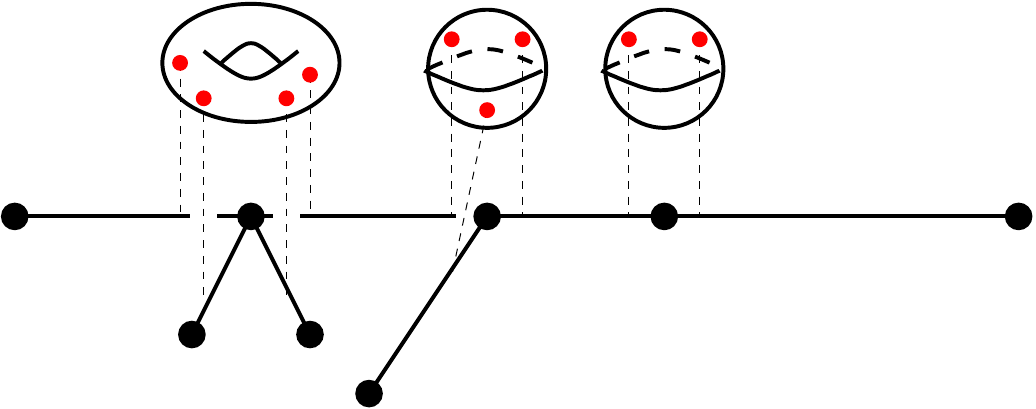}
\caption{}\label{fig:ex cplx}
\end{figure}
\end{eg}

\begin{defn} \label{def:hmmc}
  A \emph{harmonic morphism} of metrized complexes of curves
  consists of a harmonic 
  $(V(\Gamma'),V(\Gamma))$-morphism $\phi:\Gamma'\to\Gamma$ of
  augmented metric graphs, and for every
  finite vertex $p'$ of $\Gamma'$ with $d_{p'}(\phi)>0$
  a finite morphism of curves
  $\phi_{p'}:C'_{p'}\to C_{\phi(p')}$, satisfying the
  following compatibility conditions:
  \begin{enumerate}
  \item For every finite vertex $p'\in V(\Gamma')$ 
    and every tangent direction $v'\in T_{p'}(\Gamma')$
    with $d_{v'}(\phi)>0$, we have 
    $\phi_{p'}(\red_{p'}(v')) = \red_{\phi(p')}(d\phi(p')(v'))$, and
    the ramification degree of $\phi_{p'}$ at
    $\red_{p'}(v')$ is equal to $d_{v'}(\phi)$.

  \item For every finite vertex $p'\in V(\Gamma')$ with 
    $d_{p'}(\phi)>0$, every tangent direction
    $v\in T_{\phi(p')}(\Gamma)$,
    and every point 
    $x'\in\phi_{p'}\inv(\red_{\phi(p')}(v))\subset C'_{p'}(k)$,
    there exists $v'\in T_{p'}(\Gamma')$ such that
    $\red_{p'}(v') = x'$.

  \item For every finite vertex $p'\in V(\Gamma')$ with 
    $d_{p'}(\phi)>0$ we have
    $d_{p'}(\phi) = \deg(\phi_{p'})$.
  \end{enumerate}
\end{defn}

A harmonic morphism of metrized complexes of curves is called \emph{finite} if the underlying harmonic morphism of (augmented) metric graphs is finite.
Metrized complexes of curves with harmonic morphisms between them form a category, with finite harmonic morphisms giving rise to a subcategory. 

If $\Gamma$ has at least one edge, then~(3) follows 
from~(1) and~(2),
since the sum of the ramification degrees of a
finite morphism of curves along any fiber is equal to the degree of the
morphism.

\smallskip 

Let $\phi: C'\to C$ be a finite morphism of smooth, proper, connected
$k$-curves.  Recall that $\phi$ is 
\emph{tamely ramified} if either $\chr(k)=0$ 
or the ramification degree of $\phi$ at every
closed point of $C'$ is prime to the characteristic of $k$.

\begin{rem} \label{rem:tamely.ram}
  \begin{enumerate}
  \item Let $\phi:C'\to C$ be a finite morphism.  Let $D$ be the
    (finite) set of branch points, let $D' = \phi\inv(D) \subset C'(k)$,
    and let $U = C\setminus D$ and $U' = C'\setminus D'$.  Then $\phi$ is
    tamely ramified if and only if $\phi|_{U'}: U'\to U$ is a tamely ramified
    cover of $U$ over $C$ relative to $D$ in the sense 
    of~\cite[XIII.2.1.3]{SGA1}: see the discussion in \S2.0 of {\it loc.\ cit.}
  \item A tamely ramified morphism is separable (i.e.\ generically \'etale).
  \item The Riemann--Hurwitz formula applies to a tamely ramified morphism.
  \end{enumerate}
\end{rem}

\begin{defn}\label{def:tame.covering.metrized.complexes}
  Let $\phi: \cC' \to \cC$ be a finite harmonic
  morphism of metrized complexes of curves.  We say that $\phi$ is a
  \emph{tame harmonic morphism} if either $\chr(k) = 0$,
  or $\chr(k) = p > 0$ and
  $\phi_{p'}$ is tamely ramified for all finite vertices $p'\in V(\Gamma')$. 
  We call $\phi$ a {\em tame covering} if in addition
  $\phi : \Gamma'\to\Gamma$ is a generically \'etale finite morphism of
  augmented metric graphs. 
\end{defn}

\smallskip
Let $\phi:\cC'\to\cC$ be a tame harmonic morphism.  Then
each $\phi_{p'}:C_{p'}\to C_{\phi(p')}$ is a separable morphism.
Note that if $\chr(k) = p$ then $d_{e'}(\phi)$ is prime to $p$ for all
edges $e'$ of $\Gamma'$.  We have the following kind of converse
statement: 

\begin{prop} \label{prop:tame.visible.ram}
  Let $\phi: \cC'\to\cC$ be a finite harmonic morphism of metrized complexes of
  curves.  Suppose that $\phi$ is generically \'etale and that
  $\phi_{p'}:C_{p'}\to C_p$ is a separable morphism of curves
  for all finite vertices $p'$ of $\Gamma'$.  Let $p'$ be a
  finite vertex of $\Gamma'$, and suppose that 
  $\phi_{p'}$ is ramified at a point $x'\in C'_{p'}(k)$.  Then 
  there exists $v'\in T_{p'}(\Gamma')$ such that
  $\red_{p'}(v') = x'$.  
  In particular, if $\chr(k) = p > 0$ and
  $d_{e'}(\phi)$ is prime to $p$ for all edges $e'$ of
  $\cC'$, then $\phi_{p'}$ is tamely ramified and $\phi$ is a tame
  covering. 
\end{prop}

\pf Let $p = \phi(p')$.
Since $\phi$ is generically \'etale we have
\[ 0 = R_{p'} = d_{p'}(\phi) (2-2g(p)) - (2-2g(p')) - 
\sum_{v'\in T_{p'}(\Gamma')} (d_{v'}(\phi) - 1). \]
Since $g(p) = g(C_{p})$, $g(p') = g(C'_{p'})$, and 
$d_{p'}(\phi)$ is the degree of $\phi_{p'}: C'_{p'}\to C_{p}$, this gives
\[ \sum_{v' \in T_{p'}(\Gamma')} (d_{v'}(\phi)-1) 
= \deg(\phi_{p'}) (2-2g(C_{p})) - (2-2g(C'_{p'})). \]
The Riemann--Hurwitz formula as applied to $\phi_{p'}$ says that the right-hand
side of the above equation is equal to the degree of the ramification
divisor of $\phi_{p'}$ (this holds whenever $\phi_{p'}$ is separable), so
the proposition follows because for all 
$v'\in T_{p'}(\Gamma')$, the ramification degree of $\phi_{p'}$ at
$\red_{p'}(v')$ is $d_{v'}(\phi)$.\qed

\begin{rem}
  Loosely speaking, Proposition~\ref{prop:tame.visible.ram} says that all
  ramification of all morphisms $\phi_{p'}$ is ``visible'' in the
  edge degrees of the underlying morphism of metric graphs.
  It follows from Proposition~\ref{prop:tame.visible.ram} that if
  $\Gamma' = \{p'\}$, then $\phi_{p'}$ is \'etale.
\end{rem}

\section{Analytic curves and their skeleta}\label{sec:skeleta}

In this section we recall the definition and basic properties of the
skeleton of the analytification of an algebraic curve as outlined
in~\cite{temkin:berkovich_spaces} and elaborated
in~\cite[\S5]{bpr:trop_curves}.  These foundational results 
are well-known to experts, and appear in various guises in the literature.
See for example 
Berkovich~\cite{berkovich:analytic_geometry,berkovich:locallycontractible1}, 
Thuillier~\cite{thuillier:thesis}, and Ducros~\cite{ducros:triangulations}.
In this section we also show that the skeleton is naturally a
metrized complex of curves, and prove that any metrized complex arises as
the skeleton of a curve.  Finally, we introduce the important notion of a
triangulated punctured curve, which is essentially a punctured curve along
with the data of a skeleton with a distinguished set of vertices.

Recall that $K$ is an algebraically closed field which is
complete with respect to a nontrivial non-Archimedean valuation 
$\val:K\to\R\cup\{\infty\}$.  Its valuation ring is $R$, its
maximal ideal is $\fm_R$, its residue field is $k = R/\fm_R$, and
its value group is $\Lambda = \val(K^\times)\subset\R$.

By an \emph{analytic curve} we mean a strictly $K$-analytic
space $X$ of pure dimension $1$.  For our purposes the most important
example of an analytic curve is the analytification of a smooth,
connected, projective algebraic $K$-curve.
If $X$ is an analytic curve, then we define 
$\bH(X) = X\setminus X(K)$ as in~\cite[3.10]{bpr:trop_curves}.

\paragraph[Analytic building blocks] \label{par:balls.annuli}
We identify the set underlying the analytification of the affine line with
the set of multiplicative seminorms $\|\cdot\|:K[t]\to\R_{\geq 0}$ which
extend the absolute value on $K$.
We let $\val:\A^{1,\an}\to\R\cup\{\infty\}$ denote the valuation map
\[ \val(x) = -\log |x| = -\log \|t\|_x. \]
For $a\in K^\times$ the \emph{standard closed ball of radius $|a|$} is the
affinoid domain
\[ \B(a) = \val\inv([\val(a),\infty]) \subset \A^{1,\an} \]
and the \emph{standard open ball of radius $|a|$} is the open analytic domain
\[ \B(a)_+ = \val\inv((\val(a),\infty]) \subset \A^{1,\an}. \]
Note that scaling by $a$ gives isomorphisms $\B(1)\isom\B(a)$ and
$\B(1)_+\isom\B(a)_+$.  We call $\B(1)$ (resp.\ $\B(1)_+$) the 
\emph{standard closed} (resp.\ \emph{open}) \emph{unit ball}.
For $a\in K^\times$ with $\val(a)\geq 0$ the
\emph{standard closed annulus of modulus $\val(a)$} is the affinoid
domain
\[ \bS(a) = \val\inv([0,\val(a)]) \subset \A^{1,\an} \]
and if $\val(a) > 0$ the \emph{standard open annulus of modulus $\val(a)$}
is the open analytic domain
\[ \bS(a)_+ = \val\inv((0,\val(a))) \subset \A^{1,\an}. \]
The \emph{standard punctured open ball} is the open analytic domain
\[ \bS(0)_+ = \val\inv((0,\infty)) = \B(1)_+\setminus\{0\} \subset
\A^{1,\an}. \]

By a \emph{closed unit ball} (resp.\ \emph{open unit ball}, resp.\
\emph{closed annulus}, resp.\ \emph{open annulus}, 
resp.\ \emph{punctured open ball}) we will mean a $K$-analytic space
which is isomorphic to the standard closed unit ball 
(resp.\ the standard open unit ball, resp.\ the standard closed annulus of
modulus $\val(a)$ for some $a\in R\setminus\{0\}$,
resp.\ the standard open annulus of modulus $\val(a)$ for some
$a\in\fm_R\setminus\{0\}$, resp.\ the standard punctured open ball).  A
\emph{generalized open annulus} is a $K$-analytic space which is either an
open annulus or a punctured open ball.  
It is a standard fact that the modulus of a closed (resp.\ open) annulus
is an isomorphism invariant; see for
instance~\cite[Corollary~5.6]{bpr:trop_curves}.  We define
the modulus of a punctured open ball to be $\infty$.

Let $a\in\fm_R$ (so $a=0$ is allowed). There is a natural section 
$\sigma:(0,\val(a))\to\bS(a)_+$ of the valuation map 
$\val:\bS(a)_+\to(0,\val(a))$ sending $r$ to the maximal point of the
affinoid domain $\val\inv(r)$ if $r\in\Lambda$, resp.\ the only point of
$\val\inv(r)$ if $r\notin\Lambda$.
The \emph{skeleton} of $\bS(a)_+$ is defined to be the image of $\sigma$
and is denoted $\Sigma(\bS(a)_+) = \sigma((0,\val(a)))$.  
We will always identify the skeleton of the standard generalized open
annulus $\bS(a)_+$ with the interval or ray $(0,\val(a))$.
It follows from~\cite[Proposition~5.5]{bpr:trop_curves} that the skeleton
of a standard open annulus or standard punctured open ball is an
isomorphism invariant, so if $A$ is a generalized open annulus, then we can
define the \emph{skeleton $\Sigma(A)$} of $A$ to be the image of the
skeleton of a standard open annulus or standard punctured open ball
$\bS(a)_+$ under any isomorphism $\bS(a)_+\isom A$.

Let $T$ be a metric space and let $x,y\in T$. 
A \emph{geodesic segment from $x$ to $y$} is the image of a locally
isometric embedding $[a,b]\inject T$ with $a\mapsto x$ and $b\mapsto y$.
We often identify a geodesic segment with its image in $T$.  Recall that
an \emph{$\R$-tree} is a metric space $T$ with the following properties:
\begin{enum}
\item For all $x,y\in T$ there is a unique geodesic segment $[x,y]$ from
  $x$ to $y$.
\item For all $x,y,z\in T$, if $[x,y]\cap[y,z]=\{y\}$, then
  $[x,z] = [x,y]\cup[y,z]$.
\end{enum}
See~\cite[Appendix~B]{baker_rumely:book}.  It is proved in \S1.4 of loc.\
cit.\ that $\bH(\B(1))$ is naturally an $\R$-tree.  Since any
path-connected subspace of an $\R$-tree is an $\R$-tree as well, if $X$ is
a standard open annulus or standard (punctured) open ball then
$\bH(X)$ is an $\R$-tree.
It is proved in~\cite[Corollary~5.61]{bpr:trop_curves} that the metric
structure on $\bH(X)$  is an isomorphism invariant, so the same statement
applies to open balls and generalized open annuli.
For $a\in\fm_R$ the section $\sigma:(0,\val(a))\to\bS(a)_+$ is an isometry.

It also follows from the results in~\cite[\S1.4]{baker_rumely:book} that
for any $x\in\bH(\B(1))$ and any type-$1$ point $y\in\B(1)$, there is a
unique continuous injection $\alpha:[0,\infty]\inject\B(1)$ with
$\alpha(0) = x$ and $\alpha(\infty) = y$, 
$\alpha([0,\infty))\subset\bH(\B(1))$, and such that $\alpha$ is an
isometry when restricted to $[0,\infty)$.  We let $[x,y]$ denote the image
of $\alpha$ and we call $[x,y]$ the 
\emph{geodesic segment from $x$ to $y$}.  Similarly, if $x$ and $y$ are
both type-$1$ points then there is a unique continuous injection
$\alpha:[-\infty,\infty]\to\B(1)$ with $\alpha(-\infty)=x$ and
$\alpha(\infty)=y$, $\alpha(\R)\subset\bH(\B(1))$, and such that
$\alpha|_\R$ is an isometry; the image of $\alpha$ is called the
\emph{geodesic segment from $x$ to $y$} and is denoted $[x,y]$.
Restricting to a suitable analytic subdomain of $\B(1)$ allows us to
define geodesic segments between any two points of an open ball or
generalized open annulus.

\paragraph[Open balls]
The closure of $\B(1)_+$ in $\B(1)$ consists of $\B(1)_+$ and a single
type-$2$ point $x$, called the \emph{end} of $\B(1)_+$.  The end is the
Shilov boundary point of $\B(1)$: see for
instance the proof of~\cite[Lemma~5.16]{bpr:trop_curves}.  The valuation
map $\val:\B(1)_+\to(0,\infty]$ extends to a continuous map
$\val:\B(1)_+\cup\{x\}\to[0,\infty]$.  The set
$\B(1)_+\cup\{x\}$ is path-connected and compact, being a closed subspace
of the compact space $\B(1)$.  (In fact $\B(1)_+\cup\{x\}$ is the one-point
compactification of $\B(1)_+$.)  For any $y\in\B(1)_+$ the geodesic segment
$[x,y]\subset\B(1)$ is contained in $\B(1)_+\cup\{x\}$.  

\begin{lem} \label{lem:extend.end.ball}
  Let $X = \Spec(A)$ be an irreducible affine $K$-curve and let
  $\phi:\B(1)_+\to X^\an$ be a morphism with finite fibers.  Then $\phi$ 
  extends in a unique way to a continuous map
  $\B(1)_+\cup\{x\}\to X^\an$, and the image of $x$ is a type-$2$ point of
  $X^\an$. 
\end{lem}

\pf Let $f\in A$ be nonzero and define $F: (0,\infty)\to\R$ by
$F(r) = -\log |f\circ\phi(\sigma(r))|$, where
$\sigma:(0,\infty)\to\bS(0)_+\subset\B(1)_+$ is the canonical section of
$\val$.  By~\cite[Proposition~5.10]{bpr:trop_curves}, 
$F$ is a piecewise affine function which is differentiable away from
finitely many points $r\in\Lambda\cap(0,\infty)$, and for any point 
$r\in(0,\infty)$ at which $F$ is differentiable, the derivative of
$F$ is equal to the number of zeros $y$ of $f\circ\phi$ with 
$\val(y) > r$.  It follows that 
$F(\Lambda\cap(0,\infty))\subset\Lambda$.  Since $f\circ\phi$ has finitely
many zeros, the limit $\lim_{r\to 0} F(r)$ exists and is contained in
$\Lambda$.  We define $\|f\|_0 = \exp(-\lim_{r\to 0} F(r))$.  It is easy to
see that $f\mapsto\|f\|_0$ is a multiplicative seminorm on $A$ extending
the absolute value on $K$, so $\|\cdot\|_0$ is a point in $X^\an$.  Define
$\phi(x) = \|\cdot\|_0$.  One shows as in the proof
of~\cite[Lemma~5.16]{bpr:trop_curves} that $\phi$ is continuous on
$\B(1)_+\cup\{x\}$.

It remains to show that $\phi(x) = \|\cdot\|_0$ is a type-$2$ point of
$X^\an$.  Let $f\in A$ be a non-constant function that has a zero on
$\phi(\B(1)_+)$.  Since
$-\log \|f\|_0\in\Lambda$, there exists $\alpha\in K^\times$ such that
$\|\alpha f\|_0 = 1$; replacing $f$ by $\alpha f$, we may and do assume
that $\|f\|_0 = 1$.  Since $\phi\circ f$ has a zero, by the above we have that
$F = -\log|f\circ\phi\circ\sigma|$ is monotonically increasing, so
$F(r)>0$ for all $r\in(0,\infty)$.  For any
$r\in(0,\infty)$, the maximal point of $\val\inv(r)\subset\B(1)_+$ is
equal to $\sigma(r)$, so for any $y\in\B(1)_+$ such that $\val(y) = r$ we
have $F(r)\leq -\log |f\circ\phi(y)|$.  It follows that 
$|f\circ\phi(y)| < 1$ for all $y\in\B(1)_+$, so
$\phi(f(\B(1)_+))\subset\B(1)_+$.  Since $\B(1)_+$ 
is dense in $\B(1)_+\cup\{x\}$, the image of $x$ under $\phi\circ f$ is
contained in the closure $\B(1)_+\cup\{x\}$ of $\B(1)_+$ in $\B(1)$.
Since $1 = \|f\|_0 = |f\circ\phi(x)|$, the point
$f\circ\phi(x)\notin\B(1)_+$, so $f\circ\phi(x) = x$, which is a type-$2$
point of $\A^{1,\an}$.  Therefore $\phi(x)$ is a type-$2$ point of
$X^\an$.\qed 

\smallskip
Applying Lemma~\ref{lem:extend.end.ball} to a morphism 
$\phi:\B(1)_+\to\B(1)_+\subset\A^{1,\an}$, we see that any automorphism of
$\B(1)_+$ extends (uniquely) to a homeomorphism
$\B(1)_+\cup\{x\}\to\B(1)_+\cup\{x\}$, so it makes sense to speak of the
\emph{end} of any open ball.  If $B$ is an open ball with end $x$, we let
$\bar B$ denote $B\cup\{x\}$.

Let $B$ be an open ball with end $x$.  We define a partial ordering on 
$\bar B$ by declaring $y\leq z$ if $y\in[x,z]$; again
see~\cite[\S1.4]{baker_rumely:book}.  Equivalently, for $y,z\in B$ we have
$y\leq z$ if and only if $|f(y)|\geq|f(z)|$ for all analytic functions $f$
on $B$.  The following lemma is proved as
in~\cite[Proposition~5.27]{bpr:trop_curves}.  

\begin{lem} \label{lem:ball.subtract.point}
  Let $B$ be an open ball and let $y\in B$ be a type-$2$
  point.  Then $B\setminus\{y\}$ is a disjoint union of the open annulus
  $A = \{ z\in B~:~z\not\geq y \}$ with infinitely many open balls, and
  $B_y = \{ z\in B~:~z\geq y \}$ is an affinoid subdomain of $B$ 
  isomorphic to the closed ball $\B(1)$.
\end{lem}

\paragraph[Open annuli]
Let $a\in\fm_R\setminus\{0\}$.  The closure of $\bS(a)_+$ in $\B(1)$
consists of $\bS(a)_+$ and the two type-$2$ points
$x=\sigma(0)$ and $y=\sigma(\val(a))$, called the \emph{ends} of
$\bS(a)_+$: again see the proof of~\cite[Lemma~5.16]{bpr:trop_curves}.
The valuation map $\val:\bS(a)_+\to(0,\val(a))$ extends to
a continuous map $\val:\bS(a)_+\cup\{x,y\}\to[0,\val(a)]$, and for any
$z\in\bS(a)_+$ the geodesic segments $[x,z]$ and $[y,z]$ are contained in
$\bS(a)_+\cup\{x,y\}$. 
The following lemma is proved in the same way as the first part of
Lemma~\ref{lem:extend.end.ball}.

\begin{lem} \label{lem:extend.ends.annulus}
  Let $X = \Spec(A)$ be an irreducible affine $K$-curve, let
  $a\in\fm_R\setminus\{0\}$, and let 
  $\phi:\bS(a)_+\to X^\an$ be a morphism with finite fibers.  
  Let $x = \sigma(0)$ and $y = \sigma(\val(a))$ be the ends of
  $\bS(a)_+$.  Then $\phi$ 
  extends in a unique way to a continuous map
  $\bS(a)_+\cup\{x,y\}\to X^\an$.
\end{lem}

It follows from Lemma~\ref{lem:extend.ends.annulus} that any automorphism
of $\bS(a)_+$ extends to a homeomorphism 
$\bS(a)_+\cup\{x,y\}\isom\bS(a)_+\cup\{x,y\}$, so it makes sense to speak
of the \emph{ends} of any open annulus.  If $A$ is an open annulus with
ends $x,y$, then we let $\bar A$ denote the compact space $A\cup\{x,y\}$.

The closure of the punctured open ball $\bS(0)_+$ in $\B(1)$ is equal to 
$\bS(0)_+\cup\{0,x\}$, where $x$ is the end of $\B(1)_+$.  We define $x$
to be the \emph{end} of $\bS(0)_+$ and $0$ to be the \emph{puncture}.  As
above, the end and puncture of $\bS(0)_+$ are isomorphism invariants, so
it makes sense to speak of the \emph{end} and the \emph{puncture} of any
punctured open ball.  If $A$ is a punctured open ball with end $x$ and
puncture $y$, then we let $\bar A$ denote the compact space $A\cup\{x,y\}$.

Let $A$ be a generalized open annulus. Each connected
component of $A\setminus\Sigma(A)$ is an open ball~\cite[Lemma~5.12]{bpr:trop_curves}, and if 
$B$ is such a connected
component 
with end $x$, then the inclusion
$B\inject A$ extends to an inclusion $\bar B\inject A$ with $x$ mapping
into $\Sigma(A)$; the image of $\bar B$ is the closure of $B$ in $A$.  We
define the retraction to the skeleton $\tau: A\to\Sigma(A)$ by fixing
$\Sigma(A)$ and sending each connected component of $A\setminus\Sigma(A)$
to its end.  If $A = \bS(a)_+$, then the retraction $\tau$
coincides with $\sigma\circ\val: \bS(a)_+\to\Sigma(\bS(a)_+)$; in
particular, if $(r,s)\subset(0,\val(a))=\bS(a)_+$, then 
$\tau\inv((r,s)) = \val\inv((r,s))$.

\paragraph[The skeleton of a curve]
Let $X$ be a smooth, connected, proper algebraic $K$-curve and let
$D\subset X(K)$ be a finite set of closed points.  The set
$\bH(X^\an)$ is natually a metric space~\cite[Corollary~5.62]{bpr:trop_curves},
although the metric topology on $\bH(X^\an)$ is much finer than the
topology induced by the topology on the $K$-analytic space $X^\an$.  The
metric on $\bH(X^\an)$ is locally modeled on an $\R$-tree~\cite[Proposition~5.63]{bpr:trop_curves}. 

\begin{defn}\label{def:semistable vertex set}
  A \emph{semistable vertex set of $X$} is a finite set $V$ of
  type-$2$ points of $X^\an$ such that $X^\an\setminus V$ is a
  disjoint union of open balls and finitely many open annuli.  A
  \emph{semistable vertex set of $(X,D)$} is a semistable
  vertex set of $X$ such that the points of $D$ are contained in 
  distinct open ball connected components of $X^\an\setminus V$.
\end{defn}

\smallskip
It is a consequence of the semistable reduction theorem of Deligne-Mumford
that there exist semistable vertex sets of $(X,D)$: see 
\cite[Theorem~5.49]{bpr:trop_curves}.
In the sequel it will be convenient to consider a curve along with a
choice of semistable vertex set, so we give such an object a name.

\begin{defn}\label{def:triangulated curve}
  A \emph{triangulated punctured curve} $(X,V\cup D)$ is a smooth,
  connected, proper algebraic $K$-curve $X$ equipped with a finite set
  $D\subset X(K)$ of \emph{punctures} and a semistable vertex set $V$ of
  $(X,D)$.   
\end{defn}

\begin{rem}
  This terminology is loosely based on that used in~\cite{ducros:triangulations}
  as well as the forthcoming book of Ducros on analytic
  spaces and analytic curves.  Strictly speaking, what we have defined
  should be called a \emph{semistably triangulated punctured curve}, but
  as these are the only triangulations that we consider, we will not need
  the added precision.  
\end{rem}

\paragraph \label{par:semistable.decomposition}
Let $(X,V\cup D)$ be a triangulated punctured curve, so
\begin{equation} \label{eq:semistable.decomposition}
  X^\an\setminus(V\cup D)
  = A_1\cup\cdots\cup A_n\cup\bigcup_\alpha B_\alpha, 
\end{equation}
where each $A_i$ is a generalized open annulus and $\{B_\alpha\}$ is an
infinite collection of open balls.  
The \emph{skeleton of $(X,V\cup D)$} is the subset
\[ \Sigma(X,V\cup D) = V\cup D\cup
\Sigma(A_1)\cup\cdots\cup\Sigma(A_n). \]
(This set is denoted $\hat\Sigma(X,V\cup D)$ in~\cite{bpr:trop_curves}.)
For each $i$ and each $\alpha$ the inclusions
$\bH(A_i)\inject\bH(X^\an)$ and $\bH(B_\alpha)\inject\bH(X^\an)$ are local
isometries~\cite[Proposition~5.60]{bpr:trop_curves}.  For each open ball
$B_\alpha$ the map $B_\alpha\inject X^\an$ extends to an inclusion 
$\bar B_\alpha\inject X^\an$ sending the end of $B_\alpha$ to a point
$x_\alpha\in V$.  We say that $B_\alpha$ is \emph{adjacent} to $x_\alpha$.
For each open annulus $A_i$ the map $A_i\inject X^\an$ extends to a
continuous map $\bar A_i\inject X^\an$ sending the ends of $A_i$
to points $x_i,y_i\in V$.  We say that
$A_i$ is \emph{adjacent} to $x_i$ and $y_i$.  The
length $\ell(e_i)\in\Lambda$ of the geodesic segment
$e_i = \Sigma(A_i)\cup\{x_i,y_i\}$ is then the modulus of $A_i$.
We say that $V$ is \emph{strongly semistable} if $x_i\neq y_i$ for each
open annulus $A_i$.  For each generalized
open annulus $A_i$ the map $A_i\inject X^\an$ extends to a continuous map
$\bar A_i\inject X^\an$ sending the end of $A_i$ to a point
$x_i\in V$ and sending the puncture
to a point $y_i\in D$.  We say that $A_i$ is \emph{adjacent} to
$x_i$ and $y_i$. We define the length of 
$e_i = \Sigma(A_i)\cup\{x_i,y_i\}$ to be $\ell(e_i) = \infty$.

The skeleton $\Sigma = \Sigma(X,V\cup D)$ naturally has the structure of a
$\Lambda$-metric graph with distinguished finite vertex set
$V_f(\Sigma) = V$, infinite vertex set 
$V_\infty(\Sigma) = D$, and edges $\{e_1,\ldots,e_n\}$ as above.  
Note that the $\Lambda$-points of $\Sigma$ are exactly the type-$2$ points
of $X^\an$ contained in $\Sigma$.
For $x\in V$ the residue field $\td\sH(x)$ of the
completed residue field $\sH(x)$ at the type-$2$ point $x$ is a finitely
generated field extension of $k$ of transcendence degree $1$; we let $C_x$
be the smooth $k$-curve with function field $\td\sH(x)$.  For $x\in V$ we
let $g(x)$ be the genus of $C_x$, and for $x\in\Sigma\setminus V$ we set
$g(x) = 0$.  These extra data give $\Sigma$ the structure of an augmented
$\Lambda$-metric graph.  
(In \parref{par:complex.from.curves} we will see that $\Sigma$ is in fact naturally a metrized complex of $k$-curves.)
By the genus formula~\cite[(5.45.1)]{bpr:trop_curves}, the
genus of $X$ is equal to the genus of the augmented $\Lambda$-metric graph
$\Sigma$. 

The open analytic domain
$X^\an\setminus\Sigma$ is isomorphic to an infinite disjoint union of
open balls~\cite[Lemma~5.18(3)]{bpr:trop_curves}.  If $B$ is a connected component of $X^\an\setminus\Sigma$,
then the inclusion $B\inject X^\an$ extends to a map $\bar B\inject X^\an$
sending the end of $B$ to a point of $\Sigma$.  We define the
\emph{retraction} $\tau = \tau_\Sigma: X^\an\to\Sigma$ by fixing
$\Sigma$ and sending a point $x\in X^\an$ not in $\Sigma$ to the end of
the connected component of $x$ in $X^\an\setminus\Sigma$.  This is a
continuous map.  If $x\in B_\alpha$ is in an open ball connected component
of $X^\an\setminus(V\cup D)$, then $\tau(x)\in V$ is the end of $B_\alpha$,
and if $x\in A_i$, then $\tau(x)$ coincides with the image of $x$ under the
retraction map $\tau: A_i\to\Sigma(A_i)$.

Here we collect some additional facts about skeleta
from~\cite[\S5]{bpr:trop_curves}.

\begin{prop} \label{prop:skeleton.properties}
  Let $(X,V\cup D)$ be a triangulated punctured curve with skeleton
  $\Sigma = \Sigma(X, V\cup D)$.
  \begin{enumerate}
  \item The skeleton $\Sigma$ is the set of points in $X^\an$ that do not
    admit an open neighborhood isomorphic to $\B(1)_+$ and disjoint from
    $V\cup D$.

  \item Let $V_1$ be a semistable vertex set of $(X,D)$ such
    that $V_1\supset V$.  Then 
    $\Sigma(X,V_1\cup D)\supset\Sigma(X,V\cup D)$
    and $\tau_{\Sigma(X,V\cup D)} = \tau_{\Sigma(X,V\cup D)} \circ
    \tau_{\Sigma(X,V_1\cup D)}$. 

  \item Let $W\subset X^\an$ be a finite set of type-$2$ points.
    Then there exists a semistable vertex set of $(X,D)$
    containing $V\cup W$.

  \item Let $W\subset\Sigma$ be a finite set of type-$2$ points.  Then
    $V\cup W$ is a semistable vertex set of $(X,D)$ and
    $\Sigma(X,V\cup D) = \Sigma(X,V\cup W\cup D)$.

  \item Let $x,y\in\Sigma\cap\bH(X^\an)$.  Then any geodesic segment from
    $x$ to $y$ in $\bH(X^\an)$ is contained in $\Sigma$.

  \item Let $x,y\in X^\an$ be points of type-$2$ or $3$ and let
    $[x,y]$ be a geodesic segment from $x$ to $y$ in $\bH(X^\an)$.  Then
    there exists a semistable vertex set $V_1$ of $(X,D)$
    such that $V_1\supset V$ and  $[x,y]\subset\Sigma(X,V_1\cup D)$.
  \end{enumerate}
\end{prop}

\begin{defn}
  A \emph{skeleton of $(X,D)$} is a subset of $X^\an$ of the form
  $\Sigma = \Sigma(X,V\cup D)$ for some semistable vertex set $V$ of
  $(X,D)$.  Such a semistable vertex set $V$ is called a
  \emph{vertex set for $\Sigma$}.  A \emph{skeleton of $X$} is a skeleton
  of $(X,\emptyset)$.
\end{defn}

The augmented $\Lambda$-metric graph structure of 
a skeleton $\Sigma$ of $(X,D)$ does not depend on the choice of vertex set
for $\Sigma$.

\paragraph[Modifying the skeleton] 
Let $X$ be a smooth, proper, connected $K$-curve and let $D\subset X(K)$
be a finite set.  Let $\Sigma$ be a skeleton of $(X,D)$, let 
$y\in X^\an\setminus\Sigma$, let $B$ be the connected component of
$X^\an\setminus\Sigma$ containing $y$, and let $x = \tau(y)\in\Sigma$ be
the end of $B$.  If $y\in\bH(X^\an)$, then the geodesic segment
$[x,y]\subset\bar B$ is the unique geodesic segment in $\bH(X^\an)$
connecting $x$ and $y$.  If 
$y\in X(K)$, then we define $[x,y]$ to be the geodesic segment
$[x,y]\subset\bar B$ as in~\parref{par:balls.annuli}.
The following strengthening of
Proposition~\ref{prop:skeleton.properties}(3) will be important in the
sequel. 

\begin{lem} \label{lem:larger.vertex.set}
  Let $V$ be a semistable vertex set of $(X,D)$, let
  $\Sigma = \Sigma(X,V\cup D)$, let $W\subset X^\an$ be a finite set
  of type-$2$ points, and let $E\subset X(K)$ be a finite set of type-$1$
  points. 
  \begin{enumerate}
  \item There exists a minimal semistable vertex set $V_1$ of $(X,D\cup E)$
    which contains $V\cup W$, in the sense that any other
    such semistable vertex set contains $V_1$.
  \item Let $B$ be a connected component of $X^\an\setminus\Sigma$ with
    end $x$ and let $y_1,\ldots,y_n$ be the points of $(W\cup E)\cap B$.
    Then 
    \[ \Sigma(X,V_1\cup D\cup E)\cap\bar B = [x,y_1]\cup\cdots\cup[x,y_n] \]
    if $n > 0$, and $\Sigma(X,V_1\cup D\cup E)\cap\bar B = \{x\}$ otherwise.
    Therefore 
    \[ \Sigma(X,V_1\cup D\cup E) = \Sigma\cup\bigcup_{y\in W\cup E}
    [\tau_\Sigma(y),\,y]. \]
  \item The skeleton $\Sigma_1 = \Sigma(X,V_1\cup D\cup E)$ is minimal in
    the sense that any other skeleton containing $\Sigma$ and $W\cup E$
    must contain $\Sigma_1$.
  \end{enumerate}
\end{lem}

\pf To prove the first part we may assume that $W = \{y\}$ consists of a
single type-$2$ point not contained in $\Sigma$ and $E=\emptyset$, 
or $E = \{y\}$ is a single type-$1$ point not contained in $D$ and
$W=\emptyset$.  In the first case one sees using 
Lemma~\ref{lem:ball.subtract.point} that $V_1 = V\cup\{y,\tau(y)\}$ is the
minimal semistable vertex set of $(X,D)$
containing $V$ and $y$, and in the second case $V_1 = V\cup\{\tau(y)\}$
is the minimal semistable vertex set of $(X,D\cup\{y\})$
containing $V$.  For $W$ and $E$ arbitrary, let $V_1$ be the minimal 
semistable vertex set of $(X,D\cup E)$ containing $V\cup W$
and let $\Sigma_1 = \Sigma(X,V_1\cup D)$.  
If $E = \emptyset$, then it is clear from
Proposition~\ref{prop:skeleton.properties}(5) that 
$[x,y_i]\subset\Sigma_1$ for each $i$; the other inclusion is proved by
induction on $n$, adding one point at a time as above.
The case $E\neq\emptyset$ is similar and is left to the reader.
The final assertion follows easily from the first two and
Proposition~\ref{prop:skeleton.properties}.\qed

\begin{rem} \label{rem:skeleta.modifications}
The above lemma shows in particular that
 the metric graph $\Sigma(X,V_1\cup D\cup E)$ is 
obtained from $\Sigma(X,V\cup D)$ by a sequence of tropical
modifications and their inverses (see Definition~\ref{defn:modification} and Example \ref{ex:add finite edge}). 
It is easy to see that any tropical modification of $\Sigma(X,V\cup D)$ is of the form  
$\Sigma(X,V_1\cup D\cup E)$ for a semistable vertex set $V_1$ which contains $V$ and a finite subset $E \subset X(K)$.
\end{rem}

\paragraph[The minimal skeleton]
The \emph{Euler characteristic} of $X\setminus D$ is defined to be
$\chi(X\setminus D) = 2 - 2g(X) - \#D$, where $g(X)$ is the genus of $X$.
We say that $(X,D)$ is \emph{stable} if $\chi(X\setminus D) < 0$.
A semistable vertex set $V$ of $(X,D)$ is \emph{minimal} if
there is no semistable vertex set of $(X,D)$ properly
contained in $V$, and $V$ is \emph{stable} if there is no $x\in V$ of
genus zero (resp.\ one) and valency less than $3$ (resp.\ one)
in $\Sigma(X,V\cup D)$.  
It is clear that minimal semistable vertex sets exist. A skeleton
$\Sigma$ of $(X,D)$ is \emph{minimal} if there is no skeleton
of $(X,D)$ properly contained in $\Sigma$.  

The following consequence of the stable reduction theorem can be found
in~\cite[Theorem~5.49]{bpr:trop_curves}.

\begin{prop} \label{prop:minimal.skeleton}
  Let $V$ be a minimal semistable vertex set of $(X,D)$ and
  let $\Sigma = \Sigma(X,V\cup D)$.
  \begin{enumerate}
  \item If $\chi(X\setminus D)\leq 0$, then $\Sigma$ is the unique minimal
    skleton of $(X,D)$; moreover, $\Sigma$ is equal to the set
    of points of $X^\an\setminus D$ that do not admit an open neighborhood
    isomorphic to $\B(1)_+$ and disjoint from $D$.

  \item If $\chi(X\setminus D) < 0$, 
    then $V$ is the unique minimal semistable
    vertex set of $(X,D)$, $V$ is stable, and 
    \[ V = \{ x\in\Sigma~:~ x\text{ has valency} \geq 3\text{ or genus}
    \geq 1 \}. \]
  \end{enumerate}
\end{prop}

\begin{rem} \label{rem:not.stable}
  We have $\chi(X\setminus D) > 0$ if and only if $X\cong\P^1$ and $D$
  contains at most one point; in this case, any type-$2$ point of $X^\an$
  serves as a minimal semistable vertex set of $(X,D)$, and
  there does not exist a unique minimal skeleton.  If $\chi(X\setminus D) = 0$,
  then either $X\cong\P^1$ and $D$ consists of two points, or 
  $X$ is an elliptic curve and $D$ is empty.  In the first case, the
  minimal skeleton $\Sigma$ of $(X,D)$ is the extended line connecting
  the points of $D$ and any type-$2$ point on $\Sigma$ is a
  minimal semistable vertex set.  If $X$ is an elliptic curve and $X^\an$
  contains a type-$2$ point $x$ of genus $1$, then $\{x\}$ is both the
  unique minimal semistable vertex set and the minimal skeleton.
  Otherwise $X$ is a Tate curve, $\Sigma$ is a circle, and any type-$2$
  point of $\Sigma$ is a minimal semistable vertex set.
  See~\cite[Remark~5.51]{bpr:trop_curves}.
\end{rem}

\paragraph[Tangent directions and the slope formula] \label{par:slope.formula}
As above we let $X$ be a smooth, proper, connected algebraic $K$-curve.
A continuous function $F:X^\an\to\R\cup\{\pm\infty\}$ is called
\emph{piecewise affine} provided that $F(\bH(X^\an))\subset\R$ and
$F\circ\alpha:[a,b]\to\R$ is a piecewise affine function for every
geodesic segment $\alpha:[a,b]\inject\bH(X^\an)$.
To any point $x\in X^\an$ is associated a set $T_x$ of 
\emph{tangent directions} at $x$, defined as the set of germs of geodesic
segments in $X^\an$ beginning at $x$.  If $F:X^\an\to\R\cup\{\pm\infty\}$ is a
piecewise affine function and $v\in T_x$ we denote by $d_v F(x)$ the
outgoing slope of $F$ in the direction $v$.  We say that $F$ is
\emph{harmonic} at a point $x\in X^\an$ provided that there are only
finitely many $v\in T_x$ with $d_v F(x)\neq 0$, and 
$\sum_{v\in T_x} d_v F(x) = 0$.  See~\cite[5.65]{bpr:trop_curves}.

Let $x\in X^\an$ be a type-$2$ point, let $\sH(x)$ be the completed residue
field of $X^\an$ at $x$, and let $\td\sH(x)$ be its residue field, as
in~\parref{par:semistable.decomposition}.  Then 
$\td\sH(x)$ is a finitely generated field extension of $k$ of
transcendence degree $1$, and there is a canonical bijection between the
set $T_x$ of tangent directions 
to $X^\an$ at $x$ and the set $\DV(\td\sH(x)/k)$ 
of discrete valuations on $\td\sH(x)$ which
are trivial on $k$.  For $v\in T_x$ we let 
$\ord_v:\td\sH(x)\to\Z$ denote the corresponding valuation.  
Let $f$ be a nonzero analytic function on $X^\an$ defined on a
neighborhood of $x$, let $c\in K^\times$ be such that
$|c| = |f(x)|$, and let $\td f_x\in\td\sH(x)$
denote the residue of $c\inv f$.  Then $\td f_x$ is defined up to
multiplication by $k^\times$, so for any $v\in T_x$ the integer
$\ord_v(\td f_x)$ is well-defined.

The following theorem is called the slope formula
in~\cite[Theorem~5.69]{bpr:trop_curves}:

\begin{thm} \label{thm:slope.formula}
  Let $f\in K(X)^\times$ be a nonzero rational function on $X$ and let
  $F = -\log |f|: X^\an\to\R\cup\{\pm\infty\}$.  Let $D\subset X(K)$
  contain the zeros and poles of $f$ and let $\Sigma$ be a skeleton of
  $(X,D)$.  Then:
  \begin{enumerate}
  \item $F = F\circ\tau$, where $\tau:X^\an\to\Sigma$ is the retraction.
  \item $F$ is piecewise affine with integer slopes, and $F$ is affine on
    each edge of $\Sigma$ (with respect to a choice of vertex set $V$).
  \item If $x$ is a type-$2$ point of $X^\an$ and $v\in T_x$, then
    $d_v F(x) = \ord_v(\td f_x)$.
  \item $F$ is harmonic at all points of $X^\an\setminus D$.
  \item Let $x\in X(K)$ and let $v$ be the unique tangent direction at $x$.
    Then $d_v F(x) = \ord_x(f)$.
  \end{enumerate}
\end{thm}

\paragraph[The skeleton as a metrized complex of curves] 
\label{par:complex.from.curves}
Let $(X,V\cup D)$ be a triangulated punctured curve with skeleton
$\Sigma = \Sigma(X, V\cup D)$.  Recall that $\Sigma$ is an
augmented $\Lambda$-metric graph with infinite vertices $D$.  We
enrich $\Sigma$ with the structure of a $\Lambda$-metrized complex of
$k$-curves as follows.  For $x\in V$ let $C_x$ be the smooth, proper, connected
$k$-curve with function field $\td\sH(x)$ as
in~\parref{par:semistable.decomposition}.  By definition $C_x$ has genus 
$g(x)$.  We have natural bijections
\begin{equation} \label{eq:tangent.DV.closedpt}
  T_x \cong \DV(\td\sH(x)/k) \cong C_x(k) 
\end{equation}
where the first bijection $v\mapsto\ord_v$ is defined
in~\parref{par:slope.formula} and the second associates to a closed point
$\xi\in C_x(k)$ the discrete valuation $\ord_\xi$ on the function field
$\td\sH(x)$ of $C_x$.  
Let $v\in T_x$ be a tangent direction and define
$\red_x(v)$ to be the point of $C_x$ corresponding to the discrete
valuation $\ord_v\in\DV(\td\sH(x)/k)$.  These data make
$\Sigma$ into a $\Lambda$-metrized complex of $k$-curves.

\paragraph[Lifting metrized complexes of curves]
We now prove that every metrized complex of curves over $k$ arises as the
skeleton of a smooth, proper, connected $K$-curve.
This fact appears in the literature in various contexts (over
discretely-valued fields): see for
instance~\cite[Appendix~B]{baker:specialization} 
and~\cite[Lemme~6.3]{Sai96}.  For this reason we only sketch a proof using
our methods.

\begin{thm} \label{thm:graph.to.curve}
  Let $\cC$ be a $\Lambda$-metrized complex of $k$-curves.  
  There exists a triangulated punctured curve $(X,V\cup D)$ such that
  the skeleton $\Sigma(X,V\cup D)$ is isomorphic to $\cC$.
\end{thm}

\pf
Let $C$ be a smooth, proper, connected $k$-curve.  By elementary
deformation theory, there is a smooth, proper admissible formal $R$-scheme
$\fC$ with special fiber $C$.  By GAGA the analytic generic fiber $\fC_K$
is the analytification of a smooth, proper, connected $K$-curve $\sC$.
There is a reduction map
$\red:\fC_K\to C$ from the analytic generic fiber of $\fC$ to (the set
underlying) $C$, under which the inverse image of the generic point of $C$
is a single distinguished point $x$.  The set $\{x\}$ is a semistable
vertex set of $\sC$, with associated skeleton also equal to $\{x\}$.
Moreover, there is 
a canonical identification of $\td\sH(x)$ with the field of rational
functions on $\bar C$
by~\cite[Proposition~2.4.4]{berkovich:analytic_geometry}. 
See~\parref{par:models.skeleta} for a more detailed discussion of the
relationship between semistable vertex sets and admissible formal models.

Let $\Gamma$ be the metric graph underlying $\cC$.
Let $V$ be the vertices of $\Gamma$ and let $E$ be its edges.
By adding valence-$2$ vertices we may assume that $\Gamma$ has no loop
edges.  Assume for the moment that $\Gamma$ has no infinite edges.
For a vertex $x\in V$ let $C_x$ denote the smooth, proper, connected
$k$-curve associated to $x$, and choose an admissible formal curve
$\fC_x$ with special fiber isomorphic to $C_x$ as above.  
For clarity we let $\red_{\fC_x}$ denote the reduction map
$(\fC_x)_K\to C_x$.
By~\cite[Proposition~2.2]{bl:fmI}, for every $\bar x\in C_x(k)$ the formal
fiber $\red_{\fC_x}^{-1}(\bar x)$ is isomorphic to $\bB(1)_+$.  Let $e$
be an edge of $\Gamma$ with endpoints $x,y$, and let 
$\bar x_e\in C_x(k),\bar y_e\in C_y(k)$
be the reductions of the tangent vectors in the direction of $e$ at $x,y$,
respectively.  Remove closed balls from 
$\red_{\fC_x}\inv(\bar x_e)$ and $\red_{\fC_y}\inv(\bar y_e)$ whose radii are such that
the remaining open annuli have modulus equal to $\ell(e)$.  We form a new
analytic curve $X^\an$ by gluing these annuli together using some
isomorphism of annuli
for each edge $e$.  The resulting curve is proper, hence is the
analytification of an algebraic curve $X$.  By construction the image of
the set of distinguished points of the curves $(\fC_x)_K$ in $X^\an$ is a
semistable vertex set, and the resulting skeleton (considered as a metric graph) is isomorphic to
$\Gamma$. 

If $\Gamma$ does have infinite edges, then we apply the above procedure to
the finite part of $\Gamma$, then puncture the curve $X$ in the formal
fibers over the smooth points of the residue curves which correspond to
the directions of the infinite tails of $\Gamma$.\qed

\medskip

As an immediate application, we obtain the following property of the
``abstract tropicalization map'' from the moduli space of stable marked
curves to the moduli space of stable abstract tropical curves:

\begin{cor} \label{cor:graph.to.curve}
  If $g$ and $n$ are nonnegative integers with $2-2g-n < 0$, the natural map
  ${\rm trop}: M_{g,n} \to M^{\rm trop}_{g,n}$ ({\rm
    see~\cite[Remark~5.52]{bpr:trop_curves}}) is surjective. 
\end{cor}

\smallskip
See~\cite{ACP12} for an interpretation of the above map as a contraction
from $M_{g,n}^\an$ onto its skeleton (which also implies
surjectivity).

\section{Morphisms between curves and their skeleta}
\label{sec:morphisms.skeleta}

This section is a relative version of the previous one, in that we propose
to study the behavior of semistable vertex sets and skeleta under finite
morphisms of curves.  We introduce finite morphisms of triangulated
punctured curves and we prove that any finite morphism of punctured curves
can be enriched to a finite morphism of triangulated punctured curves. 
This powerful result can be used to prove the simultaneous semistable
reduction theorems of Coleman, Liu--Lorenzini, and Liu, which we do in
Section~\ref{sec:simultaneous.ss.reduction}.  It is interesting to note that we
use only analytic methods on analytic $K$-curves, 
making (almost) no explicit reference to semistable models; hence our
proofs of the results of Liu--Lorenzini and Liu are very different from theirs.
Using a relative version of the slope formula, we also show that a finite morphism of triangulated
punctured curves induces (by restricting to skeleta) a finite harmonic morphism of metrized complexes
of curves, which will be crucial in the sequel.

\paragraph[Morphisms between open balls and generalized open annuli]
The main results of this section rest on a careful study of the behavior
of certain morphisms between open balls and generalized open annuli.
Some of the lemmas and intermediate results appear in some form in the
literature --- see in particular~\cite{berkovich:analytic_geometry}
and~\cite{berkovich:etalecohomology} --- although they are hard to find
in the form we need them.  As the proofs are not difficult, for the
convenience of the reader we include complete arguments.

\begin{lem} \label{lem:maps.annuli}
  Let $a,a'\in\fm_R$ and let $\phi:\bS(a')_+\to\bS(a)_+$ be a morphism.
  \begin{enumerate}
  \item If $\val\circ\phi$ is constant on $\Sigma(\bS(a')_+)$, then
    $\phi(\bS(a')_+)$ is contained in $\bS(a)_+\setminus\Sigma(\bS(a)_+)$.
  \item If $\val\circ\phi$ is not constant on $\Sigma(\bS(a')_+)$, then 
    $\phi(\Sigma(\bS(a')_+))\subset\Sigma(\bS(a))$ and the restriction of $\phi$ to
    $\Sigma(\bS(a')_+)$ has the form $x\mapsto m\cdot x + b$ for some nonzero
    integer $m$ and some $b\in\Lambda$.
  \item If $\phi$ extends to a continuous map 
    $\bar\phi:\bar{\bS(a')}_+\to\bar{\bS(a)}_+$
    taking an end of $\bS(a')_+$ to an end of $\bS(a)_+$, then
    $\val\circ\phi$ is not constant on $\Sigma(\bS(a')_+)$.
  \end{enumerate}
\end{lem}

\pf Part~(2) is exactly~\cite[Proposition~5.5]{bpr:trop_curves}, so
suppose that $\val\circ\phi$ is constant on $\Sigma(\bS(a))$.  
The map $\phi:\bS(a')_+\to\bS(a)_+\subset\G_m^\an$ is given by a
unit $f$ on $\bS(a')_+$.  By~\cite[Proposition~5.2]{bpr:trop_curves}, we can
write $f = \alpha\,(1+g)$, where $\alpha\in K^\times$ and $|g(x')|<1$ for
all $x'\in \bS(a')_+$.  Hence $|\phi(x')-\alpha| < |\alpha|$ for all 
$x'\in \bS(a')_+$; here $|\phi(x')-\alpha|$ should be interpreted as the absolute
value of the function $f-\alpha$ in $\sH(x')$, or
equivalently as the absolute value of the function $t-\alpha$
in $\sH(\phi(x))$, where $t$ is a parameter on $\G_m$.  
But $|t(x)-\alpha|\geq|\alpha|$ for every $x\in\Sigma(\bS(a))$, so
$\phi(\bS(a')_+)\cap\Sigma(\bS(a)) = \emptyset$.

Let $x'$ be an end of $\bS(a')_+$, let $x$ be an end of $\bS(a)_+$, and
suppose that $\bar\phi(x') = x$.  Since $x'$ is a limit point of
$\Sigma(\bS(a')_+)$, there exists a sequence of points
$x_1',x_2',\ldots\in\Sigma(\bS(a')_+)$ converging to $x'$.  Then
$\phi(x_1'),\phi(x_2'),\ldots$ converge to $x$, so
$\val(\phi(x_1')),\val(\phi(x_2')),\ldots$ converges to $\val(x)$.  Since
each $\val(\phi(x_i'))\in(0,\val(a))$ but
$\val(x)\notin(0,\val(a))$, $\val\circ\phi$ is not constant on
$\Sigma(\bS(a')_+)$.\qed

\begin{lem} \label{lem:map.ball.annulus}
  Let $B'$ be an open ball, let $A$ be a generalized open
  annulus, and let $\phi: B'\to A$ be a morphism.  Then
  $\phi(B')\cap\Sigma(A)=\emptyset$, and $\phi$ extends to a continuous
  map $\bar B\p\to A$.
\end{lem}

\pf The assertions of the Lemma are clear if $\phi$ is constant, so assume
that $\phi$ is non-constant.
We identify $A$ with $\bS(a)_+$ for $a\in\fm_R$ and $B'$ with
$\B(1)_+$.  Any unit on $B'$ has constant absolute value, so
$\val\circ\phi$ is constant on $B'$.  Since $\B(1)_+\setminus\{0\}$ is a
generalized open annulus, and since the type-$1$ point $0$ maps to a
type-$1$ point of $A$, it follows from Lemma~\ref{lem:maps.annuli} that
$\phi(B')\cap\Sigma(A)=\emptyset$.  Let $B$ be the connected component of
$A\setminus\Sigma(A)$ containing $\phi(B')$.  Then the morphism 
$B'\to B\inject A$ extends to a continuous map 
$\bar B\p\to\bar B\inject A$.\qed

\begin{defn} \label{defn:emb.int.expansion}
  Let $T$ be a metric space and let $m\in\R_{>0}$.  A continuous injection
  $\phi:[a,b]\inject T$ is an \emph{embedding with expansion factor $m$}
  provided that $r\mapsto\phi(r/m):[ma,mb]\inject T$ is a geodesic
  segment.
  A continuous injection $\phi:[a,b]\inject T$ is 
  \emph{piecewise affine} provided that
  there exist $a = a_0 < a_1 < \cdots < a_r = b$  and 
  $m_1,\ldots,m_r\in\R_{>0}$ such that $\phi|_{[a_{i-1},a_i]}$ is an
  embedding with expansion factor $m_i$ for each $i=1,\ldots,r$.
\end{defn}

\begin{lem} \label{lem:inject.geodesic}
  Let $B,B'$ be open balls, let $x'$ be the end of $B$, and let 
  $\phi:B'\to B$ be a morphism with finite fibers.  Let
  $y'\in B'$, let $y = \phi(y')$, and let $x = \phi(x')$.  Then the
  restriction of $\phi$ to the geodesic segment $[x',y']$ is injective,
  and $\phi([x',y'])$ is equal to the geodesic segment $[x,y]$.
  If in addition there exists $N>0$ such that all fibers of $\phi$ have
  fewer than $N$ elements,  then $\phi|_{[x',y']}$ is piecewise affine.
\end{lem}

\pf First suppose that $y'$ does not have type $4$.
Choose isomorphisms $B'\cong\B(1)_+$ such that $0 > y'$ and
$B\cong\B(1)_+$ such that $\phi(0) = 0$.  
Since $y'\in[x',0]$ we may replace $[x',y']$ by the larger geodesic segment
$[x',0]$ to assume that $y' = 0$.  
Define $F:\Sigma(\bS(0)_+)\to\R$ 
by $F(x) = \val\circ\phi(x) = -\log|\phi(x)|$.
By~\cite[Proposition~5.10]{bpr:trop_curves}, 
$F$ is a piecewise affine function, and for any point
$z'\in\Sigma(\bS(0)_+)$ at which $F$ is differentiable, the derivative of
$F$ is equal to the number of zeros $q'$ of $\phi$ with 
$\val(q') > \val(z')$.  It follows that $F$ is monotonically increasing, so
$\phi$ is injective on $[x',0]$.  Let 
\[ Z = \{ \tau(q')~:~q'\in\bS(0)_+,\,\phi(q')=0 \}\subset
\Sigma(\bS(0)_+). \]
If $C\subset\Sigma(\bS(0)_+)\setminus Z$ is a connected component, then
$A' = \tau\inv(C)$ is a generalized open annulus mapping to $\bS(0)_+$.
By the above, $\val\circ\phi$ is not constant on 
$\Sigma(A') = \Sigma(\bS(0)_+)\cap A'$, so by Lemma~\ref{lem:maps.annuli},
$\phi(\Sigma(A'))\subset\Sigma(\bS(0)_+)$.  Since 
$\Sigma(\bS(0)_+)\setminus Z$ is dense in $\Sigma(\bS(0)_+)$, we have that
$[x',0] = \Sigma(\bS(0)_+)\cup\{0\}$ maps bijectively onto $[x,0]$.  It is
clear in this case that the restriction of $\phi$ to the closure of any
connected component of $\Sigma(\bS(0)_+)\setminus Z$ is an embedding with
integer expansion factor, so $\phi|_{[x',0]}$ is piecewise affine.  Note
that $\phi|_{[x',0]}$ changes expansion factor at most $\#\phi\inv(0)$
times. 

Now suppose that $y'$ has type $4$.  Suppose that $z',w'\in[x',y']$ are two
distinct points such that $\phi(z') = \phi(w')$.  Assume without loss of
generality that $z' < w'$.  Since $w'$ and $z'$ have the same image under
$\phi$, they both have the same type, so $w'\neq y'$ because $y'$ is the
only type-$4$ point in $[x',y']$.  Applying the above to the geodesic
$[x',w']$ gives a contradiction.  Therefore, $\phi$ is injective on
$[x',y']$, so $\phi([x',y']) = [x,y]$
since $B$ is uniquely path-connected (\cite[Corollary~1.14]{baker_rumely:book}).

Assume now that all fibers of $\phi$ have size bounded by $N$.  The
above argument proves that for all $z' < y'$, the restriction of $\phi$ to
$[x',z']$ is piecewise affine, and $\phi|_{[x',z']}$ changes expansion
factor at most $N$ times.  Therefore there exists $z_0' < y'$ and
$m\in\Z_{>0}$ such that for all $z'\in[z_0',y']$, the restriction of $\phi$
to $[z_0',z']$ is an embedding with expansion factor $m$.  It follows that 
$\phi|_{[z_0',y']}$ is an embedding with expansion factor $m$, so 
$\phi|_{[x',y']}$ is piecewise affine.\qed

\begin{lem} \label{lem:map.ball.ball}
  Let $B,B'$ be open balls and let $\phi: B'\to B$ be a morphism.  Suppose
  that $\phi$ is open, separated, and has finite fibers, so it extends
  canonically to a continuous map $\phi:\bar B\p\to\bar B$
  by Lemma~\ref{lem:extend.end.ball}.  Let
  $x'$ be the end of $B'$ and let $x = \phi(x')\in\bar B$.
  Then $B_1 = \phi(B')$ is an open ball connected component of 
  $\bar B\setminus\{x\}$ and $\phi:B'\to B_1$ is finite and order-preserving. 
  In particular, if $x$ is the end of $B$, then $\phi: B'\to B$ is finite
  and order-preserving.
\end{lem}

\pf Let $y'\in B'$.  By Lemma~\ref{lem:inject.geodesic}, the restriction
of $\phi$ to the geodesic $[x',y']$ is injective, so $\phi(y')\neq x$.
Therefore $x\notin B_1$.
Let $C$ be the connected component of $B\setminus\{x\}$ containing
$B_1$.  Since $x$ is an end of $C$ and $x'$ maps to $x$, 
Lemma~\ref{lem:map.ball.annulus} implies that $C$ is not an open annulus,
so it is an open ball.  By hypothesis, $B_1$ is an open subset of $C$;
since $B'\cup\{x'\}$ is compact, $\phi(B'\cup\{x'\}) = B_1\cup\{x\}$ is
closed in $\bar B$, and therefore $B_1 = (B_1\cup\{x\})\cap C$ is closed
in $C$.  Since $C$ is connected, we have $B_1 = C$.  

Since $\phi$ has finite fibers, $\phi:B'\to B_1$ is finite if and only if it
is proper by~\cite[Corollary~3.3.8]{berkovich:analytic_geometry}.  
If $D\subset B_1$ is compact, then $\phi\inv(D)$ is compact since 
$\phi\inv(D)$ is closed as a subset of the compact space $\bar B\p$.
One has $\Int(B'/B_1) = B'$
by~\cite[Proposition~3.1.3(i)]{berkovich:analytic_geometry}
since $B'$ is a boundaryless $K$-analytic space.  Therefore
$\phi$ is proper.

The fact that $\phi: B'\to B_1$ is order-preserving follows immediately
from Lemma~\ref{lem:inject.geodesic}.\qed

\paragraph 
Let $\phi: B'\to B$ be a finite morphism of open balls.  Since an open
ball is a smooth curve, $\phi$ is flat in the sense that if $\sM(\cA)\subset
B$ is an affinoid domain and $\phi\inv(\sM(\cA)) = \sM(\cA')$, then $\cA'$ is
a (finite) flat $\cA$-algebra.  For $x\in\sM(\cA)$ the fiber over $x$ is
\[ \phi\inv(x) = \sM\big(\cA'\hat\tensor_\cA \sH(x)\big) 
= \sM\big(\cA'\tensor_\cA\sH(x)\big), \]
where the second equality holds because $\cA'$ is a finite $\cA$-algebra.
It follows that for $x\in B$ the quantity
$\sum_{x'\mapsto x} \dim_{\sH(x)} \OO_{\phi\inv(x),x'}$ is independent of
$x$; we call this number the \emph{degree} of $\phi$.

\begin{prop} \label{prop:skel.inverse.ball}
  Let $B,B'$ be open balls and let $\phi: B'\to B$ be a finite morphism.
  Let $x'$ be the end of $B'$, let $x = \phi(x')$ be the end of $B$, let
  $y\in B$, and let
  $y_1',\ldots,y_n'\in B'$ be the inverse images of $y$ in $B'$.   Then 
  \[ \phi\inv([x,y]) = [x',y_1']\cup\cdots\cup[x',y_n']. \]
\end{prop}

\pf Let $T' = [x',y_1']\cup\cdots\cup[x',y_n']$.
The inclusion $T'\subset\phi\inv([x,y])$ follows from
Lemma~\ref{lem:inject.geodesic}.  First we claim that for $z\in[x,y]$ near
enough to $x$, there is only one preimage of $z$ in $B'$.  Shrinking
$[x,y]$ if necessary, we may assume that $y$ has type $2$.  
Choose a type-$1$ point $w\in B$ such that $w > y$, so 
$(x,y]\subset(x,w) = \Sigma(B\setminus\{w\})$.  Let $w'\in B'$
be a preimage of $w$ and choose $z'\in\Sigma(B'\setminus\{w'\})$ such that 
$z'\leq\tau(q')$ for all $q'\in\phi\inv(w)$, where 
$\tau:B'\setminus\{w'\}\to\Sigma(B'\setminus\{w'\})$ is the retraction.
Let $A'$ be the open annulus $\tau\inv((x',z'))\subset B'$.  Then 
$\phi(A')\subset B\setminus\{w\}$, and $\phi$ takes the end $x'$ of $A'$
to the end $x$ of $B\setminus\{w\}$, so
by Lemma~\ref{lem:maps.annuli}(2,3), 
$\phi(\Sigma(A'))\subset\Sigma(B\setminus\{w\})$ and the map
$\Sigma(A')\to\Sigma(B\setminus\{w\})$ is injective.
If $u'\in A'\setminus\Sigma(A')$, then $u'$ is contained in an open ball in
$A'$, so $\phi(u')\notin\Sigma(B\setminus\{w\})$ by Lemma~\ref{lem:map.ball.annulus}.
Therefore every point of $\Sigma(B\setminus\{w\})$ has at most one preimage in
$A'$.  Let $u'\in B'\setminus A'$, so $u'\geq z'$.  Then
$\phi(u')\geq\phi(z')$, so every point $u\in(x,y]$ with
$u <\phi(z')$ has exactly one preimage in $B'$ (note that
$\phi(z')\in\Sigma(B\setminus\{w\})$ by Lemma~\ref{lem:inject.geodesic}).

Let $d$ be the degree of the finite morphism $\phi$, so for every
$z\in B$ we have
\[ d = \sum_{z'\in\phi\inv(z)} \dim_{\sH(z)} \OO_{\phi\inv(z),z'}. \]
For $y\in B$ as in the statement of the Proposition,
define a function $\delta:(x,y]\to\Z$ by
\[ \delta(z) = \sum_{\substack{z'\in\phi\inv(z)\\z'\in T'}}
\dim_{\sH(z)} \OO_{\phi\inv(z),z'}. \]
Clearly $\delta(z)\leq d$ for all $z\in(x,y]$, and 
$\delta(z) = d$ if and only if $\phi\inv(z)\subset T'$.  By definition 
of $T'$ we have $\delta(y) = d$, and by the above, $\delta(z) = d$ for $z\in(x,y]$
close enough to $x$ (any geodesic segment $[x',y_i']$ surjects onto $[x,y]$, hence
contains the unique preimage of $z$).  Therefore it is enough to
show that $\delta(z_1)\geq\delta(z_2)$ if $z_1\leq z_2$.

If $z\in B$ is a point of type $2$ or $3$, then $B_z\coloneq\{w\in B~:~w\geq z\}$ is
a (not necessarily strict) affinoid subdomain of $B$.  If $w'\in B'$ is such that 
$\phi(w')\geq z$, then $\phi([x',w']) = [x,\phi(w')]$ is a geodesic containing $z$, so 
there exists $z'\in[x',w']$ mapping to $z$.  It follows that
\[ \phi\inv(B_z) = \Djunion_{z'\mapsto z} \{ w'\in B'~:~w'\geq z'\}
= \Djunion_{z'\mapsto z} B_{z'}' \]
is a disjoint union of affinoid domains.  Hence each map
$B_{z'}\to B_z$ is finite, and its degree is equal to 
$\dim_{\sH(z)} \OO_{\phi\inv(z),z'}$.  
Let $z_1,z_2\in(x,y]$ and assume that $z_1 < z_2$, so $z_1$ has type $2$
or $3$.  For $z_1'\in\phi\inv(z_1)$ we have 
\[ \dim_{\sH(z_1)} \OO_{\phi\inv(z_1),z_1'} = 
\sum_{\substack{z_2'\mapsto z_2\\z_2'\geq z_1'}}
\dim_{\sH(z_2)} \OO_{\phi\inv(z_2),z_2'}. \]
Summing over all $z_1'\in\phi\inv(z_1)\cap T'$, we obtain
\[\begin{split} \delta(z_1) 
&= \sum_{\substack{z_1'\mapsto z_1\\z_1'\in T'}} \dim_{\sH(z_1)}\OO_{\phi\inv(z_1),z_1'}
= \sum_{\substack{z_1'\mapsto z_1\\z_1'\in T'}} 
\sum_{\substack{z_2'\mapsto z_2\\z_2'\geq z_1'}} \dim_{\sH(z_2)}\OO_{\phi\inv(z_2),z_2'}\\
&\geq \sum_{\substack{z_1'\mapsto z_1\\z_1'\in T'}} 
\sum_{\substack{z_2'\mapsto z_2\\z_2'\geq z_1'\\z_2'\in T'}}\dim_{\sH(z_2)}\OO_{\phi\inv(z_2),z_2'}
= \delta(z_2),
\end{split}\]
where the final equality holds because if $z_2'\in\phi\inv(z_2)\cap T'$,
then there exists $z_1'\in\phi\inv(z_1)\cap T'$ such that $z_2'\geq z_1'$,
namely,  the unique point of $[x',z_2']$ mapping to $z_1$.\qed

\begin{prop} \label{prop:image.of.edge}
  Let $a\in\fm_R\setminus\{0\}$, let $A' = \bS(a)_+$, let $B$ be an open
  ball, and let $\phi: A'\to B$ be a morphism with finite fibers.  Suppose
  that each end of $A'$ maps to the end of $B$ or to a type-$2$ point of $B$ 
  under the induced map $\bar A\p\to\bar B$.  
  Let $\alpha = \phi\circ\sigma: [0,\val(a)]\to\bar B$.
  Then there exist finitely many numbers
  $r_0,r_1,r_2,\ldots,r_n\in\Lambda$ with
  $0 = r_0 < r_1 < r_2 < \cdots < r_n = \val(a)$ such that 
  $\alpha$ is an embedding with nonzero integer expansion factor when
  restricted to each interval $[r_i,r_{i+1}]$.
  In other words, $\alpha$ is piecewise affine with integer expansion
  factors.
  Moreover, the image of $\alpha$ is a geodesic segment between type-$2$
  points of $B$.
\end{prop}

\pf Let $x'$ be the end $\sigma(0)$ of $A'$ and let $x = \phi(x')$.
Suppose first that $x$ is the end of $B$.  Choose an identification
$B\cong\B(1)_+$, and let $r_+ = \min\{\val(y')~:~\phi(y')=0\} > 0$.
Let $A'_+$ be the open annulus $\val\inv((0,r_+))\subset A'$.  Then
$\phi(A'_+)\subset\B(1)_+\setminus\{0\}=\bS(0)_+$ and the end $x'$ of
$A'_+$ maps to the end $x$ of $\bS(0)_+$, so by
Lemma~\ref{lem:maps.annuli}(2,3), $\alpha$ is an embedding
with (nonzero) integer expansion factor when restricted to $[0,r_+]$.

Now suppose that $x$ is not the end of $B$.  Let 
$r = \min\{\val(y')~:~\phi(y') = x\} > 0$ and let 
$A'' = \val\inv((0,r))$.  Then $\phi(A'')$ is contained in a connected
component $C$ of $B\setminus\{x\}$, and $\phi$ takes the end $x'$ of $A''$
to the end $x$ of $C$.  If $C$ is an open annulus, then $\alpha$ is an
embedding with integer expansion factor when restricted to $[0,r]$
by Lemma~\ref{lem:maps.annuli}(2,3), and if $C$ is an open ball,
then we proceed as above to find $r_+\in(0,r)\cap\Lambda$ such that
$\alpha$ is an embedding with integer expansion factor when restricted to
$[0,r_+]$. 

Applying the above argument to the morphism $\phi$ composed with the
automorphism $t\mapsto a/t$ of $\bS(a)_+$ (which interchanges the two
ends), we find that there exists $r_-\in[0,\val(a))\cap\Lambda$ such that
$\alpha$ is an embedding with integer expansion factor when restricted to
$[r_-,\val(a)]$.  Let $s\in(0,\val(a))\cap\Lambda$.  Replacing
$A'$ with the annulus $\val\inv((0,s))$ (resp.\ $\val\inv((s,\val(a)))$),
the above arguments then provide us with $s_+\in(s,\val(a)]\cap\Lambda$ (resp.\
$s_-\in[0,s)\cap\Lambda$) such that 
$\alpha$ is an embedding with integer expansion factor when restricted to 
$[s,s_+]$ (resp.\ $[s_-,s]$).  The first assertions now follow because there is a
finite subcover of the open covering 
\[ \{ [0,r_+) \} \cup \{ (r_-,\val(a)] \} 
\cup \{ (s_-,s_+)~:~s\in(0,r)\cap\Lambda \} \]
of the compact space $[0,r]$.

As for the final assertion, choose
$0 = r_0 < r_1 < r_2 < \cdots < r_n = \val(a)$
such that $\alpha$ is an embedding with integer expansion factor on each
$[r_i,r_{i+1}]$.  Let $i_0\in\{0,1,\ldots,n\}$ be the largest integer
such that $\alpha(r_{i_0}) > \alpha(r_i)$ (in the canonical partial
ordering on $B$) for 
all $i < i_0$.  If $i_0 = n$, then we are done, so assume that $i_0<n$.
Let $y = \alpha(r_{i_0})\in B$, and choose an identification
$B\cong\B(1)_+$ such that $0 > y$.  If $F = -\log|\phi|$, then 
by~\cite[Proposition~5.10]{bpr:trop_curves}, at every point
$r\in(0,\val(a))$ the change in slope of $F$ at $r$ is equal to 
the negative of the number of zeros of $\phi$ with valuation $r$
(cf.\ the proof of Lemma~\ref{lem:inject.geodesic}); in particular, the
slope of $F$ can only \emph{decrease}.  By construction $F$
is monotonically increasing on $[0,r_{i_0}]$.  Since
$\alpha([r_{i_0},r_{i_0+1}])$ is a geodesic segment, it meets $y$ only at
$\alpha(r_{i_0})$.  The image of $(r_{i_0},r_{i_0+1}]$ under $\alpha$
is not contained in an open ball connected
component of $B\setminus\{y\}$ because $\alpha(r_{i_0+1})\not> y$; hence
$F$ is \emph{decreasing} on an interval $[r_{i_0},r_{i_0}+\epsilon]$ for some
$\epsilon>0$.  Since the slope of $F$ can only decrease, it follows that
$F$ is monotonically decreasing on $[r_{i_0},r_n]$.  It follows 
immediately from this that $\alpha([0,\val(a)]) = [x,y]$ or
$\alpha([0,\val(a)]) = [\alpha(\val(a)),y]$, whichever segment is larger.
\qed

\paragraph[Morphisms between curves and skeleta]
In what follows we fix smooth, connected, proper algebraic
$K$-curves $X,X'$ and a finite morphism $\phi: X'\to X$.  Let 
$D\subset X(K)$ and
$D'\subset X'(K)$ be finite sets of closed points.
The map on analytifications $\phi:X'^\an\to X^\an$ is finite and
open by~\cite[Lemma~3.2.4]{berkovich:analytic_geometry}.

\begin{prop} \label{prop:img.skeleton}
  Let $\Sigma'$ be a skeleton of $(X',D')$ and let $\Sigma$ be
  a skeleton of $(X,D)$.  There exists a skeleton
  $\Sigma_1$ of $(X,D\cup\phi(D'))$ containing
  $\Sigma\cup\phi(\Sigma')$, and there is a minimal such $\Sigma_1$
  with respect to inclusion.
\end{prop}

\pf First we will prove the Proposition in the case
$D = D' = \emptyset$.  Let $V$ be a vertex set for $\Sigma$ and
let $V'$ be a vertex set for $\Sigma'$ containing 
$\{\tau(y')~:~y'\in\phi\inv(V)\}$.  Let $A'$ be an open annulus connected
component of $X'^\an\setminus V'$ and let $e'\subset\Sigma'$ be the
associated edge.  We claim that $\phi(e')$ is a geodesic
segment between type-$2$ points of $X^\an$.  
Let $C$ be the connected component of $X^\an\setminus V$
containing $\phi(A')$.
\begin{enum}
\item If $C$ is an open ball then the claim follows immediately from
  Proposition~\ref{prop:image.of.edge}. 

\item If $C$ is an open annulus and $\phi(A')\cap\Sigma(C)=\emptyset$, then
  $A'$ is contained in an open ball in $C$ because each connected
  component of $C\setminus\Sigma(C)$ is an open ball, so the claim is
  true as in~(i).

\item If $C$ is an open annulus and $\phi(A')\cap\Sigma(C)\neq\emptyset$,
  then $\phi(e')$ is a geodesic segment in $X^\an$ by
  Lemma~\ref{lem:maps.annuli}. 
\end{enum}

Applying the above to each edge $e'$ of $\Sigma'$, we find that there
exists a finite set of type-$2$ points $x_1,y_1,x_2,y_2,\ldots,x_n,y_n\in
X^\an$ such that 
$\phi(\Sigma') = \bigcup_{i=1}^n [x_i,y_i]$, where $[x_i,y_i]$ denotes a
geodesic segment from $x_i$ to $y_i$.  Let
\[ \Sigma_1 = \Sigma \cup\bigcup_{i=1}^n 
\big([x_i,\,\tau_\Sigma(x_i)]\cup[y_i,\,\tau_\Sigma(y_i)]\big). \]
By Lemma~\ref{lem:larger.vertex.set}(2) as applied to
$W = \{ x_1,y_1,\ldots,x_n,y_n \}$, we have that $\Sigma_1$ is a skeleton
of $X$, so by Proposition~\ref{prop:skeleton.properties}(5), $\Sigma_1$
contains $\Sigma\cup\phi(\Sigma')$ since $\Sigma_1$ contains each geodesic
segment $[x_i,y_i]$.  Any other skeleton of $X$ containing
$\Sigma$ and all $x_i,y_i$ also contains the geodesics 
$[x_i,\tau_\Sigma(x_i)]$ and $[y_i,\tau_\Sigma(y_i)]$, so $\Sigma_1$ is
the minimal such skeleton.

Now assume that $D\djunion D'\neq\emptyset$.  By
Lemma~\ref{lem:larger.vertex.set}(3) there is a minimal skeleton
$\Sigma_1$ of $X$ containing $\Sigma$ and $\phi(D')$; replacing $\Sigma$
by $\Sigma_1$ and $D$ by $D\cup\phi(D')$, we may assume without loss of
generality that $\phi(D')\subset D$.
We will proceed by induction on the size of $D\djunion D'$.  
As above we let $V$ be a vertex set for $\Sigma$ and we
let $V'$ be a vertex set for $\Sigma'$ containing 
$\{\tau(y')~:~y'\in\phi\inv(V\cup D)\}$.
\begin{enum}
\item Suppose that $D'$ is not empty.  Let $x'\in D'$ and let
  $\Sigma_1$ be the minimal skeleton of $(X,D)$ containing
  $\Sigma\cup\phi(\Sigma(X',V'\cup D'\setminus\{x'\}))$.
  Let $A'$ be the connected component of $X'^\an\setminus(V'\cup D')$ whose
  closure contains $x'$, let $x = \phi(x')\in D$, and let $A$ be the
  connected component of $X^\an\setminus(V\cup D)$ whose closure contains $x$.
  Note that 
  $\Sigma' = \Sigma(X',V'\cup D'\setminus\{x'\})\cup\Sigma(A')\cup\{x'\}$.
  Clearly $\phi(A')\subset A$, and $\phi\inv(x)\cap A' = \emptyset$ by
  construction. Since $\phi$ takes the puncture of $A'$ to the puncture of
  $A$, one shows as in the proof of Lemma~\ref{lem:maps.annuli} that
  $\phi(\Sigma(A'))\subset\Sigma(A)$.  Therefore $\Sigma_1$ 
  contains $\Sigma\cup\phi(\Sigma')$.

\item Now suppose that $D' = \emptyset$ and $D\neq\emptyset$.  Let
  $x\in D$ and let $\Sigma_1$ be the minimal skeleton of
  $(X,D\setminus\{x\})$ containing  
  $\Sigma(X,V\cup D\setminus\{x\})\cup\phi(\Sigma')$.
  Let $\Sigma_2 = \Sigma_1\cup[x,\tau(x)]$.  It follows from
  Lemma~\ref{lem:larger.vertex.set}(3) that $\Sigma_2$ is the minimal
  skeleton of $(X,D)$ containing $\Sigma_1$.
\end{enum}\qed

Let $V$ be a semistable vertex set of $(X,D)$.
Recall from~\parref{par:semistable.decomposition} that a connected
component $C$ of $X^\an\setminus(V\cup D)$ is \emph{adjacent} to a vertex
$x\in V$ provided that the closure of $C$ in $X^\an$ contains $x$.
Proposition~\ref{prop:stable.skel.inverse}(1) below is 
exactly~\cite[Theorem~4.5.3]{berkovich:analytic_geometry}
when $D = D' = \emptyset$.

\begin{prop} \label{prop:stable.skel.inverse}
  Suppose that $D' = \phi\inv(D)$ and that one of the following two
  conditions holds:
  \begin{enumerate}
  \item $\chi(X\setminus D)\leq 0$ (hence also
    $\chi(X'\setminus D') \leq 0$), or
  \item $\phi\inv(V)\subset V'$.
  \end{enumerate}
  In the situation of~(1) let $\Sigma$ (resp.\ $\Sigma'$) be the minimal
  skeleton of $(X,D)$ (resp.\ $(X',D')$), and in~(2) let 
  $\Sigma = \Sigma(X,V\cup D)$ and $\Sigma' = \Sigma(X',V'\cup D')$.
  Then $\phi\inv(\Sigma)\subset\Sigma'$.
\end{prop}

\pf First suppose that~(2) holds.
Let $x'\in X'^\an\setminus\Sigma'$ and
let $B'$ be the connected component of $X'^\an\setminus\Sigma'$ containing
$x'$, so $B'$ is an open ball.  By hypothesis, $\phi(B')$ is contained in a connected
component $C$ of $X^\an\setminus(V\cup D)$.  If $C$ is an open ball, then
$\phi(B')\cap\Sigma = \emptyset$ by the definition of $\Sigma(X,V\cup D)$.  If
$C$ is a generalized open annulus, then 
$\phi(B')\cap\Sigma=\phi(B')\cap\Sigma(C) = \emptyset$ by 
Lemma~\ref{lem:map.ball.annulus}.  Therefore $\phi(x')\notin\Sigma$.

Now suppose that~(1) holds.
Let $V$ be a semistable vertex set of $(X,D)$ such that
$\Sigma=\Sigma(X,V\cup D)$.  By subdividing edges of $\Sigma$ and enlarging $V$
if necessary, we may and do assume that $\Sigma$ has no loop edges.
By Proposition~\ref{prop:minimal.skeleton}(2), no point in $V$ of genus 
zero has valence one in $\Sigma$.
First we claim that if $B'\subset X'^\an\setminus D'$
is an open analytic domain which is isomorphic to $\B(1)_+$, then 
$\phi(B')\cap V = \emptyset$.  If $\phi\inv(V)\cap B'$ contains more than
one point, then it easy to see that there exists a smaller open ball
$B''\subset B'$ such that $\phi\inv(V)\cap B''$ contains exactly one
point.  Replacing $B'$ by $B''$, we may assume that there is a unique
point $y'\in\phi\inv(V)\cap B'$.  Let $y = \phi(y')\in V$.  Since
$g(y') = 0$ we have $g(y) = 0$.  

The open analytic domain $B'\setminus\{y'\}$ is the disjoint union of an open
annulus $A'$ and an infinite collection of open balls by
Lemma~\ref{lem:ball.subtract.point}.  Each connected 
component $C'$ of $B'\setminus\{y'\}$ maps into a connected component $C$
of $X^\an\setminus(V\cup D)$ adjacent to $y$, with the end $y'$ of $C'$ mapping to
the end $y$ of $C$.  By Lemma~\ref{lem:map.ball.annulus},
no open ball connected component of $B'\setminus\{y'\}$ can map to a generalized open
annulus connected component of $X^\an\setminus(V\cup D)$.  There are at least
two generalized open annulus connected components $A$ of
$X^\an\setminus(V\cup D)$ adjacent to $y$,
so some such $A$ must satisfy $\phi(B')\cap A = \emptyset$.  
But the map $\phi:X'^\an\to X^\an$ is \emph{open} 
by~\cite[Lemma~3.2.4]{berkovich:analytic_geometry}, so $\phi(B')$ is an
open neighborhood of $y$, which contradicts the fact that
$y$ is a limit point of $A$.  This proves the claim.

Let $V'$ be a semistable vertex set of $(X',D')$ such that
$\Sigma' = \Sigma(X',V'\cup D')$.   Since $X'\setminus\Sigma'$ is a
disjoint union of open balls, by the above we have
$\phi\inv(V)\subset\Sigma'$.  Hence we may enlarge $V'$ to contain
$\phi\inv(V)$ without changing $\Sigma'$, so we are reduced to~(2).
\qed

\begin{thm} \label{thm:inv.skeleton}
  Let $\Sigma'$ be a skeleton of $(X',D')$ and let $\Sigma$ be a
  skeleton of $(X,D)$.  Suppose that
  $\phi(\Sigma')\subset\Sigma$, and if
  $X\cong\P^1$ assume in addition that there exists a type-$2$ point
  $z\in\Sigma$ such that $\phi\inv(z)\subset\Sigma'$.
  Then $\phi\inv(\Sigma)$ is a skeleton of $(X',\phi\inv(D))$.
\end{thm}

\smallskip
We will require the following lemmas.
Lemma~\ref{lem:ball.image} is similar
to~\cite[Corollary~4.5.4]{berkovich:analytic_geometry}.

\begin{lem} \label{lem:ball.image}
  Let $B'\subset X'^\an$ be an open analytic domain isomorphic to $\B(1)_+$
  and let $B = \phi(B')\subset X^\an$.  Then $B$ is an open analytic domain of
  $X^\an$, and one of the following is true:
  \begin{enumerate}
  \item $B$ is an open ball and $\phi: B'\to B$ is finite and order-preserving. 
  \item $X \cong\P^1$ and $B = X^\an$.
  \end{enumerate}
\end{lem}

\pf Suppose that the genus of $X$ is at least one, so $X'$ also has genus
at least one. 
By Lemma~\ref{lem:map.ball.ball} we only need to show that
$B$ is contained in an open ball in $X^\an$.  Let $\Sigma$
(resp.\ $\Sigma'$) be the minimal skeleton of $X$ (resp.\ $X'$).
By Proposition~\ref{prop:minimal.skeleton}(1) we have 
$B'\subset X'^\an\setminus\Sigma'$, so by
Proposition~\ref{prop:stable.skel.inverse}(1) we have 
$B\subset X^\an\setminus\Sigma$.  But
every connected component of $X^\an\setminus\Sigma$ is an open ball, so
$B$ is contained in an open ball.

In the case $X = \P^1$,
suppose first that $B'(K)\to\P^1(K)$ is not surjective, so we may
assume that $\infty\notin B$ after choosing a suitable coordinate on $\P^1$.  Let
$x'$ be the end of $B'$ and let $x = \phi(x')$.  Then $x$ is a type-$2$
point, so $\phi(B'\cup\{x'\}) = B\cup\{x\}$ is a compact subset of
$\A^{1,\an}$.  Since $\A^{1,\an}$ is covered by an increasing union of open
balls, $B$ is contained in an open ball.

Now suppose that $B'(K)\to\P^1(K)$ is surjective.  Since $\P^1(K)$ is
dense in $\P^{1,\an}$ and $B\cup\{x\}$ is closed in $\P^{1,\an}$, it
follows that $B\cup\{x\} = \P^{1,\an}$.  If $x\notin B$, then $B$ is
contained in a connected component of $\P^{1,\an}\setminus\{x\}$, which
contradicts the surjectivity of $B'(K)\to\P^1(K)$.  Therefore
$B=\P^{1,\an}$.\qed

\begin{rem}
  Case~(2) of Lemma~\ref{lem:ball.image} does occur.  For instance, let
  $\phi:\P^1\to\P^1$ be the finite morphism $t\mapsto t^2$, and let
  $B'$ be the open ball in $\P^{1,\an}$ obtained by deleting the closed
  ball of radius $1/2$ around $1\in\P^1(K)$.  If $\chr(k)\neq 2$, then for
  every point $x\in\P^1(K)$, either $x\in B'(K)$ or $-x\in B'(K)$, so 
  $B'\to\P^{1,\an}$ is surjective.
\end{rem}

\begin{lem} \label{lem:ball.skeleton}
  Let $V$ be a semistable vertex set of $(X,D)$, let
  $\Sigma=\Sigma(X,V\cup D)$, and let $B\subset X$ be an open analytic domain 
  isomorphic to $\B(1)_+$ whose end $x$ is contained in $\Sigma$.  Then 
  \[ \Sigma\cap\bar B = \bigcup_{y\in B\cap(V\cup D)} [x,y]. \]
\end{lem}

\pf We may assume without loss of generality that $x\in V$.  Let
$V_0 = V\setminus(V\cap B)$ and $D_0 = D\setminus(D\cap B)$.  Then $V_0$
is a semistable vertex set of $(X,D_0)$: indeed, every connected
component of $X^\an\setminus(V\cup D)$ is either a connected
component of $X^\an\setminus(V_0\cup D_0)$ or is contained in $B$, so 
$X^\an\setminus(V_0\cup D_0)$ is a disjoint union of open balls and
finitely many generalized open annuli.  Let 
$\Sigma_0 = \Sigma(X,V_0\cup D_0)$, so $x\in\Sigma_0$ but
$B\cap\Sigma_0=\emptyset$ by
Proposition~\ref{prop:skeleton.properties}(1).  By
Lemma~\ref{lem:map.ball.ball} we have that $B$ is a
connected component of $X^\an\setminus\Sigma_0$, so
the Lemma now follows from Lemma~\ref{lem:larger.vertex.set} as applied to 
$\Sigma = \Sigma_0$, $W = V\cap B$, and $E = D\cap B$.\qed

\pf[of Theorem~\ref{thm:inv.skeleton}]
Note that $\phi(\Sigma')\subset\Sigma$ implies $\phi(D')\subset D$.
Let $V$ (resp.\ $V'$) be a vertex set for $\Sigma$ (resp.\ $\Sigma'$).  We
claim that 
\[ \phi\inv(\Sigma) = \Sigma'\cup\bigcup_{x'\in\phi\inv(V\cup D)} [x',\,\tau_{\Sigma'}(x')], \]
which by Lemma~\ref{lem:larger.vertex.set} is the minimal skeleton
of $(X',\phi\inv(D))$ containing $\Sigma'$ and 
$\phi\inv(V)$. 

Let $B'$ be a connected component of $X'^\an\setminus\Sigma'$, let
$x'\in\Sigma'$ be its end, and let $x = \phi(x')\in\Sigma$.  
Then $B = \phi(B')$ is an open ball and $x$ is its end: if
$X\not\cong\P^1$, then this follows 
directly from Lemma~\ref{lem:ball.image}, and if $X\cong\P^1$, then 
$B\subset X^\an\setminus\{z\}$, so $\phi(B)\neq X^\an$ and therefore $B$
is an open ball in this case as well.  Hence
\[ \Sigma\cap B = \bigcup_{y\in B\cap(V\cup D)} [x,y] \]
by Lemma~\ref{lem:ball.skeleton}, so 
\[ \phi\inv(\Sigma)\cap B' = \bigcup_{y'\in B'\cap\phi\inv(V\cup D)} [x',y'] \]
by Proposition~\ref{prop:skel.inverse.ball}.\qed

\begin{rem}
  When $X\cong\P^1$, the extra hypothesis on $\Sigma'$ in
  Theorem~\ref{thm:inv.skeleton} is necessary.  Indeed, let 
  $\phi:\P^1\to\P^1$ be a finite morphism and let $x'\in\P^{1,\an}$ be a
  type-$2$ point such that $\phi\inv(\phi(x'))$ has more than one
  element.  Then $\Sigma = \{\phi(x')\}$ is a skeleton containing the
  image of the skeleton $\Sigma'=\{x'\}$, but $\phi\inv(\Sigma)$ is not a skeleton.
\end{rem}

\begin{cor} \label{cor:stable.hull.skel}
  Let $\Sigma$ (resp.\ $\Sigma'$) be a skeleton of $(X,D)$ 
  (resp.\ $(X',D')$).  There exists a skeleton
  $\Sigma_1$ of $(X,D\cup\phi(D'))$ such that
  $\Sigma_1\supset\Sigma\cup\phi(\Sigma')$ and such that
  $\phi\inv(\Sigma_1)$ is a skeleton of $(X',\phi\inv(D\cup\phi(D')))$.
  Moreover, there is a minimal such 
  $\Sigma_1$ with respect to inclusion. 
\end{cor}

\pf If $X\not\cong\P^1$ then this is an immediate consequence of
Proposition~\ref{prop:img.skeleton} and Theorem~\ref{thm:inv.skeleton}:
if $\Sigma_1$ is the minimal skeleton of $(X,D\cup\phi(D'))$ containing
$\Sigma\cup\phi(\Sigma')$, then 
$\phi\inv(\Sigma_1)$ is a skeleton of $(X',\phi\inv(D\cup\phi(D')))$.
Suppose then that $X\cong\P^1$.  Let $x'\in\Sigma'$ be a type-$2$ point, let
$x=\phi(x')$, and let $x',x_1',x_2',\ldots,x_n'$ be the points of $X'$
mapping to $x$.  Let 
\[ \Sigma_1' = \Sigma' \cup \bigcup_{i=1}^n
[x_i',\,\tau_{\Sigma'}(x_i')]. \]
This is the minimal skeleton of $(X',D')$ containing
$\Sigma'$ and $\phi\inv(x)$ by Lemma~\ref{lem:larger.vertex.set}(3).
Let $\Sigma_1$ be the minimal skeleton of $(X,D\cup\phi(D'))$ containing 
$\Sigma\cup\phi(\Sigma_1')$.  Then $\phi\inv(\Sigma_1)$ is a skeleton of
$(X',\phi\inv(D\cup\phi(D')))$ by
Theorem~\ref{thm:inv.skeleton}.  If $\Sigma_2$ is a skeleton of
$(X,D\cup\phi(D'))$ such that $\Sigma_2\supset\Sigma\cup\phi(\Sigma')$ and  
$\phi\inv(\Sigma_2)$ is a skeleton of $(X',\phi\inv(D\cup\phi(D')))$, 
then $\phi\inv(\Sigma_2)$ contains each geodesic
$[x_i',\tau_{\Sigma'}(x_i')]$, so $\phi\inv(\Sigma_2)$ contains $\Sigma_1'$ and
therefore $\Sigma_2\supset\Sigma_1$.\qed

\begin{rem} \label{rem:structure.of.image.skeleton}
  Suppose that $\phi:X'\to X$ is a finite morphism such that
  $\phi\inv(D)=D'$ and $\phi\inv(V)=V'$.  Let
  $\Sigma = \Sigma(X,V\cup D)$ and $\Sigma'=\Sigma(X',V'\cup D')$.  
  In this case we can describe the skeleton $\Sigma_1$ of
  Corollary~\ref{cor:stable.hull.skel} more explicitly, as follows.  Let
  $e'$ be an open edge of $\Sigma'$ with respect to the given choice of
  vertex set and let $A' = \tau\inv(e')$.  This is a connected component
  of $X'^\an\setminus(V'\cup D')$.  Let $A$ be the connected component of
  $X^\an\setminus(V\cup D)$ containing $\phi(A')$.  Since $A'$ is a connected
  component of $\phi\inv(A)$, the map $A'\to A$ is finite, hence
  surjective.  If $A'$ is a punctured open ball, then so is $A$, so by
  Lemma~\ref{lem:maps.annuli}, $\phi$ maps $e'$ homeomorphically with 
  nonzero integer expansion factor onto $e = \Sigma(A)$.  If $A'$ is an
  open annulus, then $A$ could be an open annulus or an open ball; in the
  former case we again have $e'$ mapping homeomorphically onto
  $e = \Sigma(A)$ with nonzero integer expansion factor.  If $A$ is an
  open ball, then by Proposition~\ref{prop:image.of.edge}, $e'$ maps onto a
  geodesic segment in $A$; clearly both vertices of $e'$ map to the end of
  $A$, so $e'$ ``goes straight and doubles back'' under $\phi$.

  It follows from this and Proposition~\ref{prop:stable.skel.inverse}(2)
  that $\phi(\Sigma')$ is equal to $\Sigma$ union a finite
  number of geodesic segments $T_1,\ldots,T_r$ contained in open balls
  $B_1,\ldots,B_r\subset X^\an\setminus(V\cup D)$ and with one endpoint at
  elements of $V$.  These are precisely the images of the edges of
  $\Sigma'$ not mapping to edges of $\Sigma$.  By
  Lemma~\ref{lem:larger.vertex.set},  
  $\Sigma\cup(T_1\cup\cdots\cup T_r)$ is again a skeleton of $X$, so
  in this case 
  \[ \Sigma_1 = \phi(\Sigma') = \Sigma\cup(T_1\cup\cdots\cup T_r) \]
  (the inverse image of $\Sigma_1$ is a skeleton by
  Theorem~\ref{thm:inv.skeleton}). 
\end{rem}

\paragraph[Tangent directions and morphisms of curves]
Our next goal is to prove that a finite morphism of curves induces a
finite harmonic morphism of metrized complexes of curves for a suitable
choice of skeleta. First we need to formulate and prove a relative version
of Theorem~\ref{thm:slope.formula}. 
Let $X,X'$ be smooth, proper, connected $K$-curves.

\begin{defn}
  A continuous function $\phi: X'^\an\to X^\an$ is \emph{piecewise affine}
  provided that, for every geodesic segment $\alpha:[a,b]\inject X'^\an$,
  the composition $\phi\circ\alpha$ is piecewise affine with integer
  expansion factors in the sense of
  Definition~\ref{defn:emb.int.expansion}. 
\end{defn}

\smallskip
Let $\phi:X'^\an\to X^\an$ be a piecewise affine function, let 
$x'\in X'^\an$, and let $x = \phi(x')\in X^\an$. 
Let $v'\in T_{x'}$ and let $\alpha:[a,b]\inject X'^\an$ be a geodesic
segment representing $v'$ (so $a = -\infty$ if $x'$ has type $1$).  Taking
$b$ small enough, we can assume that $\phi\circ\alpha$ is an embedding
with integer expansion factor $m$.  Let $v\in T_x$ be the tangent
direction represented by $\phi\circ\alpha$.  We define
\[ d\phi(x') : T_{x'}\to T_x \sptxt{by} d\phi(x')(v') = v; \]
this is independent of the choice of $\alpha$.  We call the expansion
factor $m$ the \emph{outgoing slope of $\phi$ in the direction $v'$} and we write
$m = d_{v'}\phi(x')$.  

\begin{defn}
  A piecewise affine function $\phi: X'^\an\to X^\an$ is \emph{harmonic}
  at a point $x'\in X'^\an$ provided that, for all $v\in T_{\phi(x')}$, the
  integer
  \[ \sum_{\substack{v'\in T_{x'}\\d\phi(v')=v}} d_{v'}\phi(x') \]
  is independent of the choice of $v$.
\end{defn}

\smallskip
Let $\phi: X'\to X$ be a finite morphism of smooth, proper, connected
$K$-curves, let $x'\in X'^\an$ be a type-$2$ point, and let $x = \phi(x')$.
Let $C_x$ and $C_{x'}$ denote the smooth proper connected $k$-curves with
function fields $\td\sH(x)$ and $\td\sH(x')$, respectively.
We denote by $\phi_{x'}$ the induced morphism $C_{x'}\to C_x$.  We have
the following relative version of the slope formula of
Theorem~\ref{thm:slope.formula}: 

\begin{thm} \label{thm:harmonic.skeleta}
  Let $\phi: X'\to X$ be a finite morphism of smooth, proper, connected
  $K$-curves.
  \begin{enumerate}
  \item The analytification $\phi: X'^\an\to X^\an$ is piecewise affine and
    harmonic. 
  \item Let $x'\in X'^\an$ be a type-$2$ point, let $x = \phi(x')$, let 
    $v'\in T_{x'}$, let $v = d\phi(v')\in T_{x}$, and let 
    $\xi_{v'}\in C_{x'}$ and $\xi_{v}\in C_{x}$ be the closed points
    associated~\eqref{eq:tangent.DV.closedpt} to $v'$ and $v$,
    respectively.
    Then $\phi_{x'}(\xi_{v'}) = \xi_{v}$, and the ramification degree of
    $\phi_{x'}$ at $\xi_{v'}$ is equal to $d_{v'}\phi(x')$.
  \item Let $x'\in X'^\an$ be a type-$1$ point and let $v'\in T_{x'}$ be the
    unique tangent direction.  Then $d_{v'}\phi(x')$ is the ramification
    degree of $\phi$ at $x'$.
  \end{enumerate}
\end{thm}

\pf First we claim that $\phi$ is piecewise affine.
Let $[x',y']\subset X'^\an$ be a geodesic segment.
Suppose first that $x'$ and $y'$ have type $2$ or $3$.  Then there exists
a skeleton $\Sigma'$ of $X'$ containing $[x',y']$
by~\cite[Corollaries~5.56 and~5.64]{bpr:trop_curves}.  From this we easily
reduce to the case that $[x',y']$ is an edge of $\Sigma'$ with respect to
some vertex set; now the claim follows from Lemma~\ref{lem:maps.annuli}
and Proposition~\ref{prop:image.of.edge} as in the proof of
Proposition~\ref{prop:img.skeleton}.  If $x'$ has type $1$ or $4$, then
there is an open neighborhood $B$ of $x = \phi(x')$ and an open
neighborhood $B'\subset\phi\inv(B)$ of $x'$ such that $B$ and $B'$ are
open balls.  Shrinking $B'$ if necessary we can assume that
$y'\notin B'$, so the end $z'$ of $B'$ is contained in $[x',y']$ (removing
$z'$ disconnects $X'^\an$).
The restriction of $\phi$ to $[x',z']$ is piecewise affine by
Lemma~\ref{lem:inject.geodesic}, so this proves the claim.

Next we prove~(2). By functoriality of the reduction map, for any non-zero rational function $f$ on $X$, 
we have $\phi_{x'}^*(\td f_x) = \td{\phi^*(f)}_{x'}$.
By the slope formula (Theorem~\ref{thm:slope.formula}), 
the point $\xi_v$ corresponds to the discrete valuation on $\td \sH(x)$ given by 
$\ord_{\xi_v}(\td f_x) = d_v F(x)$, where $F = -\log|f|$.
Similarly, the slope formula applied to 
$\phi^*(f)$ on $X'$ gives 
$\ord_{\xi_{v'}}(\td{\phi^*(f)}_{x'}) =  d_{v'} \phi(x') d_v F(x)$. Therefore, we get 
\[\ord_{\xi_{v'}}(\phi_{x'}^*(\td f_x)) = d_{v'} \phi(x') \ord_{\xi_v}(\td f_x).\]
Since any non-zero element of $\td \sH(x)$ is of the form
$\td f_x$ for some non-zero rational function $f$ on $X$, this shows that the center of the valuation 
$\xi_{v'}$ on $\td \sH(x) \stackrel{\phi_{x'}^*}{\lhook\joinrel\relbar\joinrel\rightarrow} \td\sH(x')$ coincides with the center of 
$\xi_v$, that is $\phi_{x'}(\xi_{v'})=\xi_v$, and at the same time,  the ramification degree of $\phi_{x'}$ at $\xi_{v'}$ is equal to 
$d_{v'} \phi(x')$.

The proof of~(3) proceeds in the same way as~(2), applying the slope formula to a non-zero rational
function $f$  on $X$ with a zero at $x = \phi(x')$, and to $\phi^*(f)$ on $X'$.\qed

\paragraph[Morphisms of curves induce morphisms of metrized complexes]\label{par:incuced-morphisms-mc}
More precisely, this occurs when the morphism of curves respects a choice
of triangulations, in the following sense.

\begin{defn} \label{defn:morphism.triang.curves}
  Let $(X,V\cup D)$ and $(X', V'\cup D')$ be triangulated punctured
  curves.  A \emph{finite morphism} 
  $\phi: (X',V'\cup D')\to (X,V\cup D)$ is a finite morphism
  $\phi: X'\to X$ such that $\phi\inv(D) = D'$, $\phi\inv(V) = V'$, and 
  $\phi\inv(\Sigma(X,V\cup D)) = \Sigma(X',V'\cup D')$ (as sets).
\end{defn}

\begin{cor} \label{cor:skel.compat.exists}
  Let $\phi:X'\to X$ be a finite morphism of smooth, connected, proper
  $K$-curves, let $D\subset X(K)$ be a finite set, and let 
  $D' = \phi\inv(D)$.  There exists a strongly semistable vertex set
  $V$ of $(X,D)$ such that $V' = \phi\inv(V)$ is a strongly semistable
  vertex set of $(X',D')$ and such that
  \[ \Sigma(X',V'\cup D') = \phi\inv(\Sigma(X,V\cup D)). \]
  In particular, $\phi$ extends to a finite morphism of triangulated punctured
  curves. 
\end{cor}

\pf By Corollary~\ref{cor:stable.hull.skel}, there is a skeleton 
$\Sigma$ of $(X,D)$ such that $\Sigma' = \phi\inv(\Sigma)$ is a
skeleton of $(X',D')$.  Let $V_0'$ be a strongly semistable vertex set of
$\Sigma'$, let $V$ be a strongly semistable vertex set of $\Sigma$
containing $\phi(V_0')$, and let $V' = \phi\inv(V)$.  Then $V'$ is again a
strongly semistable vertex set of $\Sigma'$.\qed

\paragraph
Let $\phi: (X',V'\cup D')\to(X,V\cup D)$ be a finite morphism of
triangulated punctured curves.
Let $\Sigma = \Sigma(X,V\cup D)$ and $\Sigma' = \Sigma(X',V'\cup D')$.  
Let $e'$ be an open edge of $\Sigma'$ and let $e = \phi(e')$, an open edge
of $\Sigma$.  The image of the annulus $A' = \tau\inv(e')$ is contained in
$A = \tau\inv(e)$, so the restriction of $\phi$ to $e'$ is an embedding
with integer expansion factor by Lemma~\ref{lem:maps.annuli}(2).
It is clear from this and Theorem~\ref{thm:harmonic.skeleta} that 
$\phi|_{\Sigma'}: \Sigma'\to\Sigma$ is a harmonic $(V',V)$-morphism of
$\Lambda$-metric graphs.  For $x'\in V'$ with image $x = \phi(x')$ we have
a finite morphism of $k$-curves $\phi_{x'}: C_{x'}\to C_x$.  
It now follows from Theorem~\ref{thm:harmonic.skeleta} that 
these extra data enrich $\phi|_{\Sigma'}: \Sigma'\to\Sigma$ with the
structure of a finite harmonic morphism of $\Lambda$-metrized complexes of
curves: 

\begin{cor}\label{cor:morphism.to.harmonic}
  Let $\phi:(X',V'\cup D')\to(X,V\cup D)$ be a finite morphism of
  triangulated punctured curves.  Then $\phi$ naturally
  induces a finite harmonic morphism of $\Lambda$-metrized complexes of
  curves 
  \[ \Sigma(X',\, V'\cup D') \To \Sigma(X,\, V\cup D). \]
\end{cor}

\begin{rem} \label{rem:induced.harmonic.morphism.1}
More generally, suppose that $(X,V\cup D)$ and $(X', V'\cup D')$ are
triangulated punctured $K$-curves 
and that $\phi: X'\to X$ is a finite morphism such that 
$\phi\inv(D) = D'$ and $\phi\inv(V) \subset V'$.
Then $\phi\inv(\Sigma) \subset \Sigma'$ by
Proposition~\ref{prop:stable.skel.inverse}(2).
One can show using Theorem~\ref{thm:harmonic.skeleta} that the map
$\Sigma'\to X^\an$ followed by the retraction $X^\an\to\Sigma$
is a (not necessarily finite) harmonic morphism of metric graphs
$\Sigma'\to\Sigma$.
\end{rem}

\paragraph[Tame coverings of triangulated curves] 
\label{par:tame.covers.curves} 
Many of the results in this paper involve constructing a curve $X'$
and a morphism $\phi: X'\to X$ inducing a given morphism of skeleta as in 
Corollary~\ref{cor:morphism.to.harmonic}.  For these purposes it is useful
to introduce some mild restrictions on the morphism $\phi$.
Fix a triangulated punctured curve $(X,V\cup D)$ with skeleton
$\Sigma = \Sigma(X,V\cup D)$, regarded as a metrized complex of curves as
in~\parref{par:complex.from.curves}.   

\begin{defn} 
  Let $(X',V'\cup D')$ be a triangulated punctured curve with skeleton
  $\Sigma' = \Sigma(X',V'\cup D')$ and let
  $\phi:(X',V'\cup D')\to(X,V\cup D)$ be a finite morphism.  We say that
  $\phi$ is a \emph{tame covering of $(X,V\cup D)$} provided that:
  \begin{enumerate}
  \item $D$ contains the branch locus of $\phi$, 
  \item if $\chr(k) = p >0$, then
    for every edge $e'\in E(\Sigma')$ the expansion factor
    $d_{e'}(\phi)$ is not divisible by $p$, and
  \item $\phi_{x'}$ is separable ($=$ generically \'etale) for all $x'\in V'$.
  \end{enumerate}
\end{defn}

\begin{rem} \label{rem:after.tame.cover}
  \begin{enumerate}
  \item Since $D$ is finite, it follows from~(1)
    and Theorem~\ref{thm:harmonic.skeleta}(3) that a tame covering of
    $(X,V\cup D)$ is a tamely ramified morphism of curves.
  \item If $\Sigma$ has at least one edge, then~(2) implies~(3).
  \item Let $\phi:(X',V'\cup D')\to(X,V\cup D)$ be a tame covering,
    let $S\subset\Sigma$ be a finite set of type-$2$ points, and let 
    $S' = \phi\inv(S)$.  Then $S'\subset\Sigma(X',V'\cup D')$ is also a
    finite set of type-$2$ points, so $\phi$ also defines a tame covering
    $(X',S'\cup V'\cup D')\to(X,S\cup V\cup D)$.
  \end{enumerate}
\end{rem}

\begin{lem} \label{lem:tame.implies.tame}
  Let $\phi: (X',V'\cup D')\to(X,V\cup D)$ be a tame covering
  and let $\Sigma = \Sigma(X,V\cup D)$ and
  $\Sigma' = \Sigma(X',V'\cup D')$.  Then
  $\phi|_{\Sigma'}:\Sigma'\to\Sigma$ is a tame covering of metrized
  complexes of curves.
\end{lem}

\pf Let $\Gamma$ (resp.\ $\Gamma'$) be the augmented metric graph
underlying $\Sigma$ (resp.\ $\Sigma'$).
By Proposition~\ref{prop:tame.visible.ram} we only have to show that
$\phi|_{\Gamma'}$ is generically \'etale.   
Let $R = \sum_{x'\in V'\cup D'} R_{x'}\, (x')$ be the ramification divisor
of $\phi|_{\Gamma'}$ as defined in~\eqref{eq:Rvprime}, 
and let $S = \sum_{x'\in D'} S_{x'}\, (x')$ be the ramification divisor of
$\phi:X'\to X$.  We will show that $R = S$.
Since $\phi_{x'}$ is generically \'etale for all $x'\in V'$,
we see from~\eqref{eq:Rvprime} and the Riemann--Hurwitz formula as applied
to $\phi_{x'}$ that $R_{x'} \geq 0$.  If 
$x'\in D'$, then there is a unique edge $e'$ adjacent to $x'$ and
$d_{e'}(\phi) = d_{x'}(\phi)$, so
\[ R_{x'} = 2\,d_{x'}(\phi) - 2 - d_{x'}(\phi) + 1 = d_{x'}(\phi) - 1
= S_{x'}, \]
where the final equality is Theorem~\ref{thm:harmonic.skeleta}(3).
Since $R_{x'}\geq 0$ for all $x'\in V'$,
it is enough to show that $\deg(R) = \deg(S)$.  
By~\eqref{eq:graph.RH} we have 
$K_{\Gamma'}=(\phi|_{\Gamma'})^*(K_\Gamma)+R$, so counting degrees gives 
\[ \deg(R) = \deg(\phi|_{\Gamma'})\,(2-2g(\Gamma)) - (2-2g(\Gamma')). \]
The Lemma now follows from the equalities
$\deg(\phi)=\deg(\phi|_{\Gamma'})$, $g(\Gamma)=g(X)$, $g(\Gamma')=g(X')$,
and the Riemann--Hurwitz formula as applied to $\phi:X'\to X$.\qed

\begin{rem} \label{rem:tame.visible.ram}
  We showed in the proof of Lemma~\ref{lem:tame.implies.tame} that the
  ramification divisor of $\phi|_{\Gamma'}$ coincides with the
  ramification divisor of $\phi$.
  Moreover, it follows from Proposition~\ref{prop:tame.visible.ram} that
  for every $x'\in V'$, every ramification point $\bar x\p\in C_{x'}$ is
  the reduction of a tangent direction in $\Gamma'$.  In other words, 
  for a tame covering all
  ramification points of all residue curves are ``visible'' in the
  morphism of underlying metric graphs $\Gamma'\to\Gamma$.
\end{rem}

\begin{prop} \label{prop:branch.locus}
  Let $\phi: (X',V'\cup D')\to(X,V\cup D)$ be a finite morphism of
  triangulated punctured curves.  Then $\phi$ is a tame
  covering if and only if
  \begin{enumerate}
  \item $D$ contains the branch locus of $\phi$ and
  \item $\phi_{x'}: C_{x'}\to C_{\phi(x)}$ is tamely ramified for every
    type-$2$ point $x'\in X'^\an$.
  \end{enumerate}
  Moreover, if $\phi$ is a tame covering and
  $x'\in X'^\an$ is a type-$2$ point not contained in $\Sigma'$, then
  $\phi$ is an isomorphism in a neighborhood of $x'$.
\end{prop}

\pf It is clear that conditions~(1)--(2) imply that $\phi$ is a tame
covering, so suppose that $\phi$ is a tame 
covering and $x'\in X'^\an$ is a
type-$2$ point.  We have the following cases:
\begin{bullets}
\item If $x'\in V'$, then $\phi_{x'}$ is tamely ramified  by
  Lemma~\ref{lem:tame.implies.tame} and
  Proposition~\ref{prop:tame.visible.ram}.

\item Suppose that $x'$ is contained in the interior of an edge $e'$ of
  $\Sigma'$.  Since $\phi$ is a tame covering of
  $(X,(V\cup\{\phi(y')\})\cup D)$ by Remark~\ref{rem:after.tame.cover}(3),  
  we are reduced to the case treated above.

\item 
  Suppose that $x'\notin\Sigma'$ but 
  $y' = \tau(x')$ is contained in $V'$.  Let 
  $\fX$ (resp.\ $\fX'$) be the semistable formal model of $X$ (resp.\
  $X'$) corresponding to the semistable vertex set $V$ (resp.\ $V'$) and
  let $\phi$ denote the unique finite morphism $\fX'\to\fX$ extending
  $\phi: X'\to X$ (see Theorem~\ref{thm:vsets.extensions}).
  Let $\bar x\p = \red(x')\in\fX'(k)$ and let
  $x = \phi(x')$ and $\bar x = \red(x) = \phi_k(\bar x\p)$.  
  The connected component $B'$ of $x'$ in $X'^\an\setminus\Sigma'$ is equal
  to $\red\inv(\bar x\p)$, and the connected component $B$ of $x$ in 
  $X^\an\setminus\Sigma$ is equal to $\red\inv(\bar x)$.  Since
  $\Sigma\cap B=\emptyset$, $\phi_k$ is not branched over $\bar x$ by
  Lemma~\ref{lem:tame.implies.tame} and
  Proposition~\ref{prop:tame.visible.ram}.  Therefore
  $\phi\inv(B)$ is a disjoint union of $\deg(\phi)$ open balls (one of
  which is $B'$) mapping isomorphically onto $B$.  In particular,
  $\phi$ is an isomorphism in a neighborhood of $x'$.

\item Suppose that $x'\notin\Sigma'$ but 
  $y' = \tau(x')$ is contained in the interior of an edge
  of $\Sigma'$.  We claim that $\phi$ is an isomorphism in a neighborhood
  of $x'$.  Let $y = \phi(y')\in\Sigma$.  Then $\phi$ is a tame
  covering of $(X,(V\cup\{y\})\cup D)$, so replacing
  $V$ with $V\cup\{y\}$ and $V'$ with $V'\cup\phi\inv(y)$ (see
  Remark~\ref{rem:after.tame.cover}(3)), we are reduced 
  to the previous case.
\end{bullets}\qed

\begin{rem}
  It follows from Proposition~\ref{prop:branch.locus} that 
  if $\phi$ is a tame covering then
  the set-theoretic branch locus of $\phi: X'^\an\to X^\an$ (i.e.\ the set of
  all points in $X^\an$ with fewer than $\deg(\phi)$ points in its fiber)
  is contained in $\Sigma$.
  See~\cite{faber:ramification1,faber:ramification2} for more on
  the topic of the Berkovich ramification locus in the case of
  self-maps of $\P^1$.
\end{rem}

\section{Applications to morphisms of semistable models}
\label{sec:simultaneous.ss.reduction}

In this section we show how Section~\ref{sec:morphisms.skeleta}
formally implies a large part of the results of Liu~\cite{liu:stable_hull}
on simultaneous semistable reduction of morphisms of curves
(see also~\cite{coleman:stable_maps}) over an arbitrary valued
field.  
In addition to establishing these results over a more general ground field,
we feel that the ``skeletal'' point of view on simultaneous 
semistable reduction is enlightening (see Remark~\ref{rem:skeletal.liu}).
To this end, we 
let $K_0$ be any field equipped with a nontrivial
non-Archimedean valuation $\val: K_0\to\R\cup\{\infty\}$.  Its valuation
ring will be denoted $R_0$, its maximal ideal $\fm_{R_0}$, and its
residue field $k_0$.

Let $X$ be a smooth, proper, geometrically connected algebraic
$K_0$-curve.  By a 
\emph{(strongly) semistable $R_0$-model of $X$} we mean a
flat, integral, proper relative curve $\cX\to\Spec(R_0)$ whose special
fiber $\cX_{k_0}$ is a (strongly) semistable curve (i.e.\ 
$\cX_{\bar k_0}$ is a reduced curve with at worst ordinary double point
singularities; it is 
strongly semistable if its irreducible components are smooth)
and whose generic fiber is equipped with an
isomorphism $\cX_{K_0}\cong X$.  By properness of $\cX$, any
$K_0$-point $x\in X(K_0)$ extends in a unique way to an
$R_0$-point $x\in\cX(R_0)$; the special fiber of this point is the
\emph{reduction} of $x$ and is denoted $\red(x)\in\cX(k_0)$.
Let $D\subset X(K_0)$ be a finite set.
A semistable model $\cX$ of $X$ is a 
\emph{semistable $R_0$-model of $(X,D)$} provided that the points of $D$
reduce to distinct smooth points of $\cX_{k_0}$.  The model $\cX$ is a
\emph{stable $R_0$-model of $(X,D)$} provided that every rational (resp.\
genus-$1$) component of the normalization $\cX_k$ contains at least three
points (resp.\ one point) mapping to a singular point of $\cX_k$ or to the
reduction of a point of $D$.

If $K_0 = K$ is complete and
algebraically closed, we define a 
\emph{(strongly) semistable formal $R$-model of $X$} to be an integral, proper, 
admissible formal $R$-curve $\fX$ whose analytic generic fiber $\fX_K$ is
equipped with an isomorphism to $X^\an$ and whose special fiber $\fX_k$ is
a (strongly) semistable curve.  
There is a natural map of sets $\red:X^\an\to\fX_k$, called the
\emph{reduction}; it is surjective and anti-continuous in the sense that
the inverse image of a Zariski-open subset of $\fX_k$ is a closed subset
of $X^\an$ (and vice-versa).  Using the reduction map, for a finite set
$D\subset X(K)$ we define 
\emph{semistable} and \emph{stable formal $R$-models of $(X,D)$} as above.

As we will be passing between formal and algebraic $R$-models of
$K$-curves, it is worth stating the following lemma;
see~\cite[Remark~5.30(2)]{bpr:trop_curves}.

\begin{lem} \label{lem:algebraization}
  Let $X$ be a smooth, proper, connected $K$-curve.  The $\varpi$-adic
  completion functor defines an equivalence
  from the category of semistable $R$-models of $X$ to semistable formal
  $R$-models of $X$.
\end{lem}

\smallskip
The inverse of the $\varpi$-adic completion functor of
Lemma~\ref{lem:algebraization} will be called \emph{algebraization}.

\paragraph[Descent to a general ground field]
\label{par:descent}
We will use the following lemmas to descend 
the geometric theory of Section~\ref{sec:morphisms.skeleta} to the valued field
$K_0$.  Fix an algebraic closure 
$\bar K_0$ of $K_0$ and a valuation $\val$ on $\bar K_0$ extending the
given valuation on $K_0$.  For any field $K_1\subset\bar K_0$ we consider
$K_1$ as a valued field with respect to the restriction of $\val$, and we
write $R_1$ for the valuation ring of $K_1$.
Let $K$ be the completion of $\bar K_0$ with respect to $\val$.
This field is algebraically closed
by~\cite[Proposition~3.4.1/3]{bgr:nonarch}.

\begin{lem} \label{lem:sep.pts.dense}
    Let $X_0$ be a smooth, proper, geometrically connected curve over
    $K_0$ and let $X = X_0\tensor_{K_0} K$.  Let $K_0^\sep$ be the
    separable closure of $K_0$ in $K$.  Then the image of 
    $X_0(K_0^\sep)$ under the natural inclusion
    \[ X_0(K_0^\sep)\subset X_0(K) = X(K) \subset X^\an \]
    is dense in $X^\an$.
\end{lem}

\pf Let $\hat K_0\subset K$ denote the completion of $K_0$ and let 
$(\hat K_0)^\sep\subset K$ be its separable closure.
By~\cite[Proposition~3.4.1/6]{bgr:nonarch}, 
$(\hat K_0)^\sep$ is dense in $K$.  Let 
$\hat K_1\subset(\hat K_0)^\sep$ be a finite, separable extension of 
$\hat K_0$.  By Krasner's lemma (see \S3.4.2 of loc.\ cit.), there exists
a finite, separable extension $K_1/K_0$ contained in $K_0^\sep$ which is
dense in $\hat K_1$.  It follows that $K_0^\sep$ is dense in 
$(\hat K_0)^\sep$, so $K_0^\sep$ is dense in $K$.
Since $\P^1_K(K)$ is dense in $\P_K^{1,\an}$ and the subspace topology on
$\P^1_K(K)\subset\P_K^{1,\an}$ coincides with the ultrametric topology,
this proves the lemma for $X_0 = \P^1_{K_0}$.

For general $X_0$, choose a finite, generically \'etale morphism
$\phi:X_0\to\P^1_{K_0}$.  Let $D\subset\P^1_{K}(K)$ be the branch locus,
so $\phi:X^\an\setminus\phi\inv(D)\to\P^{1,\an}_K\setminus D$ is \'etale,
hence open by~\cite[Corollary~3.7.4]{berkovich:etalecohomology}.  It 
follows that if $U\subset X^\an$ is a nonempty open set, then 
$\phi(U)\setminus D$ contains a nonempty open subset of 
$\P^{1,\an}_K\setminus D$.  Let 
$x\in\P^1_{K_0}(K_0^\sep)$ be a point contained in
$\phi(U)\setminus D$.  Then $\phi\inv(x)\subset X_0(K_0^\sep)$  and
$\phi\inv(x)\cap U\neq\emptyset$.\qed

\begin{lem} \label{lem:descend.model}
  Let $X_0$ be a smooth, proper, geometrically connected curve over $K_0$,
  let $D\subset X_0(K_0)$ be a finite subset,
  let $X = X_0\tensor_{K_0} K$, and let $\cX$ be a semistable model of 
  $(X,D)$.  There exists a finite, separable extension $K_1$ of $K_0$
  and a semistable model $\cX_1$ of $X_1 = X_0\tensor_{K_0} K_1$
  with respect to $D$ such that $\cX_1\tensor_{R_1} R\cong\cX$.
\end{lem}

\pf Let $\fX$ be the $\varpi$-adic completion of $\cX$.
For $\bar x\in\fX_k(k)$ the subset $\red\inv(\bar x)\subset X^\an$ is open
in $X^\an$.  By Lemma~\ref{lem:sep.pts.dense}, there exists a point
$x\in X_0(K_0^\sep)$ with $\red(x) = \bar x$, so after passing to a
finite, separable extension of $K_0$ if necessary, we may enlarge 
$D\subset X_0(K_0)$ to assume that $\cX$ is a \emph{stable} model of
$(X,D)$.  On the other hand, by the stable reduction theorem (i.e.\ the
properness of the Deligne-Mumford stack $\bar\cM_{g,n}$ parameterizing
stable marked curves) there exists
a stable model $\cX_0$ of $(X_0,D_0)$ after potentially passing to a
finite, separable extension of $K_0$.  By uniqueness of stable models we
have $\cX_0\tensor_{R_0} R\cong\cX$.\qed

\begin{lem} \label{lem:descend.morphism}
  Let $X_0,X'_0$ be smooth, proper, geometrically connected $K_0$-curves
  and let $\phi: X_0'\to X_0$ be a finite morphism.  Let 
  $\cX_0$ (resp.\ $\cX_0'$) be a semistable model of 
  $X_0$ (resp.\ $X_0'$), let $X = X_0\tensor_{K_0} K$ and
  $\cX = \cX_0\tensor_{R_0} R$
  (resp.\ $X' = X_0'\tensor_{K_0} K$ and 
  $\cX' = \cX_0'\tensor_{R_0} R$), and suppose that 
  $\phi_K: X'\to X$ extends to a morphism 
  $\cX'\to\cX$.  Then $\phi$ extends to a morphism
  $\cX_0'\to\cX_0$ in a unique way.
\end{lem}

\pf The morphism $\phi$ gives rise to a section
$X_0'\to X_0\times_{K_0} X_0'$ of 
$\pr_2:X_0\times_{K_0} X_0'\to X_0'$.  Let 
$\cZ$ be the schematic closure of the image of $X_0'$ in
$\cX_0\times_{R_0}\cX_0'$.  It is enough to show that
$\pr_2:\cZ\to\cX_0'$ is an isomorphism.  Since schematic closure respects
flat base change, the schematic closure of the image of
$X'$ in $\cX\times_R\cX'$ is equal to $\cZ\tensor_{R_0} R$.  
By properness of $\cX'$, the image of $\cX'$ in
$\cX\times_R\cX'$ is then equal to $\cZ\tensor_{R_0} R$, so 
$\pr_2:\cZ\tensor_{R_0} R\to\cX'$ is an isomorphism.  Since
$R_0\to R$ is faithfully flat, this implies that
$\pr_2:\cZ\to\cX_0'$ is an isomorphism.\qed

\paragraph[The relation between semistable models and skeleta]
\label{par:models.skeleta}
Let $X$ be a smooth, proper, connected $K$-curve.
The following theorem is due to Berkovich and
Bosch-L\"utkebohmert;  
see~\cite[Theorem~5.34]{bpr:trop_curves} for more precise references.

\begin{thm} \label{thm:formal.fibers}
  Let $\fX$ be a semistable formal model of $X$ and let 
  $\bar x\in\fX_k$ be a point.  Then
  \begin{enumerate}
  \item $\bar x$ is a generic point if and only if $\red\inv(\bar x)$
    consists of a single type-$2$ point of $X^\an$.
  \item $\bar x$ is a smooth closed point if and only if
    $\red\inv(\bar x)$ is an open ball.
  \item $\bar x$ is an ordinary double point if and only if
    $\red\inv(\bar x)$ is an open annulus.
  \end{enumerate}
\end{thm}

It follows from Theorem~\ref{thm:formal.fibers} that if $\fX$ is a
semistable formal model of $X$, then
\[ V(\fX) = \big\{ \red\inv(\bar\zeta)~:~\zeta\in\fX_k \text{ is a generic
  point } \big\} \]
is a semistable vertex set of $X$: indeed,
\begin{equation} \label{eq:semistable.decomp.model}
X\setminus V(\fX) 
= \Djunion_{\substack{\bar x\in\fX_k\\\text{singular}}} \red\inv(\bar x) 
\cup \Djunion_{\substack{\bar x\in\fX_k\\\text{smooth}}} \red\inv(\bar x) 
\end{equation}
is a disjoint union of open balls and finitely many open annuli.
The following folklore theorem is proved
in~\cite[Theorem~5.38]{bpr:trop_curves}.

\begin{thm} \label{thm:models.vertex.sets}
  Let $X$ be a smooth, proper, connected algebraic $K$-curve and let
  $D\subset X(K)$ be a finite set of closed points.
  The association $\fX\mapsto V(\fX)$ sets up a bijection between the set
  of (strongly) semistable formal models of $(X,D)$ (up to isomorphism)
  and the set of (strongly) semistable vertex sets of $(X,D)$.
\end{thm}

\paragraph \label{par:metrized.complex.direct}
Let $\fX$ be a semistable formal model of $(X,D)$.  One can
construct the metrized complex of curves 
$\Sigma(X,V(\fX)\cup D)$ directly from $\fX$, as follows.  Let $V$ be the set of irreducible
components of $\fX_k$, regarded as a set of vertices of a graph.  For 
$x\in V$ we will write $C_x$ to denote the corresponding component of
$\fX_k$.  To every double point 
$\bar x\in\fX_k$ we associate an edge $e$ connecting the irreducible
components of $\fX_k$ containing $\bar x$ (this is a loop edge if there is
only one such component); if $C_x$ is such a component,
then we set $\red_x(e) = \bar x$.  
The completed local ring of $\fX$
at $\bar x$ is isomorphic to $R\ps{s,t}/(st-\pi)$ for some
$\pi\in\fm_R\setminus\{0\}$; we define the length of $e$ to be
$\ell(e) = \val(\pi)$.  
We connect a point $x\in D$ to the vertex
corresponding to the irreducible component of $\fX_k$ containing
$\red(x)$.  These data define a metrized complex of curves $\Sigma(\fX,D)$.
By~\eqref{eq:semistable.decomp.model}, the graphs underlying
$\Sigma(\fX,D)$ and $\Sigma(X,V(\fX)\cup D)$ are naturally isomorphic.
The edge lengths coincide --- see for
instance~\cite[Proposition~5.37]{bpr:trop_curves} --- and the residue
curves of the two metrized complexes are naturally isomorphic 
by~\cite[Proposition~2.4.4]{berkovich:analytic_geometry}.
From this it is straightforward to verify that  
$\Sigma(\fX,D)$ and $\Sigma(X,V(\fX)\cup D)$ are identified as
metrized complexes of curves.

In particular, a semistable formal model $\fX$ of $(X,D)$ is 
stable if and only if $V(\fX)$ is a stable vertex set of
$(X,D)$.  

As another consequence of Theorem~\ref{thm:models.vertex.sets}, we obtain the
following compatibility of skeleta and extension of scalars.

\begin{prop} \label{prop:skeleton.extend.scalars}
  Let $X$ be a smooth, proper, connected $K$-curve,
  let $K'$ be a complete and algebraically closed valued field extension
  of $K$, and let $\pi: X_{K'}^\an\to X_K^\an$ be the canonical map. 
  If $\fX$ is a semistable formal model of $X$, then
  \begin{enumerate}
  \item $\pi$ maps $V(\fX_{R'})$ bijectively onto $V(\fX)$, and
  \item $\pi$ maps $\Sigma(X_{K'},V(\fX_{R'})\cup D)$ bijectively onto
    $\Sigma(X,V(\fX)\cup D)$, with 
    \[ \pi~:~\Sigma(X_{K'},V(\fX_{R'})\cup D)\isom\Sigma(X,V(\fX)\cup D) \]
    an isomorphism of augmented metric graphs.
  \end{enumerate}
\end{prop}

\pf The reduction map is compatible with extension of scalars, in that the
following square commutes
\[\xymatrix{
  {X_{K'}^\an} \ar[r]^\red \ar[d]_\pi & {\fX_{k'}} \ar[d]^{\bar\pi} \\
  {X^\an} \ar[r]_\red & {\fX_k}
}\]
with $k'$ the residue field of $K'$.  The first assertion is immediate
because the extension of scalars morphism $\bar\pi:\fX_{k'}\to\fX_k$ is a
bijection on generic points.  It is also a bijection on singular points,
and for any node $\bar x\p\in\fX_{k'}$ the open annulus
$A' = \red\inv(\bar x\p)$ is identified with the extension of scalars
of $A = \red\inv(\bar \pi(\bar x\p))$.  It is easy to see that
$\pi: A'\to A$ takes the skeleton of $A'$ isomorphically onto the skeleton
of $A$.\qed

\paragraph[Extending morphisms to semistable models: analytic criteria]
Let $X$ be a smooth, proper, connected $K$-curve and
let $\fX$ be a semistable formal model of $X$.  The 
\emph{inverse image topology} on $X^\an$ is the topology $\cT(\fX)$ whose
open sets are the sets of the form $\red\inv(\fU_k)$, where $\fU_k$ is a
Zariski-open subset of $\fX_k$.  Note that any such set is \emph{closed}
in the natural topology on $X^\an$.  If $\fU_k\subset\fX_k$
is an affine open subset, then $\fU_k$ is the special fiber of a formal
affine open $\fU = \Spf(\cA)\subset\fX$, the set underlying the 
generic fiber $\fU_K = \sM(\cA_K)\subset X^\an$ is equal to
$\red\inv(\fU_k)$, and $\cA$ is equal to the ring $\mathring\cA_K$ of power-bounded
elements in the affinoid algebra $\cA_K$.  This essentially means that $\fX$
is a \emph{formal analytic variety} in the sense
of~\cite{bosch_lutk:uniformization};
see~\cite[Remark~5.30(3)]{bpr:trop_curves} for an explanation.  
Therefore $\fX$ is constructed by gluing the canonical models 
of the affinoid subdomains of $X^\an$ which are open in the inverse image
topology $\cT(\fX)$, along the canonical models of their intersections;
here the \emph{canonical model} of an affinoid space $\sM(\cB)$ is by
definition $\Spf(\mathring\cB)$.  The following general fact about formal
analytic varieties follows from these observations and the functoriality of the
reduction map.

\begin{prop} \label{prop:morphism.fav}
  Let $X,X'$ be smooth, proper, connected algebraic curves over $K$,
  let $\phi: X'\to X$ be a finite morphism, and let $\fX$ and $\fX'$ be
  semistable formal models of $X$ and $X'$, respectively.  
  Then $\phi$ extends to a morphism
  $\fX'\to\fX$ if and only if $\phi$ is continuous with respect to the
  inverse image topologies $\cT(\fX)$ and $\cT(\fX')$, in which case there
  is exactly one morphism $\fX'\to\fX$ extending $\phi$.
\end{prop}

\smallskip
Fix smooth, proper, connected algebraic
$K$-curves $X,X'$ 
and a finite morphism $\phi: X'\to X$.  The following theorem is
well-known to experts, although no proof appears in the literature to the
best of our knowledge.

\begin{thm} \label{thm:vsets.extensions}
  Let $\fX$ and $\fX'$ be semistable formal models of $X$ and $X'$,
  respectively.  Then $\phi: X'\to X$ extends to a morphism $\fX'\to\fX$
  if and only if $\phi\inv(V(\fX))\subset V(\fX')$, and $\fX'\to\fX$ is
  finite if and only if $\phi\inv(V(\fX)) = V(\fX')$.
\end{thm}

\begin{rem} \label{rem:induced.harmonic.morphism.2}
If $\fX$ and $\fX'$ are semistable formal models of the punctured curves $(X,D)$ and $(X',D')$, respectively, and $\phi : X' \to X$ is a finite
morphism with $\phi\inv(D)=D'$ which extends to a morphism $\fX'\to\fX$,
then since $\phi\inv(V(\fX))\subset V(\fX')$ it follows from Remark~\ref{rem:induced.harmonic.morphism.1}
that there is a natural harmonic morphism of metric graphs $\Sigma(X',V(\fX') \cup D') \to \Sigma(X,V(\fX) \cup D)$.
The morphism $\fX' \to \fX$ is finite if and only if the local degree of this harmonic morphism at every $v' \in V(\fX')$ is positive.
\end{rem}

\smallskip
Before giving a proof of Theorem~\ref{thm:vsets.extensions}, we mention
the following consequences.
Let $\fX_1$ and $\fX_2$ be semistable formal models of $X$.  We say that
$\fX_1$ \emph{dominates} $\fX_2$ provided that there exists a (necessarily
unique) morphism $\fX_1\to\fX_2$ inducing the identity map on analytic
generic fibers.   
Taking $X = X'$ in Theorem~\ref{thm:vsets.extensions},
we obtain the following corollary; this gives a different proof of the
second part of~\cite[Theorem~5.38]{bpr:trop_curves}.

\begin{cor} \label{cor:domination.containment}
  Let $\fX_1$ and $\fX_2$ be semistable formal models of $X$.  Then
  $\fX_1$ dominates $\fX_2$ if and only if $V(\fX_1)\supset V(\fX_2)$. 
\end{cor}

\smallskip
In conjunction with Proposition~\ref{prop:skeleton.extend.scalars}, we
obtain the following corollary.  If $K'$ is a complete and algebraically
closed field extension of $K$, we say that a semistable formal model
$\fX'$ of $X_{K'}$ is \emph{defined over $R$} provided that it arises as
the extension of scalars of a (necessarily unique) semistable formal model
of $X$.

\begin{cor} \label{cor:dominated.means.rational}
  Let $K'$ be a complete and algebraically closed field extension of $K$,
  let $\fX'$ be a semistable formal model of $X_{K'}$, and suppose that
  there exists a semistable formal model $\fX$ of $X$ such that 
  $\fX_{R'}$ dominates $\fX'$.  Then $\fX'$ is defined over $R$.
\end{cor}

\pf Let $\pi:X_{K'}^\an\to X^\an$ be the canonical map.
By Corollary~\ref{cor:domination.containment} we have
$V(\fX_{R'})\supset V(\fX')$, and by 
Proposition~\ref{prop:skeleton.extend.scalars} the map $\pi$ defines a
bijection $V(\fX_{R'})\isom V(\fX)$ and an isomorphism
\[ \Sigma(X_{K'},V(\fX_{R'}))\isom\Sigma(X,V(\fX)). \]
Let $V = V(\fX)$ and $V' = \pi(V(\fX'))$.  We claim that $V'$ is a
semistable vertex set.  Granted this, letting $\fY$ be the semistable
formal model of $X$ associated to $V'$, we have 
$V(\fY_{R'}) = V(\fX')$ since $V(\fY_{R'})\subset V(\fX_{R'})$ (as 
$\fX_{R'}$ dominates $\fY_{R'}$) and 
$\pi(V(\fY_{R'})) = V(\fY)$, so $\fY_{R'} = \fX'$.

It remains to prove the claim that $V'$ is a semistable vertex set.  The
point is that the question of whether or not a subset of $V$ is a semistable
vertex set is intrinsic to the augmented $\Lambda$-metric graph 
$\Sigma = \Sigma(X,V)\cong\Sigma(X_{K'},V(\fX_{R'}))$.  We leave the
details to the reader.\qed

\smallskip
We will use the following lemmas in the proof of
Theorem~\ref{thm:vsets.extensions}.  
Recall from~\parref{par:semistable.decomposition} that 
if $V$ is a semistable vertex set of $X$, then a connected
component $C$ of $X^\an\setminus V$ is \emph{adjacent} to a point
$x\in V$ provided that the closure of $C$ in $X^\an$ contains $x$.

\begin{lem} \label{lem:open.in.TX}
  Let $\fX$ be a semistable formal model of $X$, let $V = V(\fX)$, let
  $\fU_k\subset\fX_k$ be a subset, and let
  $U = \red\inv(\fU_k)\subset X^\an$.  Then
  $U$ is open in the topology $\cT(\fX)$
  if and only if the following conditions hold:
  \begin{enumerate}
  \item $U$ is closed in the ordinary topology on $X^\an$, and
  \item for every  $x\in U\cap V$, all but finitely many connected components of 
    $X^\an\setminus V$ which are adjacent to $x$ are contained in $U$.
  \end{enumerate}
\end{lem}

\pf Let $\bar\zeta$ be a generic point of $\fX_k$ and let
$\bar y\in\fX_k$ be a closed point.  Let $x\in V$ be the unique point of $X^\an$
reducing to $\bar\zeta$ and let $B = \red\inv(\bar y)$.
By the anti-continuity of the reduction map, $\bar y$ is in the closure of
$\{\bar\zeta\}$ if and only if $x$ is adjacent to 
$B$.  The lemma follows easily from this and the fact
that the connected components of $X^\an\setminus V$ are exactly the
inverse images of the closed points of $\fX_k$ under $\red$.\qed

\begin{lem} \label{lem:invV.in.Vp}
  Let $V$ and $V'$ be semistable vertex sets of $X$ and $X'$,
  respectively, and suppose that
  $\phi\inv(V)\subset V'$.  Let $C'$ be a connected component of
  $X'^\an\setminus\phi\inv(V)$.  Then $C'$ has the following form:
  \begin{enumerate}
  \item If $C'$ intersects $\Sigma(X',V')$, 
    then $C' = \tau\inv\big(\Sigma(X',V')\cap C'\big)$, and
  \item otherwise $C'$ is an open ball connected component of
    $X'^\an\setminus V'$.
  \end{enumerate}
\end{lem}

\pf Let $\Sigma' = \Sigma(X',V')$.
Suppose that there exists $y'\in C'$ such that $\tau(y')\notin C'$.  
Let $B'$ be the connected component of $X'^\an\setminus\Sigma'$ containing
$y'$, so $B'$ is an open ball contained in $C'$ and $\tau(y')$ is the end of $B'$.
It follows that $\tau(y')$ is contained in the
closure of $C'$ in $X'^\an$.  Since $C'$ is a connected component of
$X'^\an\setminus\phi\inv(V)$, its closure is contained in 
$C'\cup\phi\inv(V)$, so $\tau(y')\in\phi\inv(V)$.   Therefore
$B' = (B'\cup\{\tau(y')\})\cap(X'^\an\setminus\phi\inv(V))$ is open and
closed in $X'^\an\setminus\phi\inv(V)$, so $B'$ is a 
connected component of $X'^\an\setminus\phi\inv(V)$ and hence
$B' = C'$.

Now suppose that $\tau(C')\subset C'$.  Then 
$\tau(C')\subset C'\cap\Sigma'$, so $C'\subset\tau\inv(C'\cap\Sigma')$.  
Let $y'\in\tau\inv(C'\cap\Sigma')$ and let $B'$ be the connected component
of $X'^\an\setminus\Sigma'$ containing $y'$, as above.  Since
$\bar B\p = B'\cup\{\tau(y')\}$ is a connected subset of $X'^\an\setminus\phi\inv(V)$
intersecting $C'$ we have $\bar B\p\subset C'$, so $y'\in C'$ and
therefore $C' = \tau\inv(C'\cap\Sigma')$.\qed

\pf[of Theorem~\ref{thm:vsets.extensions}]
Let $V = V(\fX)$ and $V' = V(\fX')$.
If there is an extension $\fX'\to\fX$ of $\phi$, then the square
\begin{equation} \label{eq:vsets.extension} \xymatrix{
  {X'^\an} \ar[r]^\phi \ar[d]_\red & {X^\an} \ar[d]^\red \\
  {\fX'_k} \ar[r] & {\fX_k}
}\end{equation}
commutes.  Let $x'\in\phi\inv(V)$, so $x = \phi(x')$ reduces to a generic
point $\bar\zeta$ of $\fX_k$.  Since the reduction $\bar\zeta\p$ of $x'$ maps to 
$\bar\zeta$,
the point $\bar\zeta\p$ is generic, so $x'\in V'$.  Therefore 
$\phi\inv(V)\subset V'$.  The morphism $\fX'\to\fX$ is finite
if and only if every generic point of $\fX'_k$ maps to a generic point
of $\fX_k$; as above,  this is equivalent to
$V'= \phi\inv(V)$. 

It remains to prove that if $\phi\inv(V)\subset V'$, then $\phi$ extends to
a morphism $\fX'\to\fX$.  By Proposition~\ref{prop:morphism.fav}, we must
show that $\phi$ is continuous with respect to the topologies $\cT(\fX)$
and $\cT(\fX')$.  Let $U\subset X^\an$ be $\cT(\fX)$-open and let 
$U' = \phi\inv(U)$.  Clearly $U$ is closed, so
$U'$ is closed (with respect to the ordinary topologies).  We must show
that condition~(2) of Lemma~\ref{lem:open.in.TX} holds for $U'$.  Let
$x'\in U'\cap V'$ and let $x = \phi(x')$.  If $x\notin V$, then 
let $C\subset U$ be the connected component of $X^\an\setminus V$
containing $x$ and let $C'$ be the connected component of 
$X'^\an\setminus\phi\inv(V)$ containing $x'$.  
Then $C'\subset U'$ because $\phi(C')\subset C$, and $C'$ contains every
connected component of $X'^\an\setminus V'$ adjacent to $x'$ since $C'$ is
an open neighborhood of $x'$.

Now suppose that $x\in V$.  Any connected component of 
$X'^\an\setminus V'$ which is adjacent to $x'$ maps into a connected
component of $X^\an\setminus V$ which is adjacent to $x$.
There are finitely many connected components
$X^\an\setminus V$ which are adjacent to $x$ and not contained in $U$ by 
Lemma~\ref{lem:open.in.TX}.  Let $C$ be such a component.  Since $\phi$ is
finite, there are only finitely many connected components of
$\phi\inv(C)$; 
each of these is a connected component of $X'^\an\setminus\phi\inv(V)$.
If $C'$ is such a connected component, then either $C'$ is an open ball
connected component of $X'^\an\setminus V'$ or $C'$ intersects
$\Sigma' = \Sigma(X',V')$ by Lemma~\ref{lem:invV.in.Vp}.  
There are finitely many connected components
of $X'^\an\setminus V'$ which intersect $\Sigma'$ ---
these are just the open annulus connected components of $X'^\an\setminus
V'$ --- so there are only finitely many connected components of 
$X'^\an\setminus V'$ contained in $\phi\inv(C)$.  Therefore all but
finitely many connected components of $X'^\an\setminus V'$ which are
adjacent to $x'$ map to connected components of $X^\an\setminus V$ which
are contained in $U$.\qed

\paragraph[Simultaneous semistable reduction theorems]
Recall that $K_0$ is a field equipped with a nontrivial non-Archimedean
valuation $\val: K_0\to\R\cup\{\infty\}$.
As above we fix an algebraic closure 
$\bar K_0$ of $K_0$ and a valuation $\val$ on $\bar K_0$ extending the
given valuation on $K_0$; for any field $K_1\subset\bar K_0$ we consider
$K_1$ as a valued field with respect to the restriction of $\val$, and we
write $R_1$ for the valuation ring of $K_1$.  
In what follows, $X$ and $X'$ are smooth, proper, geometrically
connected $K_0$-curves and $\phi: X'\to X$ is a finite morphism.  

\begin{prop}
  Let $\cX$ be a semistable $R_0$-model of $X$.  If there exists a
  semistable $R_0$-model $\cX'$ of $X'$ such that $\phi: X'\to X$ extends
  to a \emph{finite} morphism $\cX'\to\cX$, then there is exactly one such
  model $\cX'$ up to (unique) isomorphism.
\end{prop}

\pf Suppose first that $K = K_0$ is complete and algebraically closed.  
Let $\fX$ (resp.\ $\fX'$) be the
$\varpi$-adic completion of $\cX$ (resp.\ $\cX'$).  
In this case, the proposition follows from
Theorem~\ref{thm:models.vertex.sets} and the fact that if the morphism 
$\fX'\to\fX$ of Theorem~\ref{thm:vsets.extensions} is finite, then 
$V(\fX')$ is uniquely determined by $V(\fX)$.  The general case follows
from this case after passing to the completion of the algebraic closure
$K$ of $K_0$: if $\cX',\cX''$ are two semistable models of
$X'$ such that $\phi$ extends to finite morphisms $\cX'\to\cX$ and
$\cX''\to\cX$, then the isomorphism $\cX''_R\isom\cX'_R$ descends to an
isomorphism $\cX''\isom\cX'$ by Lemma~\ref{lem:descend.morphism}.\qed

\smallskip
Let $\cX_1,\cX_2$ be semistable $R_0$-models of $X$.  Recall that $\cX_1$
\emph{dominates} $\cX_2$ if there exists a (necessarily unique) morphism
$\cX_1\to\cX_2$ inducing the identity on $X$.

\begin{prop} \label{prop:can.extend}
  Let $\cX$ and $\cX'$ be semistable $R_0$-models of $X$ and $X'$,
  respectively.  Let $D'\subset X'(K_0)$ be a finite set of points. Then
  there exists a finite, separable extension $K_1$ of 
  $K_0$ and a semistable $R_1$-model $\cX''$ of 
  $(X'_{K_1},D')$ such that:
  \begin{enumerate}
  \item $\cX''$ dominates $\cX'_{R_1}$,
  \item $\phi_{K_1}: X'_{K_1}\to X_{K_1}$ extends to a morphism
    $\cX''\to\cX_{R_1}$, and 
  \item any other semistable $R_1$-model of $(X'_{K_1},D')$ satisfying the
    above two properties dominates $\cX''$.
  \end{enumerate}
  Moreover, the formation of $\cX''$ commutes with arbitrary valued field
  extensions $K_1\to K_1'$.
\end{prop}

\pf Suppose first that $K_0 = K$ is complete and algebraically closed.
Let $\fX$ and $\fX'$ be the $\varpi$-adic completions of $\cX$ and $\cX'$,
respectively.  Let $V = V(\fX)$ and $V' = V(\fX')$.  
By Lemma~\ref{lem:larger.vertex.set}, there is a minimal semistable vertex
set $V''$ of $(X',D')$ which contains $\phi\inv(V)\cup V'$.
Let $\fX''$ be the semistable formal model of $X'$ corresponding to $V''$.
Then $\fX''$ dominates $\fX'$ by Corollary~\ref{cor:domination.containment}
and $\phi$ extends to a morphism 
$\fX''\to\fX$ by Theorem~\ref{thm:vsets.extensions}.
Part~(3) follows from Corollary~\ref{cor:domination.containment} and the
minimality of $V''$.  Taking $\cX''$ to be the algebraization of $\fX''$
yields (1)--(3) in this case.

For a general valued field $K_0$, suppose that $K$ is the completion of
$\bar K_0$.  By Lemma~\ref{lem:descend.model} the model of $X'_K$
constructed above descends to a model $\cX''$ defined over the ring of
integers of a finite, separable extension $K_1$ of $K_0$.  Properties
(1)--(3) follow from Lemma~\ref{lem:descend.morphism} and the
corresponding properties of $\cX''_R$.

Now we address the behavior of this construction with respect to base
change.  Using Lemma~\ref{lem:descend.morphism} we immediately reduce to
the case of an extension $K\to K'$ of complete and algebraically closed
valued fields.  Let $R'$ be the ring of integers of $K'$.
Let $\cX''$ (resp.\ $\cX'''$) be the minimal $R$-model of
$X'$ dominating $\cX'$ (resp.\ $R'$-model of $X'_{K'}$ dominating
$\cX'_{R'}$) mapping to $\cX$ (resp.\ $\cX_{R'}$).  By~(3) as applied to
$\cX'''$, we have that $\cX''_{R'}$ dominates $\cX'''$.  By
Corollary~\ref{cor:dominated.means.rational}, $\cX'''$ is defined over
$R$, so by~(3) as applied to $\cX''$, we have $\cX''_{R'} = \cX'''$.\qed

\begin{thm}[Liu] \label{thm:stable.hull}
  Let $\cX$ (resp.\ $\cX'$) be a semistable $R_0$-model of $X$ (resp.\
  $X'$).  Let $D\subset X(K_0)$ and $D'\subset X'(K_0)$ be finite sets, and
  suppose that $\phi(D')\subset D$.  Then there exists a finite, separable
  extension $K_1$ of $K_0$ and semistable $R_1$-models
  $\cX_1,\cX_1'$ of $(X_{K_1},D),\,(X'_{K_1},D')$, respectively, such that
  \begin{enumerate}
  \item $\cX_1$ dominates $\cX_{R_1}$ and $\cX_1'$ dominates $\cX'_{R_1}$,
  \item the morphism $\phi_{K_1}: X'_{K_1}\to X_{K_1}$ extends to a finite
    morphism $\cX_1'\to\cX_1$, and
  \item if $\cX_2$, $\cX_2'$ are semistable formal models of 
    $(X_{K_1},D),\, (X'_{K_1},D')$, respectively, satisfying~(1) and~(2)
    above, then $\cX_2$ dominates $\cX_1$ and $\cX_2'$ dominates
    $\cX_1'$. 
  \end{enumerate}
  Moreover, the formation of $\cX_1'\to\cX_1$ commutes with arbitrary
  valued field extensions $K_1\to K_1'$.
\end{thm}

\smallskip
The morphism $\cX_1'\to\cX_1$ is called the \emph{stable marked hull} of 
$\cX'\dashrightarrow\cX$ in~\cite{liu:stable_hull}.
\smallskip

\pf First assume that $K_0 = K$ is complete and algebraically closed; let 
$\fX,\fX'$ be the $\varpi$-adic completions of $\cX,\cX'$, respectively.
Let $V = V(\fX)$ and $V' = V(\fX')$.  
By Theorem~\ref{thm:models.vertex.sets},
Theorem~\ref{thm:vsets.extensions}, and
Corollary~\ref{cor:domination.containment}, we may equivalently formulate
the existence and uniqueness of $\fX_1'\to\fX_1$ in terms of semistable
vertex sets, as follows. We must prove 
that there exists a semistable vertex set $V_1$ of $(X,D)$
such that
\begin{enumerate}
\item $V_1$ contains $V\cup\phi(V')$,
\item $\phi\inv(V_1)$ is a semistable vertex set of $(X',D')$,
  and
\item $V_1$ is minimal in the sense that if $V_2$ is another semistable
  vertex set of $(X,D)$ satisfying~(1) and~(2) above, then
  $V_2\supset V_1$.
\end{enumerate}
First we will prove the existence of $V_1$ satisfying~(1) and~(2).  By
Lemma~\ref{lem:larger.vertex.set} we may enlarge $V$ to assume that $V$ is
a semistable vertex set of $(X,D)$.  Let  
$\Sigma = \Sigma(X,V\cup D)$ and $\Sigma' = \Sigma(X',V')$.
By Corollary~\ref{cor:stable.hull.skel}, there exists a skeleton
$\Sigma_1$ of $(X,D)$ 
such that $\Sigma_1\supset\Sigma\cup\phi(\Sigma')$ and such that
$\Sigma_1' = \phi\inv(\Sigma_1)$ is a skeleton of $(X',\phi\inv(D))$.
Let $V_1'$ be a vertex set for $\Sigma_1'$.  Since any finite subset of
type-$2$ points of $\Sigma_1'$ which contains $V_1'$ is again a
vertex set for $\Sigma_1'$ by Proposition~\ref{prop:skeleton.properties}(4),
we may and do assume that $V'\subset V_1'$.
Let $V_1\subset\Sigma_1$ be the union of a vertex set for $\Sigma_1$
with $V\cup\phi(V_1')$.   Then $V_1$ is a
semistable vertex set of $(X,D)$,
and $\phi\inv(V_1)\subset\Sigma_1'$ is a finite set
of type-$2$ points containing $V_1'$, thus is a semistable vertex set of
$(X',\phi\inv(D))$ (hence of $(X',D')$ as well).

To prove that there exists a minimal such $V_1$, we make the following
recursive construction.  Let $V(0)=V$, let $V'(0) = V'$, and for each
$n\geq 1$ let $V(n)$ be the minimal semistable vertex set of 
$(X,D)$ containing $V(n-1)\cup\phi(V'(n-1))$ and let $V'(n)$ be the minimal
semistable vertex set of $(X',D')$ containing
$\phi\inv(V(n))$.  These sets exist by Lemma~\ref{lem:larger.vertex.set}.
By induction it is clear that if $V_2$ is any semistable vertex set of 
$(X,D)$ satisfying~(1) and~(2) above, then
$V(n)\subset V_2$ for each $n$.  Since $V_2$
is a finite set, for some $n$ we have $V(n) = V(n+1)$, which is to say
that $\phi(V'(n))\subset V(n)$; since $V'(n)\supset\phi\inv(V(n))$, we
have that $V'(n) = \phi\inv(V(n))$ is a semistable vertex set of 
$(X',D')$.  Then
$V_1 = V(n)$ is the minimal semistable vertex set satisfying~(1) and~(2)
above.

The case of a general ground field reduces to the geometric case handled
above exactly as in the proof of Proposition~\ref{prop:can.extend}, as
does the statement about the behavior of the stable marked hull
with respect to valued field extensions.\qed  

\begin{rem}[The skeletal viewpoint on Liu's theorem] 
  \label{rem:skeletal.liu}
  The statement of Theorem~\ref{thm:stable.hull} is strongly analogous to
  Corollary~\ref{cor:stable.hull.skel}, which is indeed the main
  ingredient in the proof.  The difference is that whereas finite morphisms of
  semistable models correspond to pairs $V,V'$ of semistable vertex sets
  such that $\phi\inv(V) = V'$, for finite morphisms of
  triangulated punctured curves one requires in addition that 
  $\phi\inv(\Sigma) = \Sigma'$.  The former condition (of Liu's theorem)
  does not imply the latter: for instance, let $X'$ be the Tate curve
  $y^2 = x^3 - x + \varpi$, let $X = \P^1$, and let
  $\phi:X'\to X$ be the cover $(x,y)\mapsto x$.
  This extends to a finite morphism of semistable models given by the
  same equations; however, the associated skeleton of $X'$ (resp.\ of
  $\P^1$) is a circle (resp.\ a point), so 
  $\phi\inv(\Sigma)\subsetneq\Sigma'$.  

  In general, if $\phi\inv(V) = V'$ and $\phi\inv(D) = D'$, then 
  by Remark~\ref{rem:structure.of.image.skeleton} the image of
  the skeleton $\Sigma(X',V'\cup D')$ is equal to $\Sigma(X,V\cup D)$
  union a finite number of geodesic segments $T_1,\ldots,T_r$, and by
  Theorem~\ref{thm:inv.skeleton} the saturation
  $\phi\inv(\phi(\Sigma(X',V'\cup D'))$ is a skeleton of
  $X'$.

  Theorem~\ref{thm:map.stable.models} below follows the same philosophy in
  deriving a simultaneous semistable reduction theorem of Liu--Lorenzini
  from Proposition~\ref{prop:stable.skel.inverse}.
\end{rem}

\begin{rem}
  Liu in fact works over an arbitrary Dedekind scheme (a connected
  Noetherian regular scheme of dimension $1$), which includes discrete
  valuation rings but not more general valuation rings.  In his statement
  of Theorem~\ref{thm:stable.hull} the given models $\cX,\cX'$ are allowed
  to be any integral, projective $\Spec(R_0)$-schemes with generic fibers
  $X$ and $X'$, respectively.  Although we restrict to semistable
  models $\cX,\cX'$ in this paper, using a more general notion of a
  triangulation our methods can be extended to treat relatively normal models
  $\cX,\cX'$.
\end{rem}

\smallskip
The following simultaneous stable reduction theorem can be found
in~\cite{liu_lorenzini:models}.  Theorem~\ref{thm:map.stable.models} is to
Proposition~\ref{prop:stable.skel.inverse}
as Theorem~\ref{thm:stable.hull} is to Corollary~\ref{cor:stable.hull.skel} 

\begin{thm}[Liu--Lorenzini] \label{thm:map.stable.models}
  Suppose that $(X,D)$ and $(X',D')$ are both stable and that 
  $\phi\inv(D) = D'$.  Assume that $(X,D)$ and $(X',D')$ admit stable
  models $\cX$ and $\cX'$, respectively, defined over $R_0$.
  Then $\phi:X'\to X$ extends to a (not necessarily finite) morphism $\cX'\to\cX$. \end{thm}

\pf By Lemma~\ref{lem:descend.morphism} we may assume that $K_0 = K$ is
complete and algebraically closed.  
Let $V$ (resp.\ $V'$) be the minimal semistable vertex set of
$(X,D)$ (resp.\ $(X',D')$).  By Theorem~\ref{thm:vsets.extensions}, we
must show that  $\phi\inv(V)\subset V'$.
Let $\Sigma = \Sigma(X,V)$ and $\Sigma' = \Sigma(X',V')$.  
Let $V_1$ be a vertex set for $\Sigma$ with respect to which $\Sigma$ has
no loop edges.  Recall that 
\[ V = \{ x\in V_1~:~ g(x)\geq 1 \text{ or } x \text{ has valency at least
} 3 \} \]
by Proposition~\ref{prop:minimal.skeleton}(2).
Let $V_1' = \phi\inv(V_1)\cup V'$; this is a vertex set for $\Sigma'$
since $\phi\inv(\Sigma)\subset\Sigma'$ by
Proposition~\ref{prop:stable.skel.inverse}(2).
If $x\in V$ has genus at least $1$, then any $x'\in\phi\inv(x)$ has
genus at least $1$, so $\phi\inv(x)\subset V'$.  Let $x\in V$ have genus
$0$, so $x$ has valency at least $3$ in $\Sigma$.    
Let $x'\in\phi\inv(x)\subset\Sigma'$ and let $U'$ be the set containing $x'$ and all of
the connected components of $X'^\an\setminus V_1'$ adjacent to $x'$.  Then
$U'$ is an open neighborhood of $x'$, so $\phi(U')$ is an open
neighborhood of $x$ by~\cite[Lemma~3.2.4]{berkovich:analytic_geometry}.  
Let $A$ be an open annulus connected component of $X^\an\setminus V_1$
adjacent to $x$.  Any open neighborhood of $x$ intersects $A$, so there
exists $y'\in U'$ with $\phi(y')\in A$.  Let $A'$ be the connected
component of $U'\setminus\{x'\}$ containing $y'$, so $\phi(A')\subset A$.
Since $x'$ is an end of $A'$ and $\phi(x')$ is an end of $A$,
Lemma~\ref{lem:map.ball.annulus} implies that $A'$ is an open annulus.
Since $x$ has valency at least $3$ in $\Sigma$ (and since $\Sigma$ has no loop edges),
there are at least $3$ distinct open annulus connected components of
$X^\an\setminus V_1$ adjacent to $x$.  By the above, the same is true of $x'$,
so $x'$ has valency at least $3$ in $\Sigma'$, so $x'\in V'$, as
desired.\qed

\section{A local lifting theorem}\label{sec:local lifting}

We now begin with the classification of lifts of harmonic morphisms of
metrized complexes of curves to finite morphisms of algebraic curves.  We
begin by working in the neighborhood of a vertex of a metrized complex.

\paragraph 
Let $X$ be a smooth, proper, connected curve over $K$ and let $x\in X^\an$
be a type-$2$ point.  Let $V\subset X^\an$ be a strongly semistable
vertex set of $X$ containing $x$, let 
$\Sigma = \Sigma(X,V)$, and let $\tau:X^\an\to\Sigma$ be the
retraction.  Let $e_1,\ldots,e_r$ be the edges of $\Sigma$ adjacent to
$x$ and let $\Sigma_0 = \{x\}\cup e_1^\circ\cup\cdots\cup e_r^\circ$,
where for a (closed) edge $e$ we let $e^\circ$ denote the corresponding
open edge, i.e.\ the edge without its endpoints.
Then $\Sigma_0$ is an open neighborhood of $x$ in $\Sigma$ and
$\tau\inv(\Sigma_0)$ is an open neighborhood of $x$ in $X^\an$.
Following~\cite[5.54]{bpr:trop_curves} and~\cite{berkovich:etalecohomology}, we define a 
\emph{simple neighborhood} of $x$ to be an open neighborhood of this form 
(for some choice of $V$).
The connected components of 
$\tau\inv(\Sigma_0)\setminus\{x\}$ are open balls and the open annuli
$\tau\inv(e_1^\circ),\ldots,\tau\inv(e_r^\circ)$.

\begin{defn}
  A \emph{star-shaped curve} is a pointed $K$-analytic space $(Y,y)$ which is
  isomorphic to $(U,x)$ where $x$ is a type-$2$ point in the
  analytification of a smooth, proper, connected curve over $K$ and
  $U$ is a simple neighborhood of $x$.
  The point $y$ is called the \emph{central vertex} of $Y$.
\end{defn}

\paragraph \label{par:star.shaped.basic}
Let $(Y,y)$ be a star-shaped curve, so $Y\setminus\{y\}$ is a disjoint
union of open balls and finitely many open annuli $A_1,\ldots,A_r$.  The
\emph{skeleton} of $Y$ is defined to be the set 
$\Sigma(Y,\{y\}) = \{y\}\cup\bigcup_{i=1}^r\Sigma(A_i)$.  A 
\emph{compatible divisor in $Y$} is a finite set $D\subset Y(K)$ whose
points are contained in distinct open ball connected components of
$Y\setminus\{y\}$, so the connected components of 
$Y\setminus(\{y\}\cup D)$ are open balls, the open annuli
$A_1,\ldots,A_r$, and (finitely many) open balls $B_1,\ldots,B_s$
punctured at a point of $D$.  The data $(Y,y,D)$ of a star-shaped curve
along with a compatible divisor is called a 
\emph{punctured star-shaped curve}.  The 
\emph{skeleton of $(Y,y,D)$} is the set 
\[ \Sigma(Y,\{y\}\cup D) = \{y\}\cup D\cup\bigcup_{i=1}^r\Sigma(A_i)
\cup\bigcup_{j=1}^s\Sigma(B_j). \]
Fix a compatible divisor $D$ in $Y$, and let 
$\Sigma_0 = \Sigma(Y,\{y\}\cup D)$.  
There is a canonical continuous retraction map $\tau: Y\to\Sigma_0$ defined exactly
as for skeleta of algebraic curves~\parref{par:semistable.decomposition}.
The connected components of $\Sigma_0\setminus\{y\}$ are called the
\emph{edges} of $\Sigma_0$; an edge is called \emph{infinite} if it contains
a point of $D$ and \emph{finite} otherwise.  

If $Y$ is the simple neighborhood $\tau\inv(\Sigma_0)$ of
$x\in X^\an$ as above, then 
$\Sigma(Y,\{x\}) = \Sigma_0 = \Sigma\cap Y$.  A finite set $D\subset Y(K)$
is compatible with $Y$ if and 
only if $V$ is a semistable vertex set for $(X,D)$, in which
case $\Sigma(Y,\{x\}\cup D) = Y\cap\Sigma(X,V\cup D)$.  The retraction
$\tau: Y\to\Sigma(Y,\{x\}\cup D)$ is the restriction of the canonical retraction
$\tau:X^\an\to\Sigma(X,V\cup D)$.

\paragraph
Let $(Y,y)$ be a star-shaped curve.  Then $Y$ is proper as a $K$-analytic 
space if and only if all connected components of $Y\setminus\{y\}$ are
open balls, i.e.\ if and only if $\Sigma(Y,\{y\}) = \{y\}$.  
If $Y$ is proper, then there is a smooth, proper, connected curve $X$ over
$K$ and an isomorphism $f: Y\isom X^\an$.  Let $x = f(y)$.  Then
$\{x\}$ is a semistable vertex set of $X$, so
by Theorem~\ref{thm:models.vertex.sets} there is a unique smooth formal model
$\fX$ of $X^\an$ such that $x$ reduces to the generic point of $\fX_k$.  Let
$D\subset Y(K)$ be a finite set.  Then $D$ is compatible with $Y$ if and
only if the points of $f(D)$ reduce to distinct closed points of $\fX_k$. 

Conversely, let $\fX$ be a smooth, proper, connected formal curve over
$\Spf(R)$.
If $x\in X^\an$ is the point reducing to the generic point of $\fX_k$, then
$(\fX_K,x)$ is a proper star-shaped curve.

\begin{prop} \label{prop:compactify.star}
  Let $(Y,y)$ be a star-shaped curve.  Then $(Y,y)$ is isomorphic to
  $(U,x)$ where $x$ is the central vertex of a proper star-shaped
  curve $X$ and $U$ is a simple neighborhood of $x$. 
\end{prop}

\pf Let $A_1,\ldots,A_r$ be the open annulus connected components of 
$Y\setminus\{y\}$.  Choose isomorphisms 
$f_i: A_i\isom\bS(a_i)_+$ with standard open annuli such that $f_i(x)$
approaches the Gauss point of $\B(1)$ as $x$ approaches $y$.  Let $X$
be the curve obtained from $Y$ by gluing an open ball $\B(1)_+$ onto each
$A_i$ via the inclusions $f_i: A_i\isom\bS(a_i)_+\subset\B(1)_+$.  Then
$X$ is a smooth, proper, connected $K$-analytic curve, and it is clear
from the construction that $(X,y)$ is star-shaped and that $Y$ is a simple
neighborhood of $y$ in $X$.\qed

\paragraph \label{par:compactifications}
A proper star-shaped curve $X$ and an inclusion $i:Y\isom U\subset X$ as in
Proposition~\ref{prop:compactify.star} is called a \emph{compactification}
of $Y$ (as a star-shaped curve).  Note that $X\setminus i(Y)$ is a
disjoint union of finitely many closed balls, one for each open annulus
connected component of $Y\setminus\{y\}$.

\paragraph \label{par:star.reduction}
Let $(Y,y)$ be a star-shaped curve.  The smooth, proper, connected
$k$-curve $C_y$ with function field $\td\sH(y)$ is called the
\emph{residue curve} of $Y$.  The tangent vectors in $T_y$ are naturally
in bijective correspondence with the connected components of
$Y\setminus\{y\}$.  We define a \emph{reduction map} 
$\red: Y\to C_y$ by sending $y$ to the generic point of $C_y$, and sending
every point in a connected component $B$ of $Y\setminus\{y\}$ to the
closed point of $C_y$ corresponding to the tangent vector determined by
$B$.  This sets up a one-to-one correspondence between the connected
components of $Y\setminus\{y\}$ and the closed points of $C_y$.

\begin{rem}
  \begin{enumerate}
  \item When $Y$ is proper, so $Y\cong\fX_K$ for a smooth, proper,
    connected formal curve over $\Spf(R)$, then $C_y\cong\fX_k$ and the
    reduction map $Y\to C_y$ coincides with the canonical reduction map
    $\fX_K\to\fX_k$.

  \item Let $D\subset Y(K)$ be a compatible divisor and let 
    $\Sigma_0 = \Sigma(Y,\{y\}\cup D)$.  Every edge $e$ of $\Sigma_0$ is
    contained in a unique connected component of $Y\setminus\{y\}$, and
    $\red(e)\in C_y(k)$ is the closed point corresponding to the tangent
    direction represented by $e$.
  \end{enumerate}
\end{rem}

\paragraph[Tame coverings]
We now study a class of morphisms of star-shaped curves 
analogous to~\parref{par:tame.covers.curves}. 
We begin with the following technical result.

\begin{lem} \label{lem:extract.roots}
  Let $A,A'$ be open annuli or punctured open balls and let
  $\phi: A'\to A$ be a finite morphism of degree $\delta$.
  Suppose that $\delta$ is prime to $\chr(k)$ if $\chr(k)>0$.  Fix an
  isomorphism $A\cong\bS(a)_+$ (where we allow $a=0$).  
  \begin{enumerate}
  \item There is an
    isomorphism $A'\cong\bS(a')_+$ such that the composition
    \begin{equation} \label{eq:extend.ball}
      \bS(a')_+ \cong A' \overset\phi\To A \cong \bS(a)_+ 
    \end{equation}
    is $t\mapsto t^\delta$.  
  \item There is an isomorphism $A'\cong\bS(a')_+$ such
    that~\eqref{eq:extend.ball} extends to a morphism 
    $\psi: \B(1)_+\to\B(1)_+$ with $\psi\inv(0)=0$.
  \item If \eqref{eq:extend.ball} extends to a morphism 
    $\psi: \B(1)_+\to\B(1)_+$ for a given isomorphism $A'\cong\bS(a')_+$,
    then the extension is unique.
  \end{enumerate}
\end{lem}

\pf Let $u$ be a parameter on $A'$, i.e.\ an isomorphism
$u:A'\isom\bS(a')_+\subset\B(1)_+$ with a standard open annulus.
By~Lemma~\ref{lem:maps.annuli}, 
$\phi$ restricts to an affine map on skeleta
$\Sigma(A')\to\Sigma(A)$ with degree $\delta$.
If $B = \bS(b,c)$ is a closed
sub-annulus of $\bS(a')_+$, then by~\cite[Lemme~2.2.1]{thuillier:thesis},
after potentially replacing $u$ by $u\inv$,
$\phi^*(t)$ has the form $\alpha u^\delta (1+g(u))$ on $B$, where
$\alpha\in R^\times$ and $|g|_{\sup} < 1$.  Since $\delta$ is not
divisible by $\chr(k)$ if $\chr(k)>0$, the Taylor expansion for 
$\sqrt[\delta]{1+g}$ has coefficients contained in $R$, hence converges to a
$\delta$th root of $1+g(u)$ on $A$.  Choosing a $\delta$th root of
$\alpha$ as well and letting
$u_B=u\,\sqrt[\delta]\alpha\sqrt[\delta]{1+g}$, we have that $u_B$ is a
parameter on $B$ such that $\phi^*(t) = u_B^\delta$.  Hence for each such
$B$ there are exactly $\delta$ choices of a parameter $u_B$ on $B$ such that
$\phi^*(t) = u_B^\delta$; choosing a compatible set of such parameters for
all $B$ yields a parameter on $A'$ satisfying~(1).

Part~(2) follows immediately from~(1).  A morphism from a $K$-analytic
space $X$ to $\B(1)_+$ is given by a unique analytic function $f$ on $X$
such that $|f(x)|<1$ for all $x\in X$, so~(3) follows from the fact that
the restriction homomorphism $\OO(\B(1)_+)\to\OO(\bS(a')_+)$ is injective.\qed

\paragraph \label{par:AutdA}
Let $A$ be an open annulus.  Let $\delta$ be a positive integer,
and assume that $\delta$ is prime to the characteristic of $k$ if
$\chr(k)>0$.  Let $A'$ be an open annulus and let 
$\phi:A'\to A$ be a finite morphism of degree $\delta$.  By
Lemma~\ref{lem:extract.roots}, the group
$\Aut_A(A')$ is isomorphic to $\Z/\delta\Z$.  Moreover, if $A''$ is
another open annulus and $\phi':A''\to A$ is a finite morphism
of degree $\delta$, then there exists an $A$-isomorphism
$\psi: A'\isom A''$; the induced isomorphism
$\Aut_A(A')\isom\Aut_A(A'')$ is independent of the choice of $\psi$ since
both groups are abelian.  Therefore the group
$\Aut_A(\delta) \coloneq \Aut_A(A')\cong\Z/\delta\Z$ is canonically
determined by $A$ and $\delta$.

\begin{defn}
  Let $(Y,y,D)$ be a punctured star-shaped curve with 
  skeleton $\Sigma_0 = \Sigma(Y,\{y\}\cup D)$
  and let $C_y$ be the residue curve of $Y$.  A \emph{tame covering} of
  $(Y,y,D)$ consists of a punctured star-shaped curve $(Y',y',D')$ and a
  finite morphism $\phi: Y'\to Y$ satisfying the following properties:
  \begin{enumerate}
  \item $\phi\inv(y)=\{y'\}$,
  \item $D' = \phi\inv(D)$, and
  \item if $C_{y'}$ denotes the residue curve of $Y'$ and 
    $\phi_{y'}:C_{y'}\to C_y$ is the morphism induced by $\phi$, then
    $\phi_{y'}$ is tamely ramified and is
    branched only over the points of $C_y$ corresponding to tangent
    directions at $y$ represented by edges in $\Sigma_0$.
  \end{enumerate}
\end{defn}

\begin{rem}
  The degree of $\phi$ is equal to the degree of $\phi_{y'}$.
\end{rem}

\begin{eg} \label{eg:curve.star.morphism}
  Let $\phi: (X',V',D')\to(X,V,D)$ be a tame covering of triangulated
  punctured curves.  Let $\Sigma = \Sigma(X,V\cup D)$ and
  $\Sigma' = \Sigma(X',V'\cup D') = \phi\inv(\Sigma)$.
  Let $x\in V$ be a finite vertex of $\Sigma$, let $e_1,\ldots,e_r$ be the
  finite edges of $\Sigma$ 
  adjacent to $x$, let $\Sigma_0 = \{x\}\cup\bigcup_{i=1}^r e_i^\circ$,
  let $Y = \tau\inv(\Sigma_0)$, and let $D_0 = D\cap Y$.  Then $(Y,x,D_0)$ is a
  punctured star-shaped curve.  Let 
  $\Sigma_0'$ be a connected component of $\phi\inv(\Sigma_0)$, let
  $x'\in\Sigma_0'$ be the unique inverse image of $x$, let
  $Y' = \tau\inv(\Sigma_0')$, and let $D_0' = D'\cap Y'$.  Then
  $(Y',x',D_0')$ is also a punctured star-shaped curve,
  and $\phi$ restricts to a finite morphism $\phi:Y'\to Y$.
  This is in fact a tame covering of punctured star-shaped curves by
  Remark~\ref{rem:tame.visible.ram} and Proposition~\ref{prop:branch.locus}.
\end{eg}

\medskip
\begin{prop} \label{prop:tame.cover}
  Let $\phi:(Y',y',D')\to(Y,y,D)$ be a degree-$\delta$ tame covering of
  punctured star-shaped curves.%
  \footnote{This $\delta$ need not be prime to $\chr(k)$.}
  Let $\Sigma_0 = \Sigma(Y,\{y\}\cup D)$, let
  $\Sigma_0' = \Sigma(Y',\{y'\}\cup D')$, let $C_y$ (resp.\ $C_{y'}$) be
  the residue curve of $Y$ (resp.\ $Y'$), and let 
  $\phi_{y'}:C_{y'}\to C_y$ be the induced morphism.
  \begin{enumerate}
  \item Let $B$ be a connected component of $Y\setminus\{y\}$ disjoint
    from $\Sigma_0$.  Then $\phi\inv$ is a disjoint union of $\delta$ open
    balls mapping isomorphically onto $B$.
  \item Let $B$ be a connected component of $Y\setminus\{y\}$ meeting
    $D$, and choose an isomorphism of $B$ with $\B(1)_+$ which identifies
    the unique point of $B\cap D$ with $0$.  Let $B'$ be a connected
    component of $\phi\inv(B)$.   Then $\phi|_{B'}:B'\to B$ is a finite
    morphism, the degree of $\phi|_{B'}$ is the ramification degree $\delta'$ of
    $\phi_{y'}$ at $\red(B')$, and there is an isomorphism 
    $B'\cong\B(1)_+$ sending the unique point of $B'\cap D'$ to $0$ such
    that the composition
    \[ \B(1)_+\cong B'\overset\phi\To B\cong\B(1)_+ \]
    is $t\mapsto t^{\delta'}$.
  \item Let $A$ be an open annulus connected component of
    $Y\setminus\{y\}$, and choose an isomorphism $A\cong\bS(a)_+$.  Let
    $A'$ be a connected component of $\phi\inv(A)$.  Then 
    $\phi|_{A'}:A'\to A$ is a finite morphism, the degree of 
    $\phi|_{A'}$ is the ramification degree $\delta'$ of 
    $\phi_{y'}$ at $\red(A')$, and there is an isomorphism of $A'$ with an
    open annulus $\bS(a')$ such that the composition
    \[ \bS(a')_+\cong A'\overset\phi\To A\cong\bS(a)_+ \]
    is $t\mapsto t^{\delta'}$.
  \item $\phi$ is \'etale over $Y\setminus D$.
  \item $\phi\inv(\Sigma_0) = \Sigma_0'$.
  \end{enumerate}
\end{prop}

\pf In the situation of~(1), let $\bar x = \red(B)$.  Since $\phi_{y'}$ is
not branched over $\bar x$, there are $\delta$ distinct points of
$C_{y'}$ mapping to $\bar x$; hence there are $\delta$ connected
components $B_1',\ldots,B_\delta'$ of $Y'\setminus\{y'\}$ mapping onto
$B$.  The restriction of $\phi$ to each $B_i'$ is finite of
degree $1$.

Next we prove~(2).  It is clear that $\phi|_{B'}:B'\to B$ is finite, hence
surjective; since $B'$ is a connected component of $Y'\setminus\{y'\}$
containing a point of $D'=\phi\inv(D)$, it is an open ball.
Let $x$ (resp.\ $x'$) be the unique point of $D$
(resp.\ $D'$) contained in $B$ (resp.\ $B'$) and let $e\subset\Sigma_0$ 
(resp.\ $e'\subset\Sigma_0'$)  be the edge adjacent to $x$ (resp.\ $x'$). 
Then $B'\setminus\{x'\}\to B\setminus\{x\}$ 
is a finite morphism of punctured open balls, so
by Lemma~\ref{lem:extend.end.ball}, Proposition~\ref{prop:skel.inverse.ball}, 
and Theorem~\ref{thm:harmonic.skeleta}(2),
the restriction of $\phi$ to $e'$ is affine morphism
$e'\to e$ of degree $\delta'$, and $\delta'$ is the
degree of $B'\to B$. So by Lemma~\ref{lem:extract.roots} there is an
isomorphism $B'\cong\B(1)_+$ as described in the statement of the
Theorem.  

In the situation of~(3), we claim that $A'$ is an open annulus.  Clearly
$\phi|_{A'}:A'\to A$ is finite, hence surjective.  If $A'$
is not an open annulus, then $A'\cong\B(1)_+$ is an open ball.  The morphism
$\B(1)_+\cong A'\to A\cong\bS(a)_+$ is given by a unit on $\B(1)_+$, which
has constant absolute value; this contradicts surjectivity, so $A'$ is in
fact an open annulus.  The proof now proceeds exactly as above.

Parts~(4) and~(5) follow immediately from parts~(1)--(3).\qed

\smallskip
The following proposition is reminiscent
of~\cite[Proposition~3.3.2]{saidi:wild_ramification}. 

\begin{cor} \label{cor:compactify.covers}
  Let $\phi:(Y',y',D_0')\to(Y,y,D_0)$ be a tame covering of punctured
  star-shaped curves and let 
  $i: (Y,y)\inject(X,x)$ be a compactification of $Y$.
  Let $D_1$ be the union of $i(D_0)$ with a choice of $K$-point from every
  connected component of $X\setminus i(Y)$.  
  Then there exists a compactification $i': (Y',y')\inject(X',x')$
  and a tame covering $\psi:(X',x',D_1')\to(X,x,D_1)$ such
  that $\psi\circ i' = i\circ\phi$.
\end{cor}

\pf We compactify $(Y',y')$ as in the proof of 
Proposition~\ref{prop:compactify.star}, gluing balls onto the annulus
connected components of $Y'\setminus\{y'\}$.  By
Proposition~\ref{prop:tame.cover}(3), if $A\cong\bS(a)_+$ is an open
annulus connected component of 
$Y\setminus\{y\}$ and $A'\subset Y'\setminus\{y'\}$ is a 
connected component mapping to $A$, then we can choose an isomorphism 
$A'\cong\bS(a')_+$ such that $\bS(a')_+\cong A'\to A\cong\bS(a)_+$ is of the form 
$t\mapsto t^\delta$; this map extends to a morphism
$\bB(1)_+\to\bB(1)_+$, and these maps glue to give a tame covering
$X'\to X$.\qed 

\paragraph[A local lifting theorem]
Let $(Y,y,D)$ be a punctured star-shaped curve with skeleton
$\Sigma_0 = \Sigma(Y,\{y\}\cup D)$ and residue curve $C_y$.
Let $C'$ be a smooth, proper, connected $k$-curve
and let $\bar\phi: C'\to C_y$ be a finite, tamely ramified morphism 
branched only over the points of $C_y$ corresponding to tangent directions
at $y$ represented by edges in $\Sigma_0$.  A 
\emph{lifting of $C'$ to a punctured star-shaped curve over $(Y,y,D)$} 
is the data of a punctured star-shaped
curve $(Y',y',D')$, a tame covering $\phi:(Y',y',D')\to(Y,y,D)$, and
an isomorphism of the 
residue curve $C_{y'}$ with $C'$ which identifies $\bar\phi$ with the
morphism $\phi_{y'}:C_{y'}\to C_y$ induced by $\phi$.  An \emph{isomorphism}
between two liftings $(Y',y')$ and $(Y'',y'')$ is a $Y$-isomorphism
$Y'\to Y''$  such that the induced morphism
$C_{y'}\isom C_{y''}$ respects the identifications 
$C_{y'}\cong C'$ and $C_{y''}\cong C'$.

\begin{thm} \label{thm:star.lifts}
  Let $(Y,y,D)$ be a punctured star-shaped curve with skeleton
  $\Sigma_0 = \Sigma(Y,\{y\}\cup D)$ and residue curve $C_y$, let
  $C'$ be a smooth, proper, connected $k$-curve,
  and let $\bar\phi: C'\to C_y$ be a finite, tamely ramified morphism 
  branched only over the points of $C_y$ corresponding to tangent directions
  at $y$ represented by edges in $\Sigma_0$.  Then there exists a lifting
  of $C'$ to a punctured star-shaped curve over $(Y,y,D)$, and
  this lifting is unique up to unique isomorphism. 
\end{thm}

\smallskip
Before giving the proof of this theorem, we need a technical lemma.
Let $\cX$ be a finitely presented, flat, separated $R$-scheme
and let $\fX$ be its $\varpi$-adic completion; this is an admissible formal
$R$-scheme by~\cite[Proposition~3.12]{bpr:trop_curves}.  There is a canonical open immersion
$i_{\cX}: \fX_K\to\cX_K^\an$ 
defined in~\cite[\S{}A.3]{conrad:irredcomps}
which is functorial in $\cX$ and respects the formation of fiber products,
and is an isomorphism when $\cX$ is proper over $R$.  (This fact is
implicitly contained in the statement of Lemma~\ref{lem:algebraization}.)

\begin{lem} \label{lem:finite.cartesian}
  Let $\cX,\cX'$ be finitely presented, flat, separated $R$-schemes and let
  $\fX,\fX'$ denote their $\varpi$-adic completions.  Let $\cX'\to\cX$ be
  a finite and flat morphism.  Then the square
  \begin{equation} \label{eq:compl.immersions}
    \xymatrix{
    {\fX'_K} \ar[r]^{i_{\cX'}} \ar[d]_f & {\cX_K'^\an} \ar[d]^g \\
    {\fX_K} \ar[r]_{i_{\cX}} & {\cX_K^\an}
  }\end{equation}
  is Cartesian.
\end{lem}

\pf The vertical arrows $f,g$ of~\eqref{eq:compl.immersions} are finite and the
horizontal arrows $i_{\cX'},i_{\cX}$ are open immersions.  Let 
$Y = \cX_K'^\an\times_{\cX_K^\an} \fX_K$, so $Y\to\cX_K'^\an$ is an open
immersion and $Y\to\fX_K$ is finite.  Let $h:\fX'_K\to Y$ be the canonical
morphism.  Then $h$ is an open immersion because its composition with
$Y\to\cX_K'^\an$ is an open immersion, and $h$ is finite because its
composition with $Y\to\fX_K$ is finite.  Therefore $h$ is an isomorphism
of $\fX'_K$ onto an open and closed subspace of $Y$.  It suffices to show
that $h$ is surjective, and since $Y(K)$ is dense in $Y$, we only need to
check that $Y(K)$ is in the image of $h$.  

Let $x\in\fX_K(K)$ and $y\in\cX_K'^\an(K)$ be points with the same image 
$z\in\cX_K^\an(K)$, so $(x,y)\in Y(K)$.  Choose an
open formal affine $\fU\subset\fX$ such that
$x\in\fU_K(K)$, and assume without loss of generality that there is an affine open 
$\cU = \Spec(A)\subset\cX$ such that $\fU = \Spf(\hat A)$, where $\hat A$
is the $\varpi$-adic completion of $A$.  Then $x$ corresponds to a
continuous $K$-homomorphism $\eta: \hat A\tensor_R K\to K$.  We have
$\eta(\hat A)\subset R$, so composing with the completion homomorphism 
$A\to\hat A$ we obtain an $R$-homomorphism $A\to R$.  Let
$\xi:\Spec(R)\to\cU\subset\cX$ denote the induced morphism.  The
corresponding morphism $\Spf(R)_K\to\fX_K$ is the point $x$.

Let $\cZ$ be the fiber product of $\cX'\to\cX$ with $\eta:\Spec(R)\to\cX$
and let $\fZ$ be its $\varpi$-adic completion.  Note that $\cZ$ is finite
and flat over $\Spec(R)$.  In particular, $\cZ$ is proper over $R$, so 
$i_\cZ:\fZ_K\to\cZ_K^\an$ is an isomorphism.  We have a commutative cube 
\[\xymatrix @=.15in{
  & {\cZ_K^\an} \ar[rr] \ar '[d] [dd] & & {\cX_K'^\an} \ar[dd] ^g \\
  {\fZ_K} \ar[ur]^\cong \ar[rr] \ar[dd] & & {\fX_K'} \ar[ur] \ar[dd] ^(.3)f & \\
  & {\Spec(K)^\an} \ar '[r] ^(.8)z [rr] & & {\cX_K^\an} \\
  {\Spf(R)_K} \ar[ur]^\cong \ar[rr]^x & & {\fX_K} \ar[ur] &
}\]
where the front and back faces are Cartesian and the diagonal arrows are
the canonical open immersions.  Since $g(y) = z$ and 
$\cZ_K^\an = f\inv(z)$, the point $y$ lifts to a point in
$\cZ_K^\an(K)$, which then lifts to a point in $\fZ_K(K)$ whose image in
$\fX'_K(K)$ maps to $x$ in $\fX_K(K)$ and $y\in\cX_K'^\an(K)$.\qed

\paragraph[Proof of Theorem~\ref{thm:star.lifts}] \label{par:algebraize}
We first prove the theorem when $Y$ is proper.  In
this case we may and do assume that
$Y$ is the analytic generic fiber of a smooth, proper, connected
formal curve $\fX$ over $\Spf(R)$.  Let $\cX\to\Spec(R)$ be the
algebraization of $\fX$ (Lemma~\ref{lem:algebraization}); this is a
smooth, proper relative curve of finite 
presentation and with connected fibers whose $\varpi$-adic completion is
isomorphic to $\fX$.  Note that $\cX_k = \fX_k = C_y$.
Let $X = \cX_K$, so $X^\an = \fX_K = Y$.  By the
valuative criterion of properness, every point $x$ of $D$ extends uniquely to
a section $\Spec(R)\to\cX$ which sends the closed point to the reduction
of $x$; hence the closure $\cD$ of $D$ in $\cX$ is a disjoint union of
sections.  Let $\cU = \cX\setminus\cD$.

The theory of the
\emph{tamely ramified \'etale fundamental group} $\pi_1^t$ of a morphism
of schemes with a relative normal crossings divisor is developed
in~\cite[Expos\'e~XIII]{SGA1}.  The finite-index 
subgroups of $\pi_1^t$ classify so-called tamely ramified \'etale covers.  
The subscheme $\cD\subset\cX$ is a relative normal crossings divisor
relative to $\Spec(R)$, so the proof 
of Corollaire~2.12 of {\it loc.\ cit.}\ shows that the specialization
homomorphism $\pi_1(\cU_K)\to\pi_1^t(\cU_k)$ is surjective (we suppress
the base points).  What this
means concretely is that every tamely ramified cover $\bar U\p\to\cU_k$ over
$\cX_k$ relative to $\cD_k$ extends to a finite \'etale
morphism $\cU'\to\cU$, unique up to unique isomorphism, and the generic
fiber of $\cU'$ is connected. 

Let $\bar D\p = \bar\phi\inv(\cD_k)$ and
let $\bar U\p = C'\setminus\bar D\p$, so $\bar U\p\to\cU_k$ is
a tamely ramified cover of $\cU_k$ over $\cX_k$ relative to $\cD_k$.  Let 
$\phi: \cU'\to\cU$ be the unique \'etale covering extending
$\bar\phi|_{\bar U\p}$;
this is equipped with an isomorphism $\cU'_k\cong\bar U\p$ identifying
$\bar\phi|_{\bar U\p}$ with $\phi_k$.  
Let $U = \cU_K = X\setminus D$, let
$U' = \cU'_K$, let $X'$ be the smooth compactification of $U'$, 
and let $\phi_K: X'\to X$ denote the unique morphism extending 
$\phi_K: U'\to U$.  
We will show that there is a simple neighborhood $W$ of $y$ in $X^\an$
such that $\phi_K\inv(W)\to W$ is a tame covering, and conclude that a
lifting exists using Corollary~\ref{cor:compactify.covers}. 

\subparagraph \label{par:what.is.xp}
First we claim that $\phi_K\inv(y)$ consists of a unique point 
$y'\in X'^\an$.  Letting $\fU$ and $\fU'$ denote the $\varpi$-adic
completions of $\cU$ and $\cU'$, respectively, we have a Cartesian square
\[\xymatrix{
  {\fU'_K} \ar[r]^{i_{\cU'}} \ar[d] & {U'^\an} \ar[d] \\
  {\fU_K} \ar[r]_{i_\cU} & {U^\an}
}\]
by Lemma~\ref{lem:finite.cartesian}.  In other words,
$\fU_K'$ is the inverse image of $\fU_K$ under
$\phi_K:U'^\an\to U^\an$.  Since the reduction map
$\fU_K\to\fU_k=\cU_k$ is surjective and $y$ is the unique point of $X^\an$
reducing to the generic point of $\fU_k\subset\fX_k$, we have $y\in\fU_K$,
so it suffices to show that there is a unique point $y'\in\fU'_K$ mapping
to $y$.  It follows easily from the functoriality of the reduction map
that the only point $y'$ mapping to $y$ is the unique point of $\fU'_K$
reducing to the generic point of $\fU'_k$.  In particular, the residue
curve of $X'^\an$ at $y'$ is identified with $C'$, which is the smooth
completion of $\fU'_k = \bar U\p$.

\subparagraph \label{par:single.edge}
Choose a strongly semistable vertex set $V'$ of $X'$ containing $y'$, let 
$\Sigma' = \Sigma(X',V')$, and let $\tau: X'^\an\to\Sigma'$ be the
retraction.  Let $B\subset X^\an$ be the formal fiber of a point 
$\bar x\in\fU_k(k)$; equivalently, $B$ is a connected component of
$X^\an\setminus\{y\}$ not meeting $D$.
Let $B'$ be a connected component of $X'^\an\setminus\{y'\}$ contained in
$\phi_K\inv(B)$.  We have
$B\subset\fU_K$; therefore $B'\subset\fU_K'$, so $B'$ is a connected
component of $\fU_K'\setminus\{y'\}$, hence $B'$ is the formal fiber of a
point $\bar x\p\in\fU_k'(k)$, so $B'\cong\B(1)_+$.  It follows that
$V'\setminus(V'\cap B')$ is again a semistable vertex set, so
we may and do assume that $\Sigma'\cap B' = \emptyset$ for all such $B'$.

Let $e'\subset\Sigma'$ be an edge adjacent to $y'$ and let $B$ be the
connected component of $X^\an\setminus\{y\}$ containing 
$\phi_K(e'^\circ)$.  By the above,
$B$ is the formal fiber of a point $\bar x\in\cD_k$.  
Fix an isomorphism $B\cong\B(1)_+$ taking the unique point in $D\cap B$
to $0$.  Let $v'\in T_{y'}$ be the tangent
vector in the direction of $e'$ and let $v = d\phi_K(y')(v')\in T_y$.  The
point in $\fX_k(k)$ corresponding to $v$ is $\bar x$; let
$\bar x\p\in\bar D\p\subset C'(k)$ be the point corresponding to
$v'$, so $\bar\phi(\bar x\p) = \bar x$.
By Theorem~\ref{thm:harmonic.skeleta}, the ramification degree $\delta$ of 
$\bar\phi$ at $\bar x\p$ is equal to $d_{v'}\phi_K(x')$.

Let $B'$ be the connected component of $X'^\an\setminus\{y'\}$ containing
$e'^\circ$. Viewing $\phi_K|_{B'}$ as a morphism
$B'\to\B(1)_+\subset\A^{1,\an}_K$ via our 
chosen isomorphism $B\cong\B(1)_+$, the map 
$x'\mapsto\log|\phi_K(x')|$ is a piecewise affine function on $e'^\circ$
which changes slope on the retractions of the zeros of $\phi_K$.
Shrinking $e'$ if necessary, we may and do assume that 
$\log|\phi_K|$ is affine on $e'^\circ$; its slope has absolute value
$\delta$, and $\log|\phi_K(x')|\to 0$ as $x'\to y'$.
By~\parref{par:incuced-morphisms-mc}
the restriction of $\phi_K$ to the open annulus
$\tau\inv(e'^\circ)$ is a finite morphism of degree $\delta$ onto an open
annulus $\bS(a)_+\subset\B(1)_+$, where $\val(a)$ is $\delta$ times the
length of $e'^\circ$.

\subparagraph 
Enlarging $\Sigma'$ if necessary, we may and do assume that every tangent
direction $v'\in T_{y'}$ corresponding to a point in $\bar D\p$ is
represented by an edge in $\Sigma'$.  Let $e_1',\ldots,e_r'$ be the edges
of $\Sigma'$ adjacent to $y'$, for $i=1,\ldots,r$ let $B_i$ be the
connected component of $X^\an\setminus\{y\}$ containing
$\phi_K(e_i'^\circ)$, and choose isomorphisms $B_i\isom\B(1)_+$ sending
the unique point of $D\cap B_i$ to $0$.  (The balls
$B_1,\ldots,B_r$ are not necessarily distinct; we mean that one should choose a single
isomorphism for each distinct ball.) Let
$\bar x\p_i\in\bar D\p$ be the point corresponding to the tangent
direction at $y'$ in the direction of $e_i'$, and let $\delta_i$ be the
ramification degree of $\bar\phi$ at $\bar x\p_i$.  Applying the procedure
of~\parref{par:single.edge} for each $e_i'$, and shrinking if necessary, we
may and do assume that for every $i$, $\phi_K$ induces a degree-$\delta_i$
morphism of $\tau\inv(e_i'^\circ)$ onto an open annulus
$\bS(a)_+\subset\B(1)_+\cong B_i$, with $a$ independent of $i$.

Let $W\subset X^\an$ be the union of $\tau\inv(y)$ with the annuli
$\bS(a)_+\subset B_i$, so $W$ is obtained from $X^\an$ by removing a closed
ball around each point of $D$.  This is a simple neighborhood of $y$ in
$X^\an$, hence is a star-shaped curve.  Let 
$\Sigma_0' = \{y'\}\cup\bigcup_{i=1}^r e_i'^\circ$ and let 
$W' = \tau\inv(\Sigma_0')$, so $W'$ is a simple neighborhood of $y'$ in 
$X'^\an$ and is hence a star-shaped curve.  We claim that
$W' = \phi_K\inv(W)$.  Clearly $\phi_K(W')\subset W$, so it suffices to
show that the fibers of $\phi_K$ have length equal to $\deg(\phi_K)$.
This is certainly the case for $\phi_K\inv(y) = \{y'\}$.  If 
$B\subset X^\an$ is the formal fiber of a point in $\fU_k(k)$, then
as in the proof of Proposition~\ref{prop:tame.cover}(1),
$\phi_K\inv(B)\subset W'$ is a disjoint union of $\deg(\phi_K)$ open balls
mapping isomorphically onto $B$.
For each $i$ the inverse image of 
$\bS(a)_+\subset\B(1)_+\cong B_i$ in $W'$ is a disjoint union of annuli,
one for each point in the fiber of $\bar\phi$ containing $\bar x\p_i$, and
the degree of $\phi_K$ restricted to each annulus is the ramification
degree of $\bar\phi$ at the corresponding point.  The sum of the
ramification degrees of $\bar\phi$ at the points in any fiber is equal to 
$\deg(\bar\phi) = \deg(\phi)$, which
proves the claim.
Therefore $\phi_K$ induces a finite morphism $W'\to W$, which is a tame
covering of $(W,y)$ (really of $(W,y,\emptyset)$).

By construction, $(Y,y)$ is a compactification of $(W,y)$, so by
Corollary~\ref{cor:compactify.covers}, the tame covering $W'\to W$ lifts
to a tame covering $Y'\to Y$ relative to $D$.  It is clear that
$(Y',y',\phi_K\inv(D))$ is a lifting of $C'$ to a punctured star-shaped
curve over $(Y,y,D)$.  (One can show that in fact $Y'\cong X'^\an$,
although this is not clear \emph{a priori}.) 

\subparagraph 
It remains to prove (still in the case when $Y$ is proper) that liftings
are unique up to unique isomorphism.  Let
$\phi_K:(Y',y',D')\to(Y,y,D)$ and $\phi'_K:(Y'',y'',D'')\to(Y,y,D)$ be two
liftings of $C'$ to punctured star-shaped curves over $(Y,y,D)$.  Then $Y'$
and $Y''$ are also proper, hence 
we may and do assume that they are the analytic generic fibers of smooth,
proper, connected formal curves $\fX'$ and $\fX''$ as
in~\parref{par:algebraize}.  Since $\phi_K(y') = \phi'_K(y'') = y$, there are
unique finite morphisms $\phi:\fX'\to\fX$ and $\phi':\fX''\to\fX$ extending
$\phi$ and $\phi'$, respectively, and $\phi_k$ and $\phi'_k$ are
identified with $\bar\phi$ under the isomorphisms
$C_{y'}=\fX'_k\cong C'$ and 
$C_{y''}=\fX''_k\cong C'$.  
Let $\cX'$ and $\cX''$ be the algebraizations of $\fX'$ and $\fX''$,
respectively, and let $\cU' = \cX'\setminus\phi\inv(\cD)$ and
$\cU'' = \cX''\setminus{\phi'}\inv(\cD)$.  Then $\cU',\cU''$ are finite
\'etale coverings of $\cU$ lifting $\bar U\p\to\cU_k$, so there is a
unique isomorphism $\psi: \cU'\isom\cU''$ over $\cU$.  Since $Y'$ (resp.\
$Y''$) is the analytification of the smooth completion of $\cU'_K$ (resp.\
$\cU''_K$), $\psi_K$ extends uniquely to a $Y$-isomorphism 
$\psi_K:Y'\isom Y''$; we have $\psi_K(y') = y''$, and 
$\psi_{y'}: C_{y'}\isom C_{y''}$ is a $C_y$-isomorphism because 
$\cU_k'\isom\cU_k''$ is a $\cU_k$-isomorphism.  Hence $\psi_K$ is an
isomorphism of liftings.

As for uniqueness, suppose that $\psi_K:(Y',y')\to(Y',y')$ is an
automorphism of $Y'$ as a lifting of $C'$.  Then $\psi_K$ extends uniquely
to an $\cX$-automorphism $\psi:\cX'\isom\cX'$ that is the identity on
$\cX'_k$, which restricts to a $\cU$-automorphism of $\cU'$ that is the identity on
$\cU'_k$.  It follows from the uniqueness of $\cU'$ up to unique isomorphism
that $\psi_K$ is the identity when restricted to 
$\cU'_K$, so since $Y'$ is the analytification of the smooth completion of
$\cU'_K$, the automorphism $\psi_K$ is the identity.  This concludes the
proof for proper $Y$.

\subparagraph 
Now suppose that $(Y,y)$ is not proper.  Let $Y\inject X$ be a
compactification~\parref{par:compactifications} of $Y$; we will identify
$Y$ with its image in $X$.  Let $D_1$ be the union of $D$ with a choice
of one $K$-point from every connected component of $X\setminus Y$.  Let
$\phi:(X',y',D_1')\to(X,y,D_1)$ be a lifting of $C'$ to a punctured
star-shaped curve over $(X,y,D_1)$ and let $Y' = \phi\inv(Y)$.  Then 
$(Y',y',\phi\inv(D))$ is a lifting of $C'$ to a punctured star-shaped
curve over $(Y,y,D)$.

Let $\psi: Y'\to Y'$ be an automorphism of $Y'$ as a lifting of $C'$.
Since $\psi$ induces the identity map $C_{y'}\to C_{y'}$ on residue
curves, $\psi$ takes each connected component of $Y'\setminus\{y'\}$ to
itself.  Recall that $X'\setminus Y'$ consists of a disjoint union of
closed balls around the points of $D'_1\setminus D'$,
where $D_1' = \phi\inv(D_1)$.  Let 
$x'\in D_1'\setminus D'$, let $B'$ be the
connected component of $X'\setminus\{y'\}$ containing $x'$, and let 
$B = \phi(B')\subset X$.  By Proposition~\ref{prop:tame.cover}(2), we can 
choose isomorphisms $B\cong\B(1)_+$ and $B'\cong\B(1)_+$ sending $x'$ and
$\phi(x')$ to $0$, such that the composition
\[ \B(1)_+\cong B'\overset\phi\To B\cong\B(1)_+ \]
is of the form $t\mapsto t^\delta$.
Let $A' = Y'\cap B'$ and let $A = Y\cap B = \phi(A')$.  The isomorphism
$B\cong\B(1)_+$ (resp.\ $B'\cong\B(1)_+$) identifies
$A$ (resp.\ $A'$) with an open annulus
$\bS(a)_+$ (resp.\ $\bS(a')_+$) in $\B(1)_+$.  
Since $\psi|_{A'}:A'\to A'$ is an $A$-morphism, the composition
$\bS(a')_+\cong A'\to A'\cong\bS(a')_+$ is of the form
$t\mapsto\zeta t$, where $\zeta\in R^\times$ is a $\delta$th root of
unity.  Therefore $\psi|_{A'}$ extends uniquely to a $B$-morphism
$\psi|_{B'}:B'\to B'$ fixing $0$.  Gluing these morphisms together,
we obtain an $X$-morphism $\psi:X'\to X'$ extending $\psi: Y'\to Y'$.  By
construction this is an automorphism of $X'$ as a lifting of $C'$ to $X$,
which is thus the identity.  Therefore $\psi:Y'\to Y'$ is the identity.

Let $\phi':(Y'',y'',D'')\to(Y,y,D)$ be another lifting of $C'$ to a
punctured star-shaped curve over $(Y,y,D)$.  By Corollary~\ref{cor:compactify.covers},
there exists a compactification $Y''\inject X''$ of $Y''$ and an extension
of $\phi'$ to a tame covering $\phi':
(X'',y'',\phi'{}\inv(D_1))\to(X,y,D_1)$.  This is another lifting of $C'$ to
a punctured star-shaped curve over $(X,y,D_1)$, so there is an
isomorphism $\psi:(X',y')\isom(X'',y'')$ of liftings.  Since 
$Y' = \phi\inv(Y)$ and $Y'' = {\phi'}\inv(Y)$, the isomorphism $\psi$
restricts to an isomorphism $(Y',y')\isom(Y'',y'')$ of liftings.\qed

\begin{cor} \label{cor:auts.lift}
  With the notation in Theorem~\ref{thm:star.lifts}, let 
  $(Y',y',D')\to(Y,y,D)$ be a lifting of $C'$ to a punctured star-shaped
  curve over $(Y,y,D)$.  Then the natural homomorphism
  \[ \Aut_Y(Y') \To \Aut_{C_y}(C') \]
  is bijective.  If $(Y'',y'',D'')\to(Y,y,D)$ is a second lifting of $C'$
  to a punctured star-shaped curve over $(Y,y,D)$ then the natural map
  \[ \Isom_Y(Y'', Y') \To \Isom_{C_y}(C_{y''}, C_{y'}) \]
  is bijective.
\end{cor}

\section{Classification of liftings of harmonic morphisms of metrized complexes}
\label{sec:global.lifting}

Fix a triangulated punctured curve $(X,V\cup D)$ with skeleton
$\Sigma = \Sigma(X,V\cup D)$ and let
$\tau:X^\an\to\Sigma$ be the canonical retraction.  Throughout this
section, we \textbf{assume that $\Sigma$ has no loop edges}.
Let $\phi:\Sigma'\to\Sigma$ be a tame covering of metrized complexes of
curves.  
A \emph{lifting of $\Sigma'$ to a tame covering of $(X,V\cup D)$}
is a tame covering of triangulated punctured curves
$\phi:(X',V'\cup D')\to(X,V\cup D)$
(see~\parref{par:tame.covers.curves}) equipped with a $\Sigma$-isomorphism
$\phi\inv(\Sigma)\cong\Sigma'$ of metrized complexes of curves.  
Since $V' = \phi\inv(V)$ and $D' = \phi\inv(D)$, we will often denote a
lifting simply by $X'$.  Let $\phi:X'\to X$ and $\phi':X''\to X$ be two
liftings of $\Sigma'$ and let $\psi: X'\isom X''$ be an $X$-isomorphism of
curves.  Then $\psi$ restricts to a $\Sigma$-automorphism
\[ \psi|_{\Sigma'}~:~ 
\Sigma' \cong \phi\inv(\Sigma)\isom\phi'{}\inv(\Sigma) \cong \Sigma'. \]
We will consider liftings up to $X$-isomorphism preserving $\Sigma'$,
i.e.\ such that $\psi_{\Sigma'}$ is the identity, and we will also
consider liftings up to isomorphism as curves over $X$.

\paragraph
For every finite vertex $x\in V$, let $\Sigma(x)$ be the connected
component of $x$ in $\{x\}\cup(\Sigma\setminus V)$, let 
$Y(x) = \tau\inv(\Sigma(x))$, and let $D(x) = D\cap Y(x)$, as in
Example~\ref{eg:curve.star.morphism}.  Then $(Y(x),x,D(x))$ is a punctured
star-shaped curve.
By Proposition~\ref{prop:tame.visible.ram}, for every finite vertex
$x'\in V(\Sigma')$ lying above $x$ the morphism $\phi_{x'}: C_{x'}\to C_x$
is finite, tamely ramified, and branched only over the points of
$C_x$ corresponding to tangent directions at $x$ represented by edges of
$\Sigma(x)$.
Let $\psi(x'):(Y'(x'),x',D(x'))\to(Y(x),x,D(x))$ be the unique lifting of 
$C_{x'}$ to a punctured star-shaped curve over $(Y(x),x,D(x))$
provided by Theorem~\ref{thm:star.lifts} and let 
$\Sigma'(x') = \psi(x')\inv(\Sigma(x)) = \Sigma(Y',\{x'\}\cup D'(x'))$.
Then $\Sigma'(x')$ is canonically identified with the connected component
of $x'$ in $\{x'\}\cup(\Sigma'\setminus V_f(\Sigma'))$ in such a way that
for every edge $e'$ of $\Sigma'(x')$, the point 
$\red_{x'}(e')\in C_{x'}(k)$ is identified with the point 
$\red(e')$ defined in~\parref{par:star.reduction}.  This induces an
identification of $D'(x')$ with $\phi\inv(D)\cap\Sigma'(x')$.  
Let $\tau_{x'}$ be the canonical retraction $Y'(x')\to\Sigma'(x')$ defined
in~\parref{par:star.shaped.basic}.

\paragraph
Let $e'\in E_f(\Sigma')$ and $e = \phi(e')$.  Choose $a,a'\in K^\times$
with $\val(a) = \ell(e)$ and $\val(a') = \ell(e')$.  Let $x',y'$ be the
endpoints of $e'$.
By Proposition~\ref{prop:tame.cover}(3), we can choose isomorphisms
$\tau\inv(e^\circ)\cong\bS(a)_+$ and 
$\tau_{x'}\inv(e'^\circ)\cong\bS(a')_+$, 
$\tau_{y'}\inv(e'^\circ)\cong\bS(a')_+$ in such a way that the
finite morphisms $\tau_{x'}\inv(e'^\circ)\to\tau\inv(e^\circ)$ and
$\tau_{y'}\inv(e'^\circ)\to\tau\inv(e^\circ)$ are given by
$t\mapsto t^{d_{e'}(\phi)}$, where $d_{e'}(\phi)$ is the degree of the edge
map $e'\to e$.  In particular, there exists a
$\tau\inv(e^\circ)$-isomorphism 
$\tau_{x'}\inv(e'^\circ) \isom \tau_{y'}\inv(e'^\circ)$.  
As explained in~\parref{par:AutdA}, there are canonical identifications
\[ \Aut_{\tau\inv(e^\circ)}(d_{e'}(\phi)) \coloneq
\Aut_{\tau\inv(e^\circ)}(\tau_{x'}\inv(e'^\circ)) = 
\Aut_{\tau\inv(e^\circ)}(\tau_{y'}\inv(e'^\circ)) \cong
\Z/d_{e'}(\phi)\Z; \]
hence the set of isomorphisms
$\tau_{x'}\inv(e'^\circ) \isom \tau_{y'}\inv(e'^\circ)$
is a principal homogeneous space under the pre- or post-composition action of
$\Aut_{\tau\inv(e^\circ)}(d_{e'}(\phi))$.

Let $E_f^\pm(\Sigma')$ denote the set of oriented finite edges of
$\Sigma'$, and for $e'\in E_f^\pm(\Sigma')$ let $\bar e\p$ denote the same
edge with the opposite orientation.  Let $\cG(\Sigma',X)$ denote the set
of tuples $(\Theta_{e'})_{e'\in E_f^\pm(\Sigma')}$ of isomorphisms
\begin{equation} \label{eq:gluing.Theta}
  \Theta_{e'} ~:~ \tau_{x'}\inv(e'^\circ) \isom
  \tau_{y'}\inv(e'^\circ) 
\end{equation}
such that $\Theta_{\bar e\p} = \Theta_{e'}\inv$,
where $e' = \overrightarrow{x'y'}$.  We call $\cG(\Sigma',X)$ the set of
\emph{gluing data} for a lifting of $\Sigma'$ to a tame covering of
$(X,V\cup D)$, and we emphasize that $\cG(\Sigma',X)$ is \emph{nonempty}.

\paragraph 
Let $\alpha\in\Aut_\Sigma(\Sigma')$, so $\alpha$ is a degree-$1$ finite
harmonic morphism $\Sigma'\isom\Sigma'$ preserving $\Sigma'\to\Sigma$.
Let $x'\in V_f(\Sigma')$ and let $x'' = \alpha(x')$ and 
$x = \phi(x') = \phi(x'')$.  Part of the data of $\alpha$ is a
$C_x$-isomorphism $\alpha_{x'}: C_{x'}\isom C_{x''}$.  By
Corollary~\ref{cor:auts.lift}
there is a unique lift of $\alpha_{x'}$ to a $Y(x)$-isomorphism
$Y'(x')\isom Y'(x'')$ of punctured star-shaped curves inducing the
isomorphism $C_{x'}\isom C_{x''}$ on residue curves.  Let 
$e'\in E_f(\Sigma')$ be an edge adjacent to $x'$, let 
$e'' = \alpha(e')$, and let $e = \phi(e') = \phi(e'')$.  Then
$\alpha_{x'}$ restricts to an isomorphism
\[ \alpha_{x'} ~:~ \tau_{x'}\inv(e'^\circ) \isom
\tau_{x''}\inv(e''^\circ). \]

Define the \emph{conjugation action} of $\Aut_\Sigma(\Sigma')$ on
$\cG(\Sigma',X)$ by the rule
\begin{equation} \label{eq:conj.action}
  \alpha\cdot(\Theta_{e'})_{e'\in E_f^\pm(\Sigma')} =
  (\alpha_{y'}\inv \circ \Theta_{\alpha(e')} \circ 
  \alpha_{x'})_{e'\in E_f^\pm(\Sigma')} 
\end{equation}
where $e' = \overrightarrow{x'y'}$.

\begin{thm}[Classification of lifts of harmonic morphisms] \label{thm:lifting2}
  Let $(X,V\cup D)$ be a triangulated punctured curve with skeleton
  $\Sigma = \Sigma(X,V\cup D)$.  Assume that $\Sigma$ has no loop edges.
  Let $\phi: \Sigma'\to\Sigma$ be a tame covering of metrized compexes of
  curves. 
  \begin{enumerate}
  \item There is a canonical bijection between the set of gluing data
    $\cG(\Sigma',X)$ and the set of liftings of $\Sigma'$ to a tame
    covering of $(X,V\cup D)$, up to $X$-isomorphism preserving
    $\Sigma'$.  In particular, there exists a lifting of $\Sigma'$.
    Any such lifting has no nontrivial automorphisms which preserve
    $\Sigma'$.

  \item Two tuples of gluing data determine $X$-isomorphic curves if and
    only if they are in the same orbit under the conjugation
    action~\eqref{eq:conj.action}. The stabilizer in
    $\Aut_\Sigma(\Sigma')$ of an element of $\cG(\Sigma',X)$ is
    canonically isomorphic to the $X$-automorphism group of the associated
    curve. 
  \end{enumerate}
\end{thm}

\pf Given $(\Theta_{e'})\in\cG(\Sigma',X)$ one can glue the local lifts
$\{Y(x')\}_{x'\in V_f(\Sigma')}$ via the
isomorphisms \eqref{eq:gluing.Theta} to obtain an analytic space which one
easily verifies is smooth and proper, hence arises as the analytification
of an algebraic curve $X'$.  Moreover the morphisms $Y'(x')\to Y(x)$ glue
to give a morphism $X'^\an\to X^\an$, which is the analytification of a
morphism $\phi: X'\to X$.  By construction, if 
$V' = \phi\inv(V)$ and $D' = \phi\inv(D)$, then
$(X',V'\cup D')\to(X,V\cup D)$ is a lifting of $\Sigma'$ to a tame
covering of $(X,V\cup D)$.

Now let $\phi:(X',V'\cup D')\to(X,V\cup D)$ be any lifting of $\Sigma'$ to
a tame covering of $(X,V\cup D)$.  As explained in
Example~\ref{eg:curve.star.morphism}, for every $x\in V$ the inverse image
$\phi\inv(Y(x))$ is a disjoint union of tame covers of the punctured star-shaped
curve $(Y(x),x,D(x))$, one for each $x'\in\phi\inv(x)$.  By
Theorem~\ref{thm:star.lifts}, we have  canonical identification 
$\phi\inv(Y(x)) = \Djunion_{x'\mapsto x} Y(x')$.  For 
$x',y'\in V'$ we have 
$Y(x')\cap Y(y') = \Djunion_{e'\ni x',y'} \tau\inv(e'^\circ)$, where
$\tau: X'^\an\to\Sigma'$ is the retraction.  Hence $X'^\an$ is obtained by
pasting the local liftings $Y(x')$ via a choice of
isomorphisms $(\Theta_{e'})\in\cG(\Sigma',X)$.

Let $\psi:X'\isom X'$ be an $X$-automorphism preserving $\Sigma'$.  Since
$\psi$ is the identity on the set $\Sigma'\subset X'^\an$, for 
$x'\in V'$ we have $\psi(Y(x')) = Y(x')$.  We also have
that for $x'\in V'$ the map of residue curves 
$\psi_{x'}:C_{x'}\to C_{x'}$ is the identity, so by
Corollary~\ref{cor:auts.lift}, $\psi$ restricts to the identity morphism
$Y(x')\to Y(x')$.  Since $X'^\an = \bigcup_{x'\in V'} Y(x')$ we have that
$\psi$ is the identity.  It follows from this that the tuple
$(\Theta_{e'})$ can be recovered from the class of $X'$ modulo
$X$-isomorphisms preserving $\Sigma'$, so different gluing data give rise
to curves which are not equivalent under such isomorphisms.  This proves~(1).

Let $\phi:(X',V'\cup D')\to(X,V\cup D)$ and 
$\phi'':(X'',V''\cup D'')\to(X,V\cup D)$ be the liftings of $\Sigma'$
associated to the tuples of gluing data
$(\Theta'_{e'}),(\Theta''_{e'})\in\cG(\Sigma',X)$, respectively.
Suppose that there exists $\alpha\in\Aut_\Sigma(\Sigma')$ such that
$\alpha\cdot(\Theta'_{e'}) = (\Theta''_{e'})$.  Then for all
$e' = \overrightarrow{x'y'}\in E_f^\pm(\Sigma')$ we have
\[ \Theta''_{\alpha(e')} \circ \alpha_{x'} = \alpha_{y'} \circ \Theta'_{e'}, \]
so the $Y(\phi(x'))$-isomorphisms $\alpha_{x'}:Y(x')\isom Y(\alpha(x'))$
glue to give an $X$-isomorphism $\alpha: X'\isom X''$.  Conversely, the
restriction of an $X$-isomorphism $\alpha:X'\isom X''$ to $\Sigma'$ is a
$\Sigma$-automorphism of $\Sigma'$.  It is easy to see that these are
inverse constructions.  Taking $X' = X''$ we have an injective homomorphism
$\Aut_X(X')\inject\Aut_\Sigma(\Sigma')$; it follows formally from the
above considerations that its image is the stabilizer of
$(\Theta'_{e'})$.\qed 

\begin{rem}\label{rem:auto}
  Let $X'$ be a lifting of $\Sigma'$ to a
  tame covering of $(X,V\cup D)$.
  It follows from Theorem~\ref{thm:lifting2}(1) that the natural
  homomorphism
  \[ \Aut_X(X') \To \Aut_\Sigma(\Sigma') \]
  is injective.  It is not in general surjective: see
  Example~\ref{eg:2.isogeny} {(where all the stabilizer groups of elements of $\cG(\Sigma',X)$ 
  are proper subgroups of $\Aut_\Sigma(\Sigma')$).} 
\end{rem}

\medskip

\begin{rem}   \label{rem:wild.lifting}
  It is worth mentioning that the question of lifting (and classification of
  all possible liftings) in the wildly ramified case is more subtle and
  cannot be guaranteed in general. Some lifting results in the wildly
  ramified case are known: e.g., Liu proves in~\cite[Proposition
  5.4]{Liu03} that any finite surjective generically \'etale admissible
  cover of proper semistable curves over an algebraically closed field $k$
  lifts to a finite morphism of smooth proper curves over $K$. However,
  the same statement for metrized complexes cannot be true in general:
  consider the metrized complex $\cC$ consisting of a single finite vertex
  $v$ and an infinite vertex $u$ attached to $v$ by an infinite edge $e$
  with $C_v \cong \P^1_k$ and $\red_v(e) = \infty$. Consider the degree
  $p$ generically \'etale morphism of metrized complexes $\cC \to \cC$
  which restricts to a degree $p$ Artin-Schreier cover $\phi: \P^1_k \to
  \P^1_k$, \'etale over $\A^1_k\to \A^1_k$ and ramified at $\infty$ (with
  ramification index $p$), and which is linear of slope $p$ on the
  infinite edge $e$. This map cannot be lifted to a degree $p$ cover
  $\P^1_K \to \P^1_K$ in characteristic zero, otherwise, we would obtain a
  (connected) \'etale degree $p$ cover of $\A^1_K$ which is impossible
  since $\A^1_K$ is simply connected.
\end{rem}

\smallskip

We now state a variant of Theorem~\ref{thm:lifting2} 
in which we allow loop edges.

\begin{thm} \label{thm:lifting}
  Let $(X,V \cup D)$ be a triangulated punctured $K$-curve,
  let $\Sigma = \Sigma(X, V\cup D)$, and let 
  $\bar{\phi}: \Sigma'\to\Sigma$ be a tame covering of $\Lambda$-metrized complexes of
  $k$-curves.  Then there exists a tame covering $\phi: (X',V'\cup D')\to(X,V\cup D)$ of
  triangulated punctured $K$-curves lifting $\phi$.
\end{thm}

\pf 
The only issue is that in Theorem~\ref{thm:lifting}, we merely require
that $V$ is a semistable vertex set whereas whereas in
Theorem~\ref{thm:lifting2} we require it to be strongly semistable.
However, we can modify $\Sigma$ by inserting a (valence $2$) vertex $v$
along each loop edge, with $C_v \cong \PP^1_k$, and with the marked points
on $C_v$ being $0$ and $\infty$.  Call the resulting metrized complex
$\tilde{\Sigma}$, and let $\bar{V}_0$ denote the set of vertices which
have been added to the vertex set $\bar{V}$ for $\Sigma$.  We construct a
new metrized complex $\tilde{\Sigma}'$ from $\Sigma'$ by adding
$\bar{V}\p_0 := \bar{\phi}^{-1}(\bar{V}_0)$ to the vertex set $\bar{V}\p$
for $\Sigma'$ and letting $C'_{v'} \cong \PP^1_k$ for $v' \in
\bar{V}\p_0$, with the marked points on $C'_{v'}$ being $0$ and $\infty$.
We can extend the tame covering $\bar{\phi} : \Sigma' \to \Sigma$ to a
tame covering $\tilde{\Sigma}' \to \tilde{\Sigma}$ by letting the map from
$C'_{v'}$ to $C_{\bar\phi(v')}$, for $v' \in \bar{V}\p_0$, be $z \mapsto
z^d$, where $d$ is the local degree of $\bar{\phi}$ along the loop edge
corresponding to $v'$.  By Theorem~\ref{thm:lifting2}, there is a tame
covering $(X',V' \cup V_0' \cup D') \to (X,V \cup V_0 \cup D)$ lifting
$\tilde{\Sigma}' \to \tilde{\Sigma}$, where $V_0'$ (resp.\ $V_0$)
corresponds to $\bar{V}'_0$ (resp.\ $\bar{V}_0$).  Removing the vertices
in $V_0'$ and $V_0$ gives a tame covering $(X', V' \cup D') \to (X, V \cup
D)$ lifting $\bar{\phi}$.  \qed

\begin{eg} \label{eg:2.isogeny}
  In this example we suppose that $\chr(k) \neq 2$.  Let 
  $E$ be a Tate curve over $K$, let $\Sigma$ be its (set-theoretic)
  skeleton, and let $\tau: E^\an\to\Sigma$ be the canonical retraction.
  Let $U:\G_m^\an/q^\Z\isom E^\an$ be the Tate uniformization of $E^\an$.  The
  $2$-torsion subgroup of $E$ is 
  $U(\{\pm 1,\pm\sqrt q\})$; choose a square root of $q$, let
  $y = U(1)$ and $z = U(\sqrt q)$, and let
  $V = \{y,z\}\subset\Sigma$.  This is a semistable vertex
  set of $E$ and $\Sigma = \Sigma(E,V)$ is the circle with circumference 
  $\val(q)$; the points $y$ and $z$ are antipodal on $\Sigma$.  Let
  $e_1,e_2$ be the edges of $\Sigma$; orient $e_1$ so that $y$ is the
  source vertex and $e_2$ so that $z$ is the source vertex.  The residue
  curves $C_y$ and $C_z$ are both isomorphic to $\P^1_k$; fix isomorphisms
  $C_y\cong\P^1_k$ and $C_z\cong\P^1_k$ such that the tangent direction at
  a vertex in the direction of the outgoing (resp.\ incoming) edge
  corresponds to $\infty$ (resp.\ $0$).

  Let $\Sigma'$ be a circle of circumference $\frac 12\val(q)$, let 
  $V' = \{y',z'\}$ be a pair of antipodal points on $\Sigma'$, and let
  $e_1' = \overrightarrow{y'z'}$ and $e_2'=\overrightarrow{z'y'}$ be the two edges of
  $\Sigma'$, with the indicated orientations.  Enrich $\Sigma'$ with the
  structure of a metrized complex of curves by setting 
  $C_{y'} = C_{z'} = \P^1_k$ and letting the tangent direction at a vertex
  in the direction of the outgoing (resp.\ incoming) edge correspond to
  $\infty$ (resp.\ $0$).  Define a morphism
  $\phi:\Sigma'\to\Sigma$ of metric graphs by 
  \[ \phi(y')=y, \quad \phi(z') = z, \quad \phi(e_1') = e_1, \quad
  \phi(e_2') = e_2, \]
  with both degrees equal to $2$.  See
  Figure~\ref{fig:twoisogeny}.  
  
  Topologically, $\Sigma'\to\Sigma$ is a
  homeomorphism.  Let $\phi_{y'}:C_{y'}\to C_y$ and 
  $\phi_{z'}:C_{z'}\to C_z$ both be the morphism
  $t\mapsto t^2:\P^1_k\to\P^1_k$.  This makes $\phi$ into a degree-$2$
  tame covering of metrized complexes of curves.

  \genericfig[ht]{twoisogeny}{This figure illustrates
    Example~\ref{eg:2.isogeny}.  The morphism $\phi$ is a homeomorphism of
    underlying sets but has a degree of $2$ on $e_1'$ and $e_2'$.
    The arrows on $e_1,e_2,e_1',e_2'$ represent the chosen orientations and
    the indicated tangent vectors at $y,z,y',z'$
    represent the tangent direction corresponding to $\infty$ in the
    corresponding $k$-curve.}

  For $x\in V$ the analytic space $Y(x)$ is an open annulus of logarithmic
  modulus $\val(q)$; fix an isomorphism $Y(x)\cong\bS(q)_+$ for each $x$.  For
  $x'\in V'$ the analytic space $Y'(x')$ is an open annulus of logarithmic
  modulus $\frac 12\val(q)$; fix isomorphisms $Y'(x')\cong\bS(\sqrt q)_+$
  such that the morphisms $Y'(x')\to Y(\phi(x'))$ are given by
  $t\mapsto t^2: \bS(\sqrt q)_+\to\bS(q)_+$.
  For $e'\in E(\Sigma')$ we have $d_{\phi}(e') = 2$ by definition. 
  Hence $\cG(\Sigma',E)$ has four elements, which we label
  $(\pm 1,\pm 1)$.  By Theorem~\ref{thm:lifting2}(1)
  there are four corresponding classes of liftings of $\Sigma'$ to a
  tame covering of $(E,V)$ up to isomorphism preserving $\Sigma'$.
  
  Any $\Sigma$-automorphism of $\Sigma'$ is the identity
  on the underlying topological space.  Hence an element of
  $\Aut_\Sigma(\Sigma')$ is a pair of automorphisms
  $(\psi_{y'},\psi_{z'})\in \Aut_{C_y}(C_{y'})\times\Aut_{C_z}(C_{z'})$ 
  such that $\psi_{y'}$, $\psi_{z'}$ fix the points 
  $0,\infty\in\P^1_k$.  Thus $\psi_{y'},\psi_{z'} = \pm 1$, so
  $\Aut_\Sigma(\Sigma') = \{\pm 1\}\times\{\pm 1\}$.  
  For $x' = y',z'$ the automorphism $-1:C_{x'}\isom C_{x'}$ lifts to the
  automorphism $-1:\bS(\sqrt q)_+\isom\bS(\sqrt q)_+$.  Therefore the
  conjugation action of $\Aut_\Sigma(\Sigma')$ on $\cG(\Sigma',E)$ is
  given as follows:
  \[\begin{split}
    (1,1)\cdot(\pm 1,\pm 1) &= (-1,-1)\cdot(\pm 1,\pm 1) = (\pm 1,\pm 1) \\
    (-1,1)\cdot(\pm 1,\pm 1) &= (1,-1)\cdot(\pm 1,\pm 1) = -(\pm 1,\pm 1). 
  \end{split}\]
  By Theorem~\ref{thm:lifting2}(2), there are two isomorphism classes of
  lifts of $\Sigma'$ to a tame covering of $(E,V)$, and each such lift has
  two automorphisms.

  These liftings can be described concretely as follows.  Fix a square
  root of $q$, let $E_\pm$ be the algebraization of the analytic elliptic curve
  $\G_m^\an/(\pm\sqrt q)^\Z$, and let $\psi_\pm:E_\pm\to E$ be the
  morphism $t\mapsto t^2$ on uniformizations. 
  Let $\Sigma_\pm = \psi_\pm\inv(\Sigma)$.  Then
  $\Sigma_\pm$ is isomorphic to $\Sigma'$ as a tame covering of $\Sigma$
  and $E_\pm$ is a lifting of $\Sigma'$ to a tame covering of $(E,D)$.
  The elliptic curves $E_\pm$ are not isomorphic (as $K$-schemes)
  because they have different $q$-invariants, so they represent the two isomorphism classes
  of liftings of $\Sigma'$.  The nontrivial
  automorphism of $E_\pm$ is given by translating by the image of
  $-1\in\G_m^\an(K)$ (this is not a homomorphism).
  In fact, since $\psi_\pm:E_\pm\to E$ is an \'etale Galois cover of
  degree $2$, this is the only nontrivial automorphism of $E_\pm$ as an
  $E$-curve, so the homomorphism 
  \[ \Aut_E(E_\pm) \To \Aut_\Sigma(\Sigma') \cong \{\pm 1\}\times\{\pm 1\} \]
  is injective but not surjective (its image is
  $\{\pm(1,1)\}$).
\end{eg}

\paragraph Consider now the automorphism group $\Aut^0_\Sigma(\Sigma')
\subset \Aut_\Sigma(\Sigma')$ consisting of all degree-$1$ finite harmonic
morphisms $\alpha: \Sigma'\rightarrow \Sigma'$ respecting
$\phi:\Sigma'\rightarrow \Sigma$ and inducing the identity on the metric graph
$\Gamma'$ underlying $\Sigma'$.  The restriction of the conjugacy action
of $\Aut_\Sigma(\Sigma')$ on $\cG(\Sigma',X)$ to the subgroup
$\Aut^0_\Sigma(\Sigma')$ admits a simplified description that we describe
now.  Combining this with arguments similar to those in the proof of
Theorem~\ref{thm:lifting2} yields a classification of the set of liftings
of $\Sigma'$ up to isomorphism as liftings of the metric graph underlying
$\Sigma'$; see~Theorem~\ref{thm:lifting3}.
 
 First, for each $x'\in V(\Sigma')$ with image 
$x = \phi(x')$, let $\Aut^0_{C_x}(C_{x'})$ be the subgroup of
$\Aut_{C_x}(C_{x'})$ that fixes every point of $C_{x'}$ of the form
$\red(e')$ for some edge $e'$ of $\Sigma'$ adjacent to $x'$.
Then
\[ \Aut^0_\Sigma(\Sigma') = 
\prod_{x'\in V_f(\Sigma')} \Aut^0_{C_{\phi(x)}} (C_{x'}) =: \cE^0. \]
 Denote by $\cE^1$ the finite abelian group
\[ \cE^1 = \prod_{e'\in E_f(\Sigma')}
\Aut_{\tau\inv(\phi(e')^\circ)}(d_{e'}(\phi)). \]
 The discussion preceding Theorem~\ref{thm:lifting2} shows that 
the set of gluing data~$\cG(\Sigma',X)$ is canonically a principal homogeneous space under $\cE^1$.

The subgroup $\Aut^0_{C_x}(C_{x'})\subset\Aut_{C_x}(C_{x'})$ corresponds
to the subgroup $\Aut^0_{Y(x)}(Y'(x'))$ of automorphisms in
$\Aut_{Y(x)}(Y'(x'))$ which act
trivially on the skeleton $\Sigma'(x')$. 
Restriction of a $Y(x)$-automorphism of $Y'(x')$ to 
$\tau_{x'}\inv(e'^\circ)$ defines a homomorphism 
\[ \rho_{x',e'}~:~ \Aut^0_{C_x}(C_{x'}) = \Aut^0_{Y(x)}(Y'(x')) \To
\Aut_{\tau\inv(e^\circ)}(d_{e'}(\phi)). \]
 Fix an orientation of each finite edge of $\Sigma'$.
For $x'\in V_f(\Sigma')$ with image $x = \phi(x')$, let
$\rho_{x'}: \Aut^0_{C_x}(C_{x'})\to\cE^1$ be the homomorphism whose
$e'$-coordinate is
\[ (\rho_{x'}(\alpha))_{e'} =
\begin{cases}
  \rho_{x',e'}(\alpha)     &\qquad\text{ if $x'$ is the source vertex of $e'$} \\
  \rho_{x',e'}(\alpha)\inv &\qquad\text{ if $x'$ is the target vertex of $e'$} \\
  1          &\qquad\text{ if $x'$ is not an endpoint of $e'$.}
\end{cases}\]
Taking the product over all $x'\in V_f(\Sigma')$
yields a homomorphism $\rho: \cE^0 \to \cE^1$.  The kernel and cokernel of
$\rho$ are independent of the choice of orientations.
Viewing $\cE^0\to\cE^1$  as a two-term complex $\cE^\bullet$ of groups, its
cohomology groups are
\[ H^0(\cE^\bullet) = \ker(\rho) \sptxt{and}
H^1(\cE^\bullet) = \coker(\rho). \]

\begin{thm}\label{thm:lifting3}
  Let $(X,V\cup D)$ be a triangulated punctured curve with skeleton
  $\Sigma = \Sigma(X,V\cup D)$.  Assume that $\Sigma$ has no loop edges.
  Let $\phi: \Sigma'\to\Sigma$ be a tame covering of metrized compexes of
  curves. 
  \begin{enumerate}
  \item $\cG(\Sigma',X)$ is canonically a principal homogeneous space under $\cE^1$ 
  and the conjugacy action of $\Aut^0_\Sigma(\Sigma')$ on $\cG(\Sigma',X)$
  is given by the action of the subgroup $\rho(\Aut^0_\Sigma(\Sigma')) = \rho(\cE^0) \subseteq \cE^1$ 
  on $\cG(\Sigma',X)$.

  \item The set of liftings of $\Sigma'$ up to isomorphism as liftings of
    the metric graph underlying $\Sigma'$ is a principal homogeneous space
    under $H^1(\cE^\bullet)$, and the group of automorphisms of a given
    lifting as a lifting of the metric graph underlying $\Sigma'$ is
    isomorphic to $H^0(\cE^\bullet)$.
  \end{enumerate}
\end{thm}

\paragraph[Descent to a general ground field] 
\label{par:descent.lifting}
Let $K_0$ be a subfield of $K$, let $X_0$ be a smooth, projective,
geometrically connected $K_0$-curve, and let $D\subset X_0(K_0)$ be a finite
set.  Let $X = X_0\tensor_{K_0} K$, let $V$ be a strongly semistable vertex
set of $(X,D)$, let $\Sigma = \Sigma(X,V\cup D)$, and let
$\phi: \Sigma'\to\Sigma$ be a tame covering of metrized complexes of
curves, as in the statement of Theorem~\ref{thm:lifting2}.
Let $\phi: X'\to X$ be a lifting of $\phi$ to a tame covering of 
$(X,V\cup D)$.  Whereas we take the data of the morphism
$\Sigma'\to\Sigma$ to be geometric, i.e.\ only defined over $K$, the
covering $X'\to X$ is in fact defined over a finite, separable extension
of $K_0$.  This follows from the fact that  
if $U = X\setminus D$ and $U' = X'\setminus\phi\inv(D)$, then
$\phi: U'\to U$ is a tamely ramified cover of $U$ over $X$ relative to
$D$ (see Remarks~\ref{rem:tamely.ram}(1)
and~\ref{rem:after.tame.cover}(1)), along with the following lemma.

\begin{lem}
  Let $K_0$ be any field, let $K$ be a separably closed field
  containing $K_0$, let $X_0$ be a smooth, projective, geometrically
  connected $K_0$-curve, and let $D \subset X_0(K_0)$ be a finite set.
  Let $X = X_0\tensor_{K_0} K$, let $\phi: X'\to X$ be a finite morphism
  with $X'$ smooth and (geometrically) connected,
  and suppose that $\phi$ is branched only over $D$, with all ramification
  degrees prime to the characteristic of $K$.  Then there exists a finite,
  separable extension $K_1$ of $K$ and a morphism 
  $\phi_1:X_1'\to X_0\tensor_{K_0} K_1$ descending $\phi$.
\end{lem}

\pf Let $U_0 = X_0\setminus D$, and let $U = X\setminus D$ and 
$U' = X'\setminus\phi\inv(D)$.
First suppose that $K_0$ is separably closed.
By \cite[Expos\'e~XIII, Corollaire~2.12]{SGA1}, the tamely ramified
\'etale fundamental groups $\pi_1^t(U_0)$ and $\pi_1^t(U)$ are
isomorphic (with respect to some choice of base point).  Since 
$\phi: U'\to U$ is a tamely ramified cover of $U$ over
$X$ relative to $D$, it is classified by a finite-index subgroup of 
$\pi_1^t(U) = \pi_1^t(U_0)$, so there exists a tamely ramified cover 
$\phi_0: U_0'\to U_0$ of $U_0$ over $X_0$ relative to $D$ descending
$\phi$.

Now we drop the hypothesis that $K_0$ is separably closed.  By the
previous paragraph we may assume that $K$ is a separable closure of
$K_0$.  By general principles the projective morphism $X'\to X$ descends
to a subfield of $K$ which is finitely generated (i.e.\ finite) over
$K_0$.\qed

\begin{rem}\label{rem:rational.lifting}
  With the notation in~\parref{par:descent.lifting}, suppose that
  $K_0$ is a complete valued field with value group
  $\Lambda_0 = \val(K_0^\times)$ and \textbf{algebraically closed} residue field $k$,
  that $\Sigma$ is ``rational over $K_0$'' in that it comes from a (split)
  semistable formal $R_0$-model of $X_0$ in the sense
  of Section~\ref{sec:simultaneous.ss.reduction}, and that
  $\Sigma'$ has edge lengths contained in $\Lambda_0$.  With some extra
  work it is possible to carry out the 
  gluing arguments of Theorem~\ref{thm:lifting} directly over the field
  $K_0$ (in this context Lemma~\ref{lem:extract.roots}(1) still holds),
  which   shows that the cover $X'\to X$ is in fact defined over 
  $K_0$.  In the case of a discrete valuation this also follows
  from~\cite{wewers:tame_covers}, or from~\cite{Sai97} if $D=\emptyset$.
\end{rem}

\paragraph[Liftings of tame harmonic morphisms] 
\label{par:tame.morphism.lifting}
Theorem~\ref{thm:lifting} implies the existing of liftings for (finite)
tame harmonic morphisms of metrized complexes which are not necessarily
tame coverings, the difference being generic \'etaleness.
See~Definition~\ref{def:tame.covering.metrized.complexes} for both
definitions.

\begin{prop} \label{prop:lifting} 
  Let $(X,V \cup D)$ be a triangulated
  punctured $K$-curve, let $\Sigma = \Sigma(X, V \cup D)$, and let $\phi:
  \Sigma'\to\Sigma$ be a tame harmonic morphism of metrized complexes of
  curves.  Then there exists a triangulated punctured $K$-curve 
  $(X',V' \cup D')$ and a finite morphism 
  $\psi : (X',V' \cup D') \to (X, V \cup D)$ lifting $\phi$.
\end{prop}

\pf For any finite vertex $x'$ of $\Sigma'$ at which $\phi$, seen as a
morphism of augmented metric graphs,  
is ramified, let $q_1,\ldots,q_r \in C_v$ be all the branch points of $\phi_{x'}$ which do not correspond 
to any edge of $\Sigma$,
and let $p_{ij}$ denote the preimages of $q_i$ under $\phi_{x'}$.
We modify $\Sigma'$ and $\Sigma$ by attaching infinite edges $e_i$ to $\Sigma$ at $x=\phi(x')$ 
for each $q_i$, infinite edges $e'_{ij}$ to $\Sigma'$ at $x'$ for each $p_{ij}$, and defining
$\red_{x}(e_i) = q_i$ and 
$\red_{x'}(e'_{ij})=p_{ij}$.  Since each $\phi_{x'}$ is tamely ramified,
the harmonic morphism $\phi$ naturally extends to a tame covering 
$\widetilde{\phi} : \widetilde{\Sigma'} \to \widetilde{\Sigma}$
between the resulting
modifications. Enlarging $D$ to $\widetilde D$ by choosing points in $X(K) \setminus D$ with reduction $q_i$, 
we can assume that $\widetilde \Sigma = \Sigma(X, V\cup \widetilde D)$. The result now follows from 
Theorem~\ref{thm:lifting} by first lifting $\widetilde \phi$ to a tame covering
$\psi : (X',V' \cup \widetilde D') \to (X, V \cup \widetilde D)$ and then taking the restriction of $\psi$ to 
$(X', V'\cup D') \to (X,V\cup D)$, where $D' = \psi\inv(D)$.
\qed

\begin{rem} \label{rem:many.lifts}
  As mentioned above, in Proposition~\ref{prop:lifting} we do not require
  that $\phi:\Sigma'\to\Sigma$ be generically \'etale.  This corresponds
  to not requiring that $D\subset X(K)$ contain the branch locus of the
  lift $\psi:X'\to X$: indeed, by Lemma~\ref{lem:tame.implies.tame}, if
  $D$ contains the branch locus then $\phi$ is a tame covering.
  In the situation of Proposition~\ref{prop:lifting}
  the set of liftings of $\Sigma'$ to a tame covering of $(X,V\cup D)$ can
  be infinite.  For example, if $\Gamma' = \{ p' \}$ and $\Gamma = \{ p
  \}$ are both points and the morphism $\phi_{p'} : C_{p'} \cong \PP^1 \to
  C_p \cong \PP^1$ is $z \mapsto z^2$, with $(X,V)$ a minimal
  triangulation of $\PP^1$ and $D=\emptyset$, then there are infinitely
  many such lifts, corresponding to the different ways of lifting the
  critical points and critical values of $\phi_{p'}$ from $k$ to $K$.
\end{rem}

\section{Application: component groups of N\'eron models}\label{sec:applications1}

\paragraph
\label{par:ribet}
In contrast to the rest of the paper, we assume in this section that $R$ is a complete {\em discrete valuation ring} with fraction field $K$ and 
algebraically closed residue field $k$.  In this case we take the value group
$\Lambda$ to be $\Z$.

There is a natural notion of harmonic $1$-forms on a metric graph 
$\Gamma$ of genus $g$ (see \cite{MikhalkinZharkov}).
The space $\Omega^1(\Gamma)$ is a $g$-dimensional real vector space which 
can be canonically identified with $H^1(\Gamma,\RR)$, but we write elements of $\Omega^1(\Gamma)$ as $\omega = \sum_{e} \omega_e \, de$
as in \cite{BakerFaber}, where the sum is over all edges with respect to a fixed vertex
set for $\Gamma$.
There is a canonical lattice $\Omega^1_{\ZZ}(\Gamma)$ of {\em integer
  harmonic $1$-forms} inside $\Omega^1(\Gamma)$; these are the
harmonic $1$-forms for which every $\omega_e$ is an integer.
A harmonic morphism $\varphi : \Gamma' \to \Gamma$ induces a natural pullback map on harmonic $1$-forms via the formula 
\[
\varphi^* (\sum_{e} \omega_{e} \, de) = \sum_{e'} 
d_{e'}(\phi) \cdot \omega_{\varphi(e')} \, de'.
\]

Recall that a {\em $\Z$-metric graph} (i.e.\ a $\Lambda$-metric graph for
$\Lambda=\Z$) with no infinite vertices is a (compact and finite) metric graph whose edge lengths are all positive integers
(or equivalently, having a vertex set with respect to which all edge lengths are $1$).
If $X/K$ is a smooth, proper, geometrically connected analytic curve and $\fX$ is a semistable $R$-model for $X$, 
the skeleton $\Gamma_{\fX}$ of $\fX$ is naturally a $\ZZ$-metric graph.
Moreover, as we have seen earlier in this paper, a finite
morphism of semistable models induces in a natural way a harmonic morphism 
of $\ZZ$-metric graphs.

\medskip

Let $\Gamma$ be a $\ZZ$-metric graph with no infinite vertices.  We define the {\em regularized Jacobian} $\Jac_{\rm reg}(\Gamma)$ of $\Gamma$ to be the group $\Jac(G)$, 
where $G$ is the ``regular model'' for $\Gamma$ (the unique graph induced
by a vertex set with all edge lengths equal to $1$).
The group $\Jac(G)$ can be described explicitly as
\[
\Jac(G) = \Omega^1(G)^\# / H_1(G,\ZZ),
\]
where $\Omega^1(G)^\#$ denotes the linear functionals $\Omega^1(G,\RR) \to \RR$ of
the form $\int_\alpha$ with $\alpha \in C_1(G,\ZZ)$ (see \cite{BakerFaber}).
There is a canonical isomorphism between $\Jac(G)$ and $\Pic^0(G)$, the group of divisors of degree $0$ on $G$ modulo the principal divisors 
(see \cite{BdlHN})
as well as a canonical isomorphism 
\[
\Jac(G) \cong H^1(G,\ZZ)^\# / \Omega^1_\ZZ(G),
\]
where 
\[
H^1(G,\ZZ)^\# = \bigg\{ \omega \in \Omega^1(G,\RR) \; : \; \int_\gamma \omega \in \ZZ \; \forall \gamma \in H_1(G,\ZZ) \bigg\}.
\]

We recall (in our own terminology) the following result of Raynaud
(cf.~\cite{raynaud:picard_specialization} 
and~\cite[Appendix~A]{baker:specialization}): 

\begin{thm}[Raynaud]
\label{thm:Raynaud}
If $X/K$ is a semistable curve, then the component group of the N{\'e}ron model of $\Jac(X)$ over $R$ is canonically isomorphic to 
$\Jac_{\rm reg}(\Gamma_{\fX})$ (for any semistable model $\fX$ of $X$).
\end{thm}

\medskip

A harmonic morphism $\varphi : \Gamma' \to \Gamma$ of $\ZZ$-metric graphs induces in a functorial way homomorphisms $\varphi_* : \Jac_{\rm reg}(\Gamma') \to \Jac_{\rm reg}(\Gamma)$,
defined by
\[
\varphi_*\bigg(\bigg[\int_{\alpha'}\bigg]\bigg) = \bigg[\int_{\varphi_*(\alpha')}\bigg],
\]
where
\[
\varphi_*\big(\sum \alpha_{e'} e'\big) = \sum \alpha_{e'} \varphi(e'),
\]
and $\varphi^* : \Jac_{\rm reg}(\Gamma) \to \Jac_{\rm reg}(\Gamma')$, defined by
\[
\varphi^*([\omega])=[\varphi^* \omega],
\]
where 
\[
\varphi^*\big(\sum \omega_{e} de\big) = \sum_{e'} d_{e'}(\varphi) \omega_{\varphi(e')} de'.
\]

We have the following elementary result, whose proof we omit:

\begin{lem}
\label{lemma:lowerstarisomorphism}
Under the canonical isomorphism between $\Jac(G)$ and $\Pic^0(G)$, the homomorphism $\varphi_*: \Jac(G') \to \Jac(G)$ corresponds to the map 
$[D'] \mapsto [\varphi_*(D')]$ from $\Pic^0(G')$ to $\Pic^0(G)$, where
$\varphi_*(\sum a_{v'} (v')) = \sum a_{v'} \varphi(v')$.
\end{lem}

\medskip

The following is a ``relative'' version of Raynaud's theorem; it implies that
the covariant functor which takes a semistable $R$-model $\fX$ for a $K$-curve $X$ 
to the component group of the N{\'e}ron model of the Jacobian of $X$
factors as the ``reduction graph functor'' $\fX \mapsto \Gamma_\fX$ 
followed by the ``regularized Jacobian functor''
$\Gamma \mapsto \Jac_{\rm reg}(\Gamma)$.  
It is a straightforward consequence of the analytic description of Raynaud's theorem given in 
\cite[Appendix A]{baker:specialization} (see also \cite{baker_rabinoff:analytic_raynaud}):

\begin{thm}
\label{thm:RelativeRaynaud}
If $f_K : X' \to X$ is a finite morphism of curves over $K$, the
induced maps $f_* : \Phi_{J(X')} \to \Phi_{J(X)}$ and
$f^* : \Phi_{J(X)} \to \Phi_{J(X')}$ on component groups coincide with the induced maps $\varphi_* : \Jac_{\rm reg}(\Gamma_{\fX'}) \to \Jac_{\rm reg}(\Gamma_{\fX})$ and
$\varphi^* : \Jac_{\rm reg}(\Gamma_{\fX}) \to \Jac_{\rm reg}(\Gamma_{\fX'})$
for any morphism $f : \fX' \to \fX$ of semistable models 
extending $f_K$.
(Here $\varphi$ denotes the harmonic morphism of skeleta induced by $f$; see
Remark~\ref{rem:induced.harmonic.morphism.2}.)
\end{thm}

\medskip

It follows easily from Lemma~\ref{lemma:lowerstarisomorphism} and Theorem~\ref{thm:RelativeRaynaud} that 
if $f_K : X' \to X$ is {\em regularizable}, i.e., if $f_K$ extends to a morphism of regular semistable models, then the induced map
$f_* : \Phi_{X'} \to \Phi_{X}$ on component groups is surjective.
Thus whenever $f_*$ is not surjective, it follows that $f_K$ does not extend to a morphism of regular semistable models. 
 One can obtain a number of concrete examples of this situation from modular curves (e.g.\ the map $X_0(33) \to E$ over $\QQ_3^{\rm unr}$, where $E$ is the optimal elliptic curve of level $33$.)

\begin{rem}
One can show that if $\varphi : \Gamma' \to \Gamma$ is a harmonic morphism of $\ZZ$-metric graphs, 
then $\varphi_* : \Phi_{X'} \to \Phi_{X}$ is surjective iff $\varphi^* : \Phi_{X} \to \Phi_{X'}$ is injective.
Indeed, it is not hard to check that the maps $\varphi_*$ and $\varphi^*$ are adjoint with respect to the {\em combinatorial monodromy pairing},
the non-degenerate symmetric bilinear form $\langle \; , \; \rangle : \Jac(G) \times \Jac(G) \to \QQ / \ZZ$ defined by
$\langle [\omega] , [\int_\alpha] \rangle = [\int_\alpha \omega]$,
where $\omega \in H_1(G,\ZZ)^\#$ and $\int_\alpha \in \Omega^1(G)^\#$.
Thus the groups ${\rm ker}(\varphi^*)$ and ${\rm coker}(\varphi_*)$
are canonically dual.  This is a combinatorial analogue of results proved by Grothendieck in SGA7 on the (usual) monodromy pairing.
\end{rem}

\medskip

As an application of Theorem~\ref{thm:RelativeRaynaud} and our results on lifting harmonic morphisms, one can
construct many examples of harmonic morphisms of $\ZZ$-metric graphs for which $\varphi_*$ is not surjective.  
For example, consider the following question posed by Ken Ribet in a 2007 email correspondence with the second author (Baker):

\begin{itemize}
\item[]
Suppose $f : X' \to X$ is a finite morphism of semistable curves over a complete discretely valued field $K$ with $g(X)\geq 2$. 
 Assume that the special fiber of the minimal regular model of $X'$ consists of two projective lines intersecting transversely.  Is the induced map $f_* : \Phi_{X'} \to \Phi_X$ on component groups of N{\'e}ron models necessarily surjective?
\end{itemize}

\medskip

We now show that the answer to Ribet's question is {\bf no}.  

\begin{eg}
\label{ex:RibetExample}
Consider the ``banana graph'' $B(\ell_1,\ldots,\ell_{g+1})$ consisting of two vertices
and $g+1$ edges of length $\ell_i$.  This is the reduction graph of a
semistable curve whose reduction has two $\PP^1$'s crossing transversely
at singular points of thickness $\ell_1,\ldots,\ell_{g+1}$.  
If we set $G'=B(1,1,1,1)$ and $G = B(1,2,2)$ and let $\Gamma',\Gamma$ be the geometric realizations of
$G'$ and $G$, respectively, then there is a
degree $2$ harmonic morphism of $\ZZ$-metric graphs $\varphi : \Gamma' \to \Gamma$ taking $e_1'$ and $e_2'$ to $e_1$, $e_3'$ to $e_2$, and $e_4'$ to
$e_3$.  The homomorphism $\varphi_*$ is non-surjective since $|\Jac_{\rm reg}(\Gamma)|=4$ and
$|\Jac_{\rm reg}(\Gamma')|=8$.  The map $\varphi$ can be enriched to a harmonic morphism $\tilde{\varphi}$ of
metrized complexes of curves by attaching a $\P^1$ to each vertex and
letting each morphism $\P^1\to\P^1$ be $z\mapsto z^2$, with $0$, $\pm 1$,
and $\infty$ being the marked points upstairs.  
By Remark~\ref{rem:rational.lifting}, 
the morphism $\tilde{\varphi}$ lifts to a morphism $\psi : X' \to X$ of curves over $K$.
Since all edges of $G'$ have length $1$, the minimal proper regular model of $X'$ has a special fiber consisting of
two $\PP^1$'s crossing transversely.  By Theorem~\ref{thm:RelativeRaynaud}, the induced covariant map on
component groups of N{\'e}ron models is not surjective.
\end{eg}

\bibliographystyle{thesis}
\bibliography{harm_etale}
\bigskip~\bigskip

\immediate\closeout\exportaux

\end{document}

%% file: Figures/ExNonHarm.pdf_t
\begin{picture}(0,0)%
\includegraphics{Figures/ExNonHarm.pdf}%
\end{picture}%
\setlength{\unitlength}{4144sp}%
\begingroup\makeatletter\ifx\SetFigFont\undefined%
\gdef\SetFigFont#1#2#3#4#5{%
  \reset@font\fontsize{#1}{#2pt}%
  \fontfamily{#3}\fontseries{#4}\fontshape{#5}%
  \selectfont}%
\fi\endgroup%
\begin{picture}(4448,2014)(4388,-6023)
\put(6976,-4336){\makebox(0,0)[lb]{\smash{{\SetFigFont{25}{30.0}{\familydefault}{\mddefault}{\updefault}{\color[rgb]{0,0,0}$d$}%
}}}}
\put(8821,-5236){\makebox(0,0)[lb]{\smash{{\SetFigFont{25}{30.0}{\familydefault}{\mddefault}{\updefault}{\color[rgb]{0,0,0}$\phi$}%
}}}}
\end{picture}%

%% file: Figures/ExHarm.pdf_t
\begin{picture}(0,0)%
\includegraphics{Figures/ExHarm.pdf}%
\end{picture}%
\setlength{\unitlength}{4144sp}%
\begingroup\makeatletter\ifx\SetFigFont\undefined%
\gdef\SetFigFont#1#2#3#4#5{%
  \reset@font\fontsize{#1}{#2pt}%
  \fontfamily{#3}\fontseries{#4}\fontshape{#5}%
  \selectfont}%
\fi\endgroup%
\begin{picture}(4628,2014)(4388,-6023)
\put(6301,-4336){\makebox(0,0)[lb]{\smash{{\SetFigFont{25}{30.0}{\familydefault}{\mddefault}{\updefault}{\color[rgb]{0,0,0}$d$}%
}}}}
\put(9001,-5326){\makebox(0,0)[lb]{\smash{{\SetFigFont{25}{30.0}{\familydefault}{\mddefault}{\updefault}{\color[rgb]{0,0,0}$\phi$}%
}}}}
\end{picture}%

%% file: Figures/EffectiveMorphi2.pdf_t
\begin{picture}(0,0)%
\includegraphics{Figures/EffectiveMorphi2.pdf}%
\end{picture}%
\setlength{\unitlength}{4144sp}%
\begingroup\makeatletter\ifx\SetFigFont\undefined%
\gdef\SetFigFont#1#2#3#4#5{%
  \reset@font\fontsize{#1}{#2pt}%
  \fontfamily{#3}\fontseries{#4}\fontshape{#5}%
  \selectfont}%
\fi\endgroup%
\begin{picture}(4943,3599)(4388,-7598)
\put(9316,-5146){\makebox(0,0)[lb]{\smash{{\SetFigFont{25}{30.0}{\familydefault}{\mddefault}{\updefault}{\color[rgb]{0,0,0}$\phi$}%
}}}}
\end{picture}%

%% file: Figures/NonEffectiveMorphis.pdf_t
\begin{picture}(0,0)%
\includegraphics{Figures/NonEffectiveMorphis.pdf}%
\end{picture}%
\setlength{\unitlength}{4144sp}%
\begingroup\makeatletter\ifx\SetFigFont\undefined%
\gdef\SetFigFont#1#2#3#4#5{%
  \reset@font\fontsize{#1}{#2pt}%
  \fontfamily{#3}\fontseries{#4}\fontshape{#5}%
  \selectfont}%
\fi\endgroup%
\begin{picture}(4673,3145)(4388,-6023)
\put(5176,-4336){\makebox(0,0)[lb]{\smash{{\SetFigFont{25}{30.0}{\familydefault}{\mddefault}{\updefault}{\color[rgb]{0,0,0}$d$}%
}}}}
\put(6976,-4336){\makebox(0,0)[lb]{\smash{{\SetFigFont{25}{30.0}{\familydefault}{\mddefault}{\updefault}{\color[rgb]{0,0,0}$d$}%
}}}}
\put(6391,-5011){\makebox(0,0)[lb]{\smash{{\SetFigFont{25}{30.0}{\familydefault}{\mddefault}{\updefault}{\color[rgb]{0,0,0}$p'$}%
}}}}
\put(9046,-5146){\makebox(0,0)[lb]{\smash{{\SetFigFont{25}{30.0}{\familydefault}{\mddefault}{\updefault}{\color[rgb]{0,0,0}$\phi$}%
}}}}
\end{picture}%